\documentclass[preprint,12pt]{elsarticle}

\usepackage{amsmath,amssymb}
\usepackage{graphicx}
\usepackage{subcaption}
\usepackage{booktabs}
\usepackage{hyperref}
\usepackage{tikz}
\usetikzlibrary{arrows.meta,positioning,shapes.geometric,fit}

\journal{Journal of Computational Physics}

\begin{document}

\title{A Massively Parallel Hybridizable Discontinuous Galerkin Solver for Direct Numerical Simulation of Compressible Flows on GPUs}

\author[mit]{Andrew Welter}
\author[mitcsail]{Thea Collin}
\author[mit]{Ngoc Cuong Nguyen}
\author[mit]{Jaime Peraire}
\address[mit]{Department of Aeronautics and Astronautics, Massachusetts Institute of Technology, \ 77 Massachusetts Avenue, Cambridge, MA 02139, USA}
\address[mitcsail]{Department of Electrical Engineering and Computer Science, Massachusetts Institute of Technology, \ 77 Massachusetts Avenue, Cambridge, MA 02139, USA}

\begin{abstract}
Direct numerical simulation (DNS) of compressible transitional and turbulent flows requires numerical methods that combine high-order accuracy, robustness, and computational efficiency to resolve a broad range of spatial and temporal scales.  This paper presents a massively parallel  hybridizable discontinuous Galerkin (HDG) solver for DNS of the compressible Navier–Stokes equations on GPU-accelerated high-performance computing systems. The proposed solver combines high-order HDG  discretization with robust shock capturing, diagonally implicit Runge–Kutta (DIRK) time integration, and an efficient Newton–GMRES solution strategy accelerated by additive Schwarz preconditioning and reduced-basis approximation. A distributed implementation of these methods based on GPU-aware MPI, Kokkos, and CUDA/HIP libraries enables scalable execution on heterogeneous computing platforms. The solver is demonstrated on three canonical benchmark problems covering a wide range of Mach-number flow regimes: subsonic transitional flow over the Eppler 387 airfoil, the supersonic Taylor–Green vortex, and hypersonic boundary-layer transition. Numerical results are compared with available experimental measurements and published DNS data, showing good agreement across distinct flow regimes. The results demonstrate the ability of the proposed solver to resolve laminar–turbulent transition, strong compressibility effects, shock-associated flow structures, and fully three-dimensional turbulent dynamics. 
\end{abstract}


\maketitle

\section{Introduction}

\subsection{Background and Motivation}

Direct numerical simulation (DNS) remains the most systematic approach for resolving transitional and turbulent compressible flows without relying on turbulence closures, but the method is demanding because the flow dynamics span a wide range of coupled spatial and temporal scales. Wall-bounded transition involves the growth, interaction, and breakdown of instability waves inside thin viscous layers; fully developed turbulence transfers energy through a broad cascade of dynamically important scales; and compressibility introduces acoustic waves, thermodynamic coupling, and Mach-number-dependent flow structures that must be represented without excessive numerical dissipation. These requirements are especially important when the flow also contains curved boundaries, laminar separation and reattachment, shock waves, shock--turbulence interactions, or strong wall-normal anisotropy in the near-wall region. In such regimes, numerical methods must deliver high resolution of smooth structures, localized nonlinear stabilization near discontinuities, geometric flexibility, and a solution procedure that remains practical on modern computing systems.

High-order discontinuous Galerkin (DG) methods are attractive in this setting because they combine element-local polynomial approximation with compact interelement coupling, local conservation, and natural compatibility with unstructured meshes \cite{CockburnShuReview,Bassi1997,NguyenPeraireCockburnHDGNS11,Vila-Perez2021}. For eddy-resolving and DNS calculations, their principal appeal lies in the possibility of reducing numerical dispersion and dissipation while retaining high-order accuracy on complex geometries \cite{Bassi2015,Gassner2013,Moura2015,Moura2017}. The element-local structure of DG discretizations is also favorable for parallel implementation and for constructing implicit time integrators whose communication pattern remains compact even at high polynomial order. These properties have made high-order DG methods a compelling alternative to structured finite-difference, compact finite-difference, spectral, and spectral-element approaches for DNS calculations.

For wall-resolved compressible flows, however, spatial approximation order alone is not sufficient. Thin boundary layers, separation bubbles, and shock-associated gradients impose severe mesh requirements, and isotropic refinement can render DNS prohibitively expensive even before the cost of implicit nonlinear solves is considered. Curved high-order elements are needed to represent aerodynamic and hypersonic geometries without degrading formal accuracy, while anisotropic mesh distributions are essential to align resolution with the strong directional stiffness of boundary-layer flows \cite{Hoskin2024,Nguyen2020gpu,Moro2017a,nguyen2023optimal,Bai2022a}. In practice, the efficiency of an unstructured high-order DNS strategy depends on whether the mesh can concentrate degrees of freedom where the solution varies rapidly, such as near walls and shocks, while avoiding unnecessary refinement in smooth regions of the domain.

Compressible transitional and turbulent flows with shocks impose an additional constraint: the discretization must preserve high-order accuracy in smooth regions while introducing dissipation only where it is needed to control nonphysical oscillations. This balance is difficult because excessive stabilization can corrupt instability growth, acoustic propagation, and turbulent spectra, whereas insufficient stabilization can destabilize the simulation in the presence of shocklets, compression waves, or stronger shock-associated structures. Shock-capturing strategies based on localized artificial viscosity and smoothness sensing provide one practical route for high-order DG discretizations \cite{persson06:_shock_capturing}. For DNS of compressible flows, the central issue is not merely to capture a discontinuity, but to do so in a way that preserves the surrounding resolved flow physics and remains compatible with implicit high-order discretizations on unstructured curved meshes.

Time integration presents a second major obstacle. The stiffness induced by viscous terms, acoustic propagation, and high polynomial order often makes explicit time stepping excessively restrictive for DNS on fine meshes. Implicit schemes alleviate the stability constraint, but they replace it with a sequence of large nonlinear algebraic systems whose cost can dominate the simulation. In a diagonally implicit Runge--Kutta (DIRK) discretization, each time step requires the solution of several closely related stage systems; each nonlinear iteration in turn generates a linearized problem whose trace-space operator changes gradually across Newton iterations and across time steps. As a result, the overall efficiency of implicit high-order DNS depends not only on the spatial discretization, but also on the nonlinear solver, the Krylov method, the preconditioner, and the quality of the initial guesses used for repeated linear solves.

These considerations motivate the present work. The paper develops a massively parallel hybridizable discontinuous Galerkin (HDG) methodology for DNS of compressible flows on distributed GPU systems. The aim is not to present GPU execution as an isolated software contribution, but to show how high-order HDG discretization, localized shock capturing, anisotropic curved meshes, DIRK time integration, Newton--GMRES solution, additive Schwarz preconditioning, and reduced-basis approximation can be organized into a single computational framework for transitional and turbulent compressible flows spanning subsonic, supersonic, and hypersonic regimes.

\subsection{Related Work}

High-order DG methods are now well established for compressible-flow simulation and have been studied extensively because their compact stencil, elementwise conservation, and high-order approximation are well suited to wave propagation and multiscale flow physics \cite{CockburnShuReview,Bassi1997,Bassi2015,NguyenPeraireCockburnHDGNS11,Vila-Perez2021}. For under-resolved turbulence and DNS, several studies have emphasized that the low-dissipation and low-dispersion characteristics of high-order discretizations can materially influence the resolved dynamics, particularly when the flow contains a broad spectral range and limited numerical damping is desired \cite{Gassner2013,Moura2015,Moura2017}. At the same time, the use of unstructured elements and curved geometry descriptions extends these benefits to configurations that are difficult to treat efficiently with structured-grid methods. These advantages explain the sustained interest in DG methods for transitional and turbulent compressible flows, but they do not by themselves remove the high cost associated with the globally coupled degrees of freedom in standard implicit DG formulations.

HDG methods address that difficulty through hybridization and static condensation. By expressing the element-interior unknowns in terms of traces on the mesh skeleton, HDG reduces the globally coupled system to interface variables while preserving the local element structure that makes high-order DG methods attractive \cite{CockburnGopalakrishnan09HDGStokes,Moro2011a,NguyenPeraireCockburnHDGNS11,Vila-Perez2021}. This reduction is important for implicit compressible-flow solvers because the condensed trace formulation changes the balance between local dense algebra and global sparse coupling. In smooth viscous-dominated regimes, HDG has also been shown to retain desirable high-order accuracy properties \cite{NguyenPeraireCockburnHDGNS11}. Related hybridized formulations have also been developed for wave propagation, compressible magnetohydrodynamics, embedded discontinuous Galerkin discretizations for computational fluid dynamics, and multiscale compressible-flow discretizations \cite{Fernandez2018a,Ciuca2020,Nguyen2015c,N2013}. The method itself is therefore not new; what remains challenging is to exploit the trace formulation effectively in the fully discrete, nonlinear, and massively parallel setting required for large-scale DNS of compressible flows.

The treatment of shocks and sharp compressive structures remains one of the central issues in high-order DG and HDG simulation of compressible flows. A wide range of approaches has been developed in the broader literature, including limiting, entropy-based stabilization, subcell strategies, artificial-viscosity formulations, and interface-aligned or shock-tracking discretizations \cite{Zahr2020,Kercher2021}.  Artificial-viscosity methods have been developed to treat shock waves \cite{persson06:_shock_capturing,Nguyen2011d,Moro2011a,Moro2016,Fernandez2018}. More recent developments in this direction include shock-capturing strategies for hypersonic non-equilibrium flow, adaptive viscosity regularization, and the combination of shock regularization with mesh adaptivity and high-order hypersonic DG simulation \cite{VanHeyningen2023,Nguyen2023c,nguyen2023optimal,Hoskin2024}. For transitional and turbulent compressible flows, this localization is especially important because shocks may coexist with instability waves, near-wall streaks, separation bubbles, and broadband turbulent fluctuations. The mesh design is closely connected to this issue: anisotropic curved elements are needed not only to represent boundaries accurately, but also to place resolution economically in thin wall layers and other strongly directional structures \cite{Hoskin2024,Nguyen2020gpu,Moro2017a,nguyen2023optimal,Bai2022a}.

Implicit solution strategies for compressible high-order discretizations have likewise received sustained attention. Newton--Krylov methods and GMRES-based linear solvers are a common choice for nonlinear compressible-flow systems because they separate the treatment of nonlinear residual reduction from the solution of large nonsymmetric linearized problems \cite{sasc86,saad93,persson:GMRESDG}. For high-order DG discretizations of the compressible Navier--Stokes equations, previous work has investigated Newton--GMRES solvers together with multigrid and other preconditioning strategies \cite{persson:GMRESDG,Shahbazi2009,Diosady2009,Franciolini2020}. On the HDG side, implicit time integration and hybridized formulations for time-dependent compressible flows have been explored, including DIRK-based strategies and applications to transitional turbulent flows \cite{Jaust2014a,Fernandez2017a}. These studies demonstrate the potential of implicit high-order discretizations, but they also make clear that preconditioning remains a major bottleneck as polynomial order, Reynolds number, and problem size increase.

Domain-decomposition preconditioners provide one of the classical pathways to scalable implicit solvers on distributed-memory systems. Additive Schwarz and restricted additive Schwarz methods have been studied extensively because they offer a flexible compromise between locality, communication cost, and convergence behavior \cite{Fischer1997,Cai1999,Efstathiou2003}. Their appeal in the HDG setting is that the condensed trace system already isolates the globally coupled degrees of freedom on faces or traces, making subdomain solves and overlap strategies natural ingredients of the linear solver. Even so, the repeated linear systems arising from DIRK stages and successive Newton iterations remain a difficult aspect of implicit DNS. Because neighboring solves are closely related, projection and recycling ideas for multiple linear systems provide useful guidance for exploiting information from previous iterations \cite{Chan1999,Frommer1998}. In the present work, that perspective motivates the use of a reduced-basis approximation to construct effective initial guesses for GMRES within a repeated-solve setting.

Progress in hardware has added another layer of complexity. GPU-accelerated high-order CFD solvers have shown that polynomial-based methods can map well to throughput-oriented architectures because much of the work is carried out in local dense tensor or element algebra \cite{Witherden2014,Vermeire2017,Kolev2021}. Portable programming models such as Kokkos provide a path toward maintaining a single implementation across evolving manycore architectures \cite{Carter2014,Trott2022}, and software systems such as Exasim reflect broader efforts to generate or organize high-order DG codes for heterogeneous machines \cite{VilaPerez2022}. These efforts involve a wide range of libraries, techniques, applications, and discretizations~\cite{africa2024deal, kolev2021efficient,andrej2024high,abdelfattah2021gpu,mills2025petsc,mayr2026trilinos, tomov2021libceed}. For HDG in particular, the condensed formulation exposes substantial local dense linear algebra, which is well matched to vendor-optimized GPU libraries, but the method also introduces a distributed trace system whose iterative solution requires sparse coupling, halo exchange, Krylov reductions, and preconditioner application. Early GPU-oriented work for HDG also examined sparse matrix--vector products for the condensed trace system \cite{Roca2011}, and more recent work has examined broader solver and preconditioning strategies for HDG discretizations \cite{Welter2026}. The remaining challenge for DNS is to combine these implementation techniques with an implicit nonlinear solution procedure that remains effective across multiple Mach-number regimes and across multi-GPU distributed-memory executions, rather than only on a single device or within a narrowly defined benchmark class.

The application space considered in this paper further sharpens that need. Transitional flow over the Eppler 387 airfoil is a demanding subsonic validation problem because it combines wall-bounded instability growth, laminar separation, transition, and reattachment on a curved aerodynamic geometry with available low-disturbance experimental data \cite{McGhee1988_2,Fernandez2017a}. The Taylor--Green vortex is a canonical setting for assessing how a high-order compressible discretization represents transition to turbulence and the resulting multiscale dynamics; prior high-order studies have used Taylor--Green-type problems to examine accuracy and dissipation properties in eddy-resolving computations \cite{Gassner2013,Fernandez2017}. Hypersonic boundary-layer transition provides a complementary test because instability amplification, nonlinear breakdown, wall heat transfer, and compressibility effects become tightly coupled, and because both experiments and DNS have shown the importance of quiet-tunnel conditions and second-mode-dominated transition mechanisms in such flows \cite{Kendall1975,Schneider2015,Juliano2015a,Marxen2014,Zhang2013,Hader2019}. Taken together, these three applications probe wall-bounded subsonic transition with separation, compressible transition and turbulence in a canonical unsteady flow, and high-Mach-number boundary-layer breakdown with fine near-wall scales and strong thermal effects.

\subsection{Contributions and Organization}

The gap is not the absence of HDG discretizations, implicit time integrators, Krylov solvers, additive Schwarz methods, or GPU programming models in isolation. Rather, the gap lies in the lack of an integrated massively parallel HDG methodology for DNS of compressible flows that couples these ingredients in a form suitable for repeated implicit solves on distributed GPU systems and demonstrates the resulting capability across subsonic, supersonic, and hypersonic flow regimes. Existing studies establish important pieces of that picture, but the combination of high-order HDG discretization, localized shock treatment, anisotropic curved meshes, Newton--GMRES solution, domain-decomposition preconditioning, and repeated-solve acceleration remains insufficiently consolidated for the class of transitional and turbulent compressible flows considered here.

The present paper addresses that gap with a unified computational methodology. Specifically, it develops a massively parallel HDG solver for  DNS of compressible flows on distributed GPU systems; combines high-order HDG spatial discretization, shock capturing, DIRK time integration, and Newton--GMRES nonlinear solution in a single fully discrete workflow; introduces an additive Schwarz preconditioner tailored to the HDG trace system; and employs a reduced-basis approximation to construct effective initial guesses for repeated GMRES solves. The implementation uses distributed-memory parallelism, portable node-level execution through Kokkos, and vendor-optimized dense linear algebra libraries for NVIDIA and AMD GPUs to organize the local and global algebra generated by HDG condensation. The methodology is designed to work with unstructured curved high-order meshes that concentrate resolution in boundary layers and shock regions, thereby reducing the number of elements required for wall-resolved simulations of complex compressible flows.

The numerical studies then validate the methodology against available experiments and published DNS results for three complementary problems: subsonic transition over the Eppler 387 airfoil, the supersonic Taylor--Green vortex, and hypersonic boundary-layer transition. Through these cases, the paper examines whether the solver can resolve laminar--turbulent transition, compressibility effects, shock-associated structures, and fully three-dimensional turbulent dynamics.

The remainder of the paper is organized as follows. Section~2  introduces the governing equations, the high-order HDG discretization, the DIRK temporal discretization, and the nonlinear algebraic formulation. Section~3 describes the solver architecture, the parallel Newton--GMRES strategy, the additive Schwarz preconditioner, the reduced-basis initialization procedure, and the heterogeneous distributed-memory implementation. Section~4 presents the subsonic transition study. Section~5 examines compressible transition and turbulence in the canonical vortex problem. Section~6 addresses the high-Mach-number boundary-layer case. Section~7 summarizes the main findings and directions for future work.
 
\section{Numerical Methods}

\subsection{Governing Equations}

We consider the compressible Navier--Stokes equations on a spatial domain $\Omega \subset \mathbb{R}^{d}$, with $d=3$, over a time interval $(0,T]$. The formulation is written in conservative form as
\begin{equation}
\label{eq:ns-strong}
\frac{\partial \boldsymbol{u}}{\partial t}
+ \nabla \cdot \left( \boldsymbol{F}^{c}(\boldsymbol{u}) - \boldsymbol{F}^{v}(\boldsymbol{u},\nabla \boldsymbol{u}) \right)
= \boldsymbol{0}
\qquad \text{in } \Omega \times (0,T],
\end{equation}
subject to appropriate initial and boundary conditions. In the present work there are no volumetric source terms in the governing equations. The vector of conservative variables is
\begin{equation}
\label{eq:state-vector}
\boldsymbol{u}
=
\begin{bmatrix}
\rho, &
\rho v_1, &
\ldots, &
\rho v_d, &
\rho E
\end{bmatrix}^{T},
\end{equation}
where $\rho$ is the density, $v_i$ are the Cartesian velocity components, and $E$ is the total specific energy. Repeated spatial indices imply summation from $1$ to $d$. The inviscid flux in the $x_j$-direction is
\begin{equation}
\label{eq:convective-flux}
\boldsymbol{f}^{c}_{j}(\boldsymbol{u})
=
\begin{bmatrix}
\rho v_j \\
\rho v_1 v_j + \delta_{1j} p \\
\vdots \\
\rho v_d v_j + \delta_{dj} p \\
\rho H v_j
\end{bmatrix},
\qquad
\boldsymbol{F}^{c}(\boldsymbol{u})
=
\left(
\boldsymbol{f}^{c}_{1},
\ldots,
\boldsymbol{f}^{c}_{d}
\right),
\end{equation}
where $p$ is the pressure, $\delta_{ij}$ is the Kronecker delta, and $H = E + p/\rho$ is the total specific enthalpy. For a calorically perfect gas,
\begin{equation}
\label{eq:thermo-closure}
E = e + \frac{1}{2} v_i v_i,
\qquad
e = c_v T,
\qquad
h = c_p T,
\qquad
p = \rho R T,
\qquad
\gamma = \frac{c_p}{c_v},
\end{equation}
so that
\begin{equation}
\label{eq:ideal-gas}
p = (\gamma - 1)\rho \left(E - \frac{1}{2} v_i v_i\right).
\end{equation}
The variables may be dimensional or nondimensional depending on the reference scales used in each application problem; in either case the constitutive and thermodynamic relations retain the form above.

The viscous flux is written as
\begin{equation}
\label{eq:viscous-flux}
\boldsymbol{f}^{v}_{j}(\boldsymbol{u},\nabla \boldsymbol{u})
=
\begin{bmatrix}
0 \\
\tau_{1j} \\
\vdots \\
\tau_{dj} \\
v_i \tau_{ij} + f_j
\end{bmatrix},
\qquad
\boldsymbol{F}^{v}(\boldsymbol{u},\nabla \boldsymbol{u})
=
\left(
\boldsymbol{f}^{v}_{1},
\ldots,
\boldsymbol{f}^{v}_{d}
\right),
\end{equation}
where the stress tensor and heat-conduction term are
\begin{equation}
\label{eq:stress-heat}
\tau_{ij}
=
\mu
\left(
\frac{\partial v_i}{\partial x_j}
+
\frac{\partial v_j}{\partial x_i}
-
\frac{2}{3}\frac{\partial v_k}{\partial x_k}\delta_{ij}
\right)
+
\beta \frac{\partial v_k}{\partial x_k}\delta_{ij},
\qquad
f_j = \kappa \frac{\partial T}{\partial x_j}.
\end{equation}
Here $\mu$ denotes the shear viscosity, $\beta$ the bulk viscosity, and $\kappa$ the thermal conductivity. In the absence of shock capturing, the physical coefficients satisfy $\beta = 0$ under Stokes' hypothesis, $\mu = \mu_f$, and
\begin{equation}
\label{eq:thermal-conductivity}
\kappa = \kappa_f = \frac{c_p \mu_f}{Pr},
\end{equation}
with $Pr$ the Prandtl number. When temperature-dependent transport is used, $\mu_f$ is evaluated with Sutherland's law,
\begin{equation}
\label{eq:sutherland}
\mu_f(T)
=
\mu_{\mathrm{ref}}
\left( \frac{T}{T_{\mathrm{ref}}} \right)^{3/2}
\frac{T_{\mathrm{ref}} + S}{T + S},
\end{equation}
or its nondimensional equivalent. The Reynolds, Mach, and Prandtl numbers associated with a given test case are specified in the corresponding application section.

The boundary data are imposed weakly through the trace variable and the boundary numerical flux. Periodic boundaries are treated by identifying paired faces and enforcing equal and opposite numerical fluxes across the pair. On inviscid or far-field boundaries, the boundary state enters through the exterior trace in the numerical flux, with all conservative variables prescribed for supersonic inflow and the interior state extrapolated for supersonic outflow. Subsonic inflow and outflow boundaries are imposed by prescribing the incoming characteristic information together with the required thermodynamic data. On solid walls, the no-penetration condition is enforced through the trace, while the viscous boundary data impose either no-slip isothermal conditions, $ \boldsymbol{v} = \boldsymbol{0}$ and $T = T_w$, or no-slip adiabatic conditions, $ \boldsymbol{v} = \boldsymbol{0}$ and $\partial T/\partial n = 0$. Symmetry boundaries enforce zero normal velocity together with zero normal gradients of the tangential velocity and temperature. The case-specific values of the free-stream and wall data are deferred to Sections~4--6.

For the shocked cases considered later, the discretization is supplemented by a scalar artificial-viscosity (AV) field that adds localized Laplacian diffusion to the conservative variables  of the form
\begin{equation}
\label{eq:av-total}
\boldsymbol{f}^{\mathrm{AV}}_{j}
=
\nu_{\mathrm{AV}}
\frac{\partial \boldsymbol{u}}{\partial x_{j}},
\end{equation}
so that the total diffusive operator is obtained by replacing $\boldsymbol{f}^{v}_{j}$ in \eqref{eq:viscous-flux} with $\boldsymbol{f}^{v}_{j} + \boldsymbol{f}^{\mathrm{AV}}_{j}$. The same scalar coefficient therefore diffuses density, momentum, and total energy through the gradients of the conservative variables. 

In the present work, the raw elemental ariticial viscosity is constructed from a bounded compression sensor with a Ducros-type vortical suppression factor \cite{persson06:_shock_capturing,Nguyen2011a,Moro2016}. Let $\nabla \cdot \boldsymbol{v}$ denote the velocity divergence, $\omega = |\nabla \times \boldsymbol{v}|$ the vorticity magnitude, and $c^{\ast}$ the local acoustic scale computed from a thermodynamic state. The raw artificial viscosity is then written as
\begin{equation}
\label{eq:av-coeffs}
\nu_{\mathrm{AV,raw}}
=
c_{\mathrm{AV}}\,
L\!\left(
\sqrt{\frac{h_{m}}{k}}\,
\frac{L(\nabla \cdot \boldsymbol{v})}{c^{\ast}}\,
\frac{L(\nabla \cdot \boldsymbol{v})^{2}}
{L(\nabla \cdot \boldsymbol{v})^{2} + L(\omega)^{2} + \varepsilon}
\right),
\end{equation}
where $h_{m}$ is the prescribed sensor length scale, $k$ is the polynomial degree, $c_{\mathrm{AV}}$ is the user-specified AV coefficient, $\varepsilon$ is a small positive constant, and $L(\cdot)$ denotes the limiting operation used to enforce non-negativity and upper bounds on the divergence, vorticity, sensor value, and regularized thermodynamic variables \cite{Nguyen2023c,nguyen2023optimal,VanHeyningen2023}. In particular, the sensor responds to compressive regions and is attenuated in predominantly vortical flow by the Ducros factor. The raw AV field \eqref{eq:av-coeffs} is passed through a smoothing operation:
\[
\nu_{\mathrm{AV}} = \mathcal{S}_{h}\!\left[\nu_{\mathrm{AV,raw}}\right],
\]
where $\mathcal{S}_{h}$ represents the discrete smoothing operator defined by averaging the raw AV field to obtain the smooth continuous AV field $\nu_{\mathrm{AV}}$ \cite{Moro2016,Fernandez2018,Ciuca2020}.


\subsection{High-Order HDG Discretization}

Let $\mathcal{T}_{h}$ be a conforming tessellation of $\Omega$ by curved high-order elements $K$, and let
\[
\partial \mathcal{T}_{h} = \bigcup_{K \in \mathcal{T}_{h}} \partial K,
\qquad
\mathcal{F}_{h} = \bigcup_{K \in \mathcal{T}_{h}} \partial K
\]
denote the mesh skeleton. The spatial discretization follows the standard HDG construction in which element-interior variables are coupled only through a single-valued trace on the skeleton \cite{CockburnGopalakrishnan09HDGStokes}. This formulation is part of a broader family of hybridized methods that also includes embedded discontinuous Galerkin variants and HDG formulations for wave propagation \cite{Nguyen2015c,Fernandez2018a}. For polynomial degree $k \ge 1$, we define the volume and trace spaces
\begin{align}
\label{eq:spaces}
\mathcal{Q}_{h}^{k}
&=
\left\{
\boldsymbol{q} \in [L^{2}(\mathcal{T}_{h})]^{m \times d}
:
\boldsymbol{q}|_{K} \in [\mathcal{P}^{k}(K)]^{m \times d},
\ \forall K \in \mathcal{T}_{h}
\right\},
\nonumber \\
\mathcal{V}_{h}^{k}
&=
\left\{
\boldsymbol{w} \in [L^{2}(\mathcal{T}_{h})]^{m}
:
\boldsymbol{w}|_{K} \in [\mathcal{P}^{k}(K)]^{m},
\ \forall K \in \mathcal{T}_{h}
\right\},
\nonumber \\
\mathcal{M}_{h}^{k}
&=
\left\{
\boldsymbol{\mu} \in [L^{2}(\mathcal{F}_{h})]^{m}
:
\boldsymbol{\mu}|_{F} \in [\mathcal{P}^{k}(F)]^{m},
\ \forall F \subset \mathcal{F}_{h}
\right\},
\end{align}
where $m=d+2$ for the compressible Navier--Stokes equations. The curved geometry mapping is assumed to be of sufficiently high order to represent the physical boundary and the boundary-layer mesh anisotropy consistently with the polynomial approximation.

To hybridize the viscous terms, we introduce the auxiliary variable
\begin{equation}
\label{eq:mixed-variable}
\boldsymbol{q} = \nabla \boldsymbol{u},
\end{equation}
which rewrites \eqref{eq:ns-strong} as a first-order system in $(\boldsymbol{q},\boldsymbol{u})$. The HDG semidiscretization then seeks
\[
(\boldsymbol{q}_{h},\boldsymbol{u}_{h},\widehat{\boldsymbol{u}}_{h})
\in
\mathcal{Q}_{h}^{k} \times \mathcal{V}_{h}^{k} \times \mathcal{M}_{h}^{k},
\]
where $\widehat{\boldsymbol{u}}_{h}$ is the single-valued trace unknown on the mesh skeleton. Using the standard element and face inner products
\[
(\cdot,\cdot)_{\mathcal{T}_{h}}
=
\sum_{K \in \mathcal{T}_{h}}(\cdot,\cdot)_{K},
\qquad
\langle \cdot,\cdot \rangle_{\partial \mathcal{T}_{h}}
=
\sum_{K \in \mathcal{T}_{h}}\langle \cdot,\cdot \rangle_{\partial K},
\]
the semidiscrete HDG formulation reads: for all $(\boldsymbol{r},\boldsymbol{w},\boldsymbol{\mu}) \in \mathcal{Q}_{h}^{k} \times \mathcal{V}_{h}^{k} \times \mathcal{M}_{h}^{k}$,
\begin{subequations}
\label{eq:hdg-weak}
\begin{align}
\label{eq:hdg-grad}
(\boldsymbol{q}_{h},\boldsymbol{r})_{\mathcal{T}_{h}}
+
(\boldsymbol{u}_{h},\nabla \cdot \boldsymbol{r})_{\mathcal{T}_{h}}
-
\langle \widehat{\boldsymbol{u}}_{h},\boldsymbol{r}\boldsymbol{\cdot}\boldsymbol{n} \rangle_{\partial \mathcal{T}_{h}}
&= 0,
\\
\label{eq:hdg-cons}
\left(\frac{\partial \boldsymbol{u}_{h}}{\partial t},\boldsymbol{w}\right)_{\mathcal{T}_{h}}
-
\left(\boldsymbol{F}^{c}(\boldsymbol{u}_{h})-\boldsymbol{F}^{v}(\boldsymbol{u}_{h},\boldsymbol{q}_{h}),\nabla \boldsymbol{w}\right)_{\mathcal{T}_{h}}
+
\langle \widehat{\boldsymbol{f}}_{h},\boldsymbol{w} \rangle_{\partial \mathcal{T}_{h}}
&= 0,
\\
\label{eq:hdg-trace}
\langle \widehat{\boldsymbol{f}}_{h},\boldsymbol{\mu} \rangle_{\partial \mathcal{T}_{h}\setminus\partial\Omega}
+
\langle \widehat{\boldsymbol{b}}_{h},\boldsymbol{\mu} \rangle_{\partial \Omega}
&= 0.
\end{align}
\end{subequations}
Equation \eqref{eq:hdg-grad} is the local weak statement of the gradient relation, \eqref{eq:hdg-cons} is the elementwise conservation law, and \eqref{eq:hdg-trace} enforces interelement transmission and boundary conditions through a global trace equation. The unknowns $\boldsymbol{q}_{h}$ and $\boldsymbol{u}_{h}$ are element-local once $\widehat{\boldsymbol{u}}_{h}$ is given, whereas $\widehat{\boldsymbol{u}}_{h}$ carries the global coupling. The numerical flux is defined as
\begin{equation}
\label{eq:num-flux-hdg}
\widehat{\boldsymbol{f}}_{h}
=
\left(
\boldsymbol{F}^{c}(\widehat{\boldsymbol{u}}_{h})
-
\boldsymbol{F}^{v}(\widehat{\boldsymbol{u}}_{h},\boldsymbol{q}_{h})
\right)\boldsymbol{\cdot}\boldsymbol{n}
+
\boldsymbol{\tau}_{h}\left(\boldsymbol{u}_{h}-\widehat{\boldsymbol{u}}_{h}\right),
\end{equation}
where $\boldsymbol{\tau}_{h}$ is a face-based stabilization matrix. In practice, $\boldsymbol{\tau}_{h}$ combines an inviscid contribution that scales with a characteristic wave speed and a viscous contribution that scales with the effective transport coefficients and the inverse element size \cite{Moro2011a,Nguyen2012,Schutz2013,Vila-Perez2021,Ciuca2020}. The artificial-viscosity augmentation described above enters \eqref{eq:num-flux-hdg} through $\boldsymbol{F}^{v}$ and through the diffusive part of $\boldsymbol{\tau}_{h}$ when shock capturing is active. On boundary faces, $\widehat{\boldsymbol{b}}_{h}$ replaces the interior transmission condition by the weak statement of the physical boundary condition, using either prescribed exterior states, wall data, or periodic face pairing as appropriate.

\subsection{Implicit Temporal Discretization}

After spatial discretization, only the conservative state carries a time derivative, whereas the gradient variable and hybrid trace remain algebraic. Denoting by $\mathbf{Q}$, $\mathbf{U}$, and $\widehat{\mathbf{U}}$ the coefficient vectors of $\boldsymbol{q}_{h}$, $\boldsymbol{u}_{h}$, and $\widehat{\boldsymbol{u}}_{h}$, respectively, the semidiscrete HDG system may therefore be written as the index-1 differential-algebraic system \cite{Nguyen2012,Jaust2014a}
\begin{subequations}
\label{eq:semi-discrete}
\begin{align}
\label{eq:semi-diff}
\mathbf{M}\frac{d\mathbf{U}}{dt} + \mathbf{N}(\mathbf{Q},\mathbf{U},\widehat{\mathbf{U}}) &= \mathbf{0},
\\
\label{eq:semi-alg1}
\mathbf{G}(\mathbf{Q},\mathbf{U},\widehat{\mathbf{U}}) &= \mathbf{0},
\\
\label{eq:semi-alg2}
\mathbf{H}(\mathbf{Q},\mathbf{U},\widehat{\mathbf{U}}) &= \mathbf{0},
\end{align}
\end{subequations}
where $\mathbf{M}$ is the element-block mass matrix associated with the conservative variables, $\mathbf{N}$ collects the elemental conservation residuals, $\mathbf{G}$ represents the discrete gradient equations, and $\mathbf{H}$ collects the global transmission and boundary residuals. Local elimination can be applied either to the semidiscrete residual or after stage-wise linearization; in the present formulation the temporal discretization is introduced at the HDG residual level and static condensation is carried out inside each nonlinear stage solve.

Time advancement from $t^{n}$ to $t^{n+1}=t^{n}+\Delta t$ is performed with an $s$-stage diagonally implicit Runge--Kutta method \cite{alexa77,Jaust2014a}. Let $(a_{ij},b_i,c_i)$ denote the Butcher coefficients, with $a_{ij}=0$ for $j>i$, and let $d_{ij}$ be the entries of the inverse Runge--Kutta matrix $(a_{ij})^{-1}$. If $\mathbf{W}^{n,i} = (\mathbf{Q}^{n,i},\mathbf{U}^{n,i},\widehat{\mathbf{U}}^{n,i})$ denotes the $i$th stage state, then the stage equations are written as
\begin{subequations}
\label{eq:dirk-stage}
\begin{align}
\label{eq:dirk-stage-diff}
\sum_{j=1}^{s} d_{ij}\mathbf{M}\left(\mathbf{U}^{n,j}-\mathbf{U}^{n}\right)
+ \Delta t\,\mathbf{N}(\mathbf{W}^{n,i}) &= \mathbf{0},
\\
\label{eq:dirk-stage-alg1}
\mathbf{G}(\mathbf{W}^{n,i}) &= \mathbf{0},
\\
\label{eq:dirk-stage-alg2}
\mathbf{H}(\mathbf{W}^{n,i}) &= \mathbf{0},
\end{align}
\end{subequations}
for $i=1,\ldots,s$. Because the scheme is diagonally implicit, stage $i$ depends only on stages $1,\ldots,i$, and the already converged lower stages enter the residual as known data. Once the stage solutions are available, the conservative variables are updated by
\begin{equation}
\label{eq:dirk-update}
\mathbf{U}^{n+1}
=
\left(1-\sum_{j=1}^{s} e_j\right)\mathbf{U}^{n}
+
\sum_{j=1}^{s} e_j \mathbf{U}^{n,j},
\qquad
e_j = \sum_{i=1}^{s} b_i d_{ij},
\end{equation}
followed by a final evaluation of the algebraic variables consistent with $\mathbf{U}^{n+1}$. The specific Butcher tableau used in the production calculations can be inserted without changing the formulation below.

\subsection{Linearizarion of the Nonlinear System}

At the fully discrete level, each stage therefore requires the solution of an HDG system with the same algebraic structure as the semidiscrete problem, but with the conservative-state residual augmented by a DIRK mass contribution. For stage $i$, the local conservation residual on each element takes the form
\begin{equation}
\label{eq:fully-discrete-local}
\mathbf{R}^{K,n,i}_{u}
=
\sum_{j=1}^{s} d_{ij}\mathbf{M}^{K}\left(\mathbf{U}^{K,n,j}-\mathbf{U}^{K,n}\right)
+
\Delta t\,\mathbf{N}^{K}\!\left(\mathbf{Q}^{K,n,i},\mathbf{U}^{K,n,i},\widehat{\mathbf{U}}^{K,n,i}\right),
\end{equation}
while the gradient and trace residuals retain their semidiscrete form. To simplify the notation, we drop the superscripts $n$ and $i$ to write the element-local equations on each $K \in \mathcal{T}_{h}$ as residuals
\begin{equation}
\label{eq:local-residuals}
\boldsymbol{R}^{K}_{q}\left(\mathbf{Q}^{K},\mathbf{U}^{K},\widehat{\mathbf{U}}^{K}\right)=\boldsymbol{0},
\qquad
\boldsymbol{R}^{K}_{u}\left(\mathbf{Q}^{K},\mathbf{U}^{K},\widehat{\mathbf{U}}^{K}\right)=\boldsymbol{0},
\end{equation}
and the global transmission equation as
\begin{equation}
\label{eq:trace-residual}
\boldsymbol{R}_{\widehat{u}}\left(\{\mathbf{Q}^{K},\mathbf{U}^{K}\}_{K \in \mathcal{T}_{h}},\widehat{\mathbf{U}}\right)=\boldsymbol{0}.
\end{equation}
The global trace residual expresses single-valued normal numerical fluxes across interior faces and incorporates all boundary conditions on $\partial\Omega$. Because the element equations are conservative and the trace variable is shared, the resulting HDG discretization is locally conservative on every element. We collect all unknowns into
\begin{equation}
\label{eq:stage-unknown}
\mathbf{W}
=
\begin{bmatrix}
\mathbf{Q} &
\mathbf{U} &
\widehat{\mathbf{U}}
\end{bmatrix},
\end{equation}
and define the stage residual
\begin{equation}
\label{eq:nonlinear-residual}
\boldsymbol{\mathcal{F}}(\mathbf{W})
=
\begin{bmatrix}
\mathbf{R}_{q}(\mathbf{Q},\mathbf{U},\widehat{\mathbf{U}}) \\
\mathbf{R}_{u}(\mathbf{Q},\mathbf{U},\widehat{\mathbf{U}}) \\
\mathbf{R}_{\widehat{u}}(\mathbf{Q},\mathbf{U},\widehat{\mathbf{U}})
\end{bmatrix}
=
\mathbf{0}.
\end{equation}
The first block contains the elemental gradient residuals, the second contains the fully discrete elemental conservation residuals, and the third contains the face transmission and boundary residuals. In shocked flows, $\boldsymbol{\mathcal{F}}$ depends on the artificial-viscosity field through the modified viscous fluxes and any viscosity-dependent stabilization terms.

This construction defines the discrete nonlinear problem solved at every DIRK stage. Newton's method is applied stage by stage in the standard Newton--Krylov setting for implicit high-order flow solvers \cite{persson:GMRESDG,Knoll2004}. Given an iterate $\mathbf{W}^{k}$, the linearized correction $\delta \mathbf{W}^{k}$ satisfies
\begin{equation}
\label{eq:newton-linearization}
\boldsymbol{\mathcal{J}}(\mathbf{W}^{k})\delta \mathbf{W}^{k}
=
-\boldsymbol{\mathcal{F}}(\mathbf{W}^{k}),
\qquad
\boldsymbol{\mathcal{J}}(\mathbf{W}^{k})
=
\frac{\partial \boldsymbol{\mathcal{F}}}{\partial \mathbf{W}}(\mathbf{W}^{k}),
\end{equation}
and the stage state is updated as
\begin{equation}
\label{eq:newton-update}
\mathbf{W}^{k+1}
=
\mathbf{W}^{k} + \alpha^{k}\delta \mathbf{W}^{k},
\qquad 0 < \alpha^{k} \le 1.
\end{equation}
At each DIRK stage, the Newton iteration is initialized from the available time-accurate information, typically the converged state at the previous time step together with already computed lower-order DIRK stage data. Convergence is assessed through the nonlinear residual norm.

\subsection{Static Condensation}

Static condensation is carried out after linearization. Let the Newton corrections associated with one element be $\delta \mathbf{Q}^{K}$, $\delta \mathbf{U}^{K}$, and $\delta \widehat{\mathbf{U}}^{K}$. The linearized element system takes the block form
\begin{equation}
\label{eq:newton-local-block}
\begin{bmatrix}
\boldsymbol{A}^{K}_{qq} & \boldsymbol{A}^{K}_{qu} & \boldsymbol{B}^{K}_{q\widehat{u}} \\
\boldsymbol{A}^{K}_{uq} & \boldsymbol{A}^{K}_{uu} & \boldsymbol{B}^{K}_{u\widehat{u}}
\end{bmatrix}
\begin{bmatrix}
\delta \mathbf{Q}^{K} \\
\delta \mathbf{U}^{K}
\end{bmatrix}
=
-
\begin{bmatrix}
\boldsymbol{R}^{K}_{q} \\
\boldsymbol{R}^{K}_{u}
\end{bmatrix}
-
\begin{bmatrix}
\boldsymbol{B}^{K}_{q\widehat{u}} \\
\boldsymbol{B}^{K}_{u\widehat{u}}
\end{bmatrix}
\delta \widehat{\mathbf{U}}^{K},
\end{equation}
from which
\begin{equation}
\label{eq:local-elimination}
\begin{bmatrix}
\delta \mathbf{Q}^{K} \\
\delta \mathbf{U}^{K}
\end{bmatrix}
=
-
\left(\boldsymbol{A}^{K}\right)^{-1}
\begin{bmatrix}
\boldsymbol{R}^{K}_{q} \\
\boldsymbol{R}^{K}_{u}
\end{bmatrix}
-
\left(\boldsymbol{A}^{K}\right)^{-1}\boldsymbol{B}^{K}\delta \widehat{\mathbf{U}}^{K},
\quad
\boldsymbol{A}^{K}
=
\begin{bmatrix}
\boldsymbol{A}^{K}_{qq} & \boldsymbol{A}^{K}_{qu} \\
\boldsymbol{A}^{K}_{uq} & \boldsymbol{A}^{K}_{uu}
\end{bmatrix}.
\end{equation}
Substitution into the linearized trace equation yields the condensed global system
\begin{equation}
\label{eq:schur-trace}
\boldsymbol{S}\,\delta \widehat{\mathbf{U}}
=
\boldsymbol{g},
\end{equation}
where $\boldsymbol{S}$ is the Schur complement assembled from elemental contributions. After the trace correction is computed, the elemental corrections $\delta \mathbf{Q}^{K}$ and $\delta \mathbf{U}^{K}$ are recovered independently from \eqref{eq:local-elimination}.  The next section describes how the condensed trace system (\ref{eq:schur-trace}) is assembled, preconditioned, and solved efficiently on distributed GPU platforms.

\section{Massively Parallel HDG Solver}

\subsection{Solver Overview}

\begin{figure}[h]
\centering
\includegraphics[width=\textwidth]{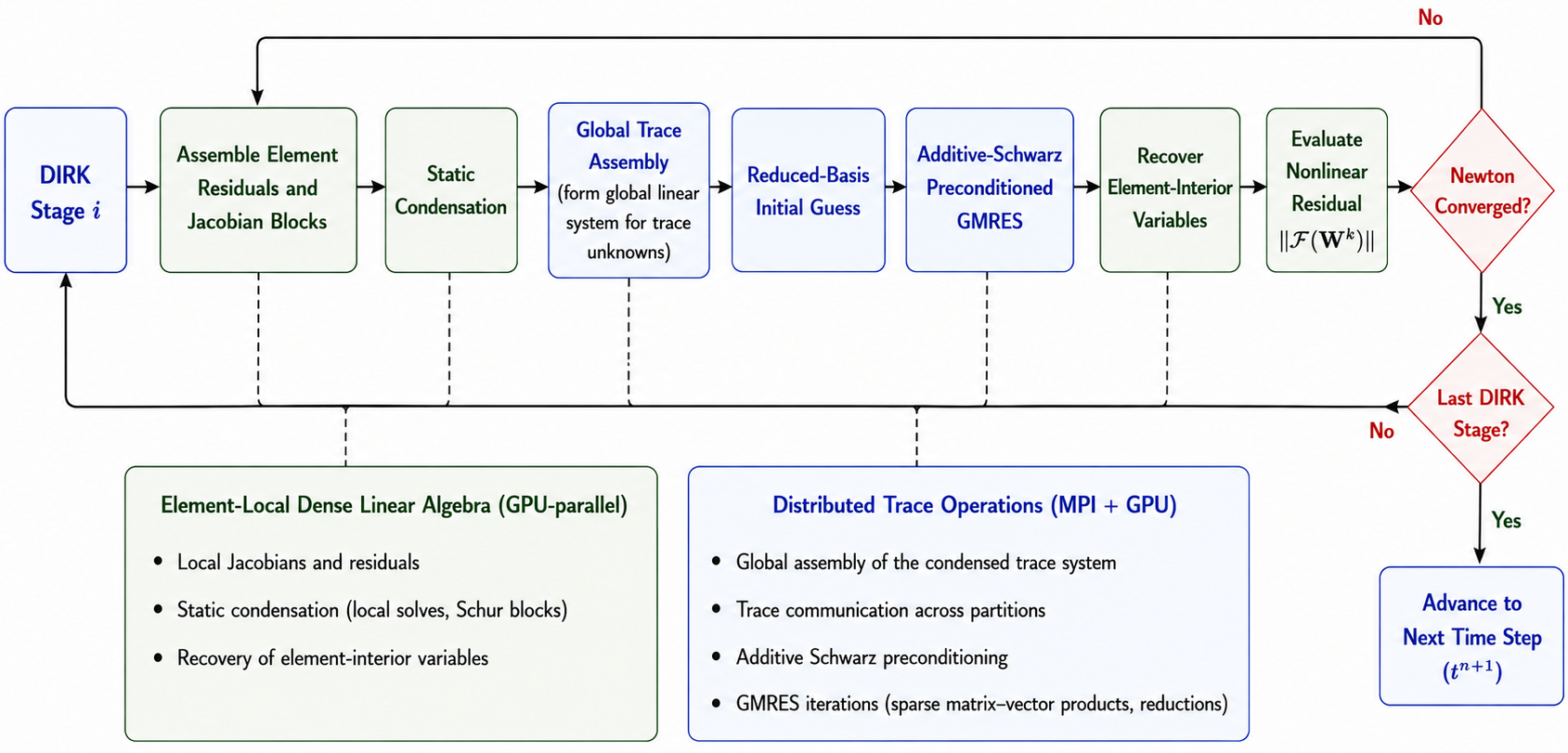}
\caption{
Computational workflow of the  massively parallel HDG solver.
At each DIRK stage, element-local residuals and Jacobian blocks are assembled and statically condensed to form the global trace system.
The condensed system is solved using a reduced-basis initialized, additive-Schwarz-preconditioned GMRES method, after which the element-interior variables are recovered and the nonlinear residual is evaluated.
Newton iterations are repeated until convergence before advancing to the next DIRK stage or time step.
The lower panels highlight the decomposition of the computation into element-local dense linear algebra, which is well suited for GPU acceleration, and distributed trace operations, which require communication across MPI partitions.
}
\label{fig:solver_overview}
\end{figure}

The implicit HDG discretization developed in Section~2 naturally decomposes the solution procedure into element-local and globally coupled computations. Following Newton linearization and static condensation, the globally coupled problem is expressed entirely in terms of the hybrid trace unknowns, while the conservative and gradient variables remain local to each element and are recovered independently after the trace correction has been computed. This decomposition forms the foundation of the proposed massively parallel solver. The resulting computational workflow is illustrated in Fig.~\ref{fig:solver_overview}.

A distinguishing feature of the HDG formulation is the separation between local dense linear algebra and global sparse communication, which is fundamental to its efficient implementation~\cite{klockner:DGGPU,Fabien2020,Karakus2019,Gibson2020, King2013}. Static condensation transforms the original coupled HDG system into a significantly smaller trace system defined on the mesh skeleton. Consequently, the dominant element-level computations consist of dense matrix factorizations, triangular solves, and matrix–matrix products associated with the local Schur complements, whereas global communication is restricted to the trace degrees of freedom shared between neighboring mesh partitions. This separation substantially reduces the size of the globally coupled problem while preserving the high-order accuracy and geometric flexibility of the HDG discretization. The proposed solver is designed to exploit this algebraic structure on distributed heterogeneous computing systems. Element-local operations are performed independently on each mesh partition and therefore expose a high degree of fine-grained parallelism suitable for batched dense linear algebra on GPUs. In contrast, the condensed trace system is distributed across MPI partitions, where communication is confined to neighboring trace unknowns. This organization naturally combines distributed-memory parallelism across compute nodes with accelerator-based execution within each node, enabling efficient utilization of modern GPU-based supercomputers.

The repeated solution of closely related condensed linear systems constitutes the dominant computational cost of the HDG solver. A restarted GMRES method provides a robust  solution strategy for the condensed HDG system. An additive Schwarz preconditioner exploits the locality of the trace coupling to accelerate Krylov convergence. And, a reduced-basis approximation constructed from previously computed solutions provides an effective initial guess for successive GMRES solves arising from Newton iterations, DIRK stages, and time steps. Finally, the entire solution procedure is implemented using distributed-memory parallelism together with portable GPU kernels and vendor-optimized dense linear algebra libraries, allowing the  solver to execute efficiently on GPU-accelerated HPC systems. The remainder of this section describes each of these components in detail.

\subsection{Parallel Newton–GMRES Solver}

The nonlinear system arising at each DIRK stage is solved using an inexact Newton–GMRES method. This choice follows the broader use of Newton–Krylov strategies for high-order DG discretizations of the compressible Navier--Stokes equations \cite{persson:GMRESDG} and the HDG-based implicit transitional-flow solver developed by Fernandez et al.~\cite{Fernandez2017a}. As described in Section~2.4, Newton linearization together with static condensation reduces each nonlinear iteration to the solution of the condensed trace system \eqref{eq:schur-trace}. Consequently, the globally coupled linear solve is performed entirely in the trace space, while the element-interior corrections are recovered locally through \eqref{eq:local-elimination} after the trace correction has been computed. The condensed trace operator is assembled once at each Newton iteration and reused throughout the corresponding Krylov solve. Assembly across partition interfaces is handled through the distributed trace representation described in Section~3.5. The condensed trace system is solved using restarted, left-preconditioned GMRES \cite{sasc86,saad93}.

\subsection{Additive Schwarz Preconditioner}
The condensed trace system is preconditioned by a one-level additive Schwarz method defined directly on the trace space \cite{Fischer1997,Cai1999,Efstathiou2003}. For HDG discretizations, this form of face-based domain decomposition is attractive because it preserves the condensed trace structure and exposes independent dense local solves; related GPU-oriented preconditioning strategies for HDG systems are discussed by Welter and Nguyen \cite{Welter2026}. Let the global trace unknowns be partitioned into subdomain trace spaces associated with element-centered face patches. If $\boldsymbol{R}_{r}$ denotes the restriction from the global trace vector to subdomain $r$, the corresponding local matrix is
\[
\boldsymbol{S}_{r}
=
\boldsymbol{R}_{r}\boldsymbol{S}\boldsymbol{R}_{r}^{T},
\]
and the preconditioner is written as
\begin{equation}
\label{eq:additive-schwarz}
\boldsymbol{M}_{\mathrm{AS}}^{-1}
=
\sum_{r=1}^{N_{s}}
\boldsymbol{R}_{r}^{T}\boldsymbol{S}_{r}^{-1}\boldsymbol{R}_{r}.
\end{equation}
The trace degrees of freedom on a shared face appear in the neighboring subdomains attached to the adjacent elements, so the preconditioner inherits the local overlap naturally induced by the HDG trace coupling. No global coarse correction is included.

This construction is closely matched to the algebra of the hybridized discretization. Each local matrix $\boldsymbol{S}_{r}$ is a dense face-level block assembled from the same condensed elemental contributions that define the global trace operator. The setup phase therefore consists of extracting the local trace blocks and factorizing them independently, whereas application of $\boldsymbol{M}_{\mathrm{AS}}^{-1}$ consists of three steps: restriction of the trace residual to the subdomain traces, independent solution of the local dense problems, and prolongation of the resulting corrections back to the global trace vector. Because the subdomain problems are independent, the expensive part of the preconditioner lies in batched dense local algebra rather than in additional global communication.

The principal role of \eqref{eq:additive-schwarz} is to reduce the number of GMRES iterations required to solve the condensed trace system while preserving the dense-block structure exposed by static condensation. This is particularly important for high-order HDG discretizations, for which the face-level blocks become increasingly expensive and the trace coupling remains compact. By defining the preconditioner on the same face patches that arise naturally in the assembly of the condensed operator, one obtains a solver in which both the global Krylov iteration and the local preconditioner application are compatible with parallel batched dense linear algebra.

\subsection{Reduced-Basis Approximation}
The condensed trace systems generated by successive Newton iterations, DIRK stages, and time steps are strongly correlated. The solver exploits that correlation through a reduced-basis approximation that provides a good initial guess for GMRES. Reduced-basis methods are widely used for low-dimensional approximation of parametrized PDEs and flow problems \cite{peterson89,prudhomme02:_reliab_real_time_solut_param,Nguyen2008,Nguyen2024,Nguyen2025}. Related reduced-basis ideas have also appeared in GPU-accelerated high-order flow solvers for hypersonic applications \cite{Nguyen2020gpu}. The present construction is closer in spirit to the reduced-subspace initialization used previously for matrix-free discontinuous Galerkin simulations of transonic buffet \cite{CuongNguyen2022}. The reduced basis is formed from previously accepted trace corrections rather than from solution snapshots of the underlying PDE, so it targets the repeated sequence of closely related linear systems that arises in the implicit HDG solve. Let
\[
\boldsymbol{W}_{n_b}
=
\begin{bmatrix}
\boldsymbol{w}_{1} & \cdots & \boldsymbol{w}_{n_b}
\end{bmatrix}
\in \mathbb{R}^{N_{\widehat{u}}\times n_b},
\qquad
n_b \le n_b^{\max},
\]
denote the current reduced basis, where each column is a previously accepted GMRES solution of the condensed trace problem, equivalently a previously accepted trace correction. Associated with this basis is the stored matrix of projected vectors
\[
\boldsymbol{Z}_{n_b}
=
\begin{bmatrix}
\boldsymbol{z}_{1} & \cdots & \boldsymbol{z}_{n_b}
\end{bmatrix}
=
\begin{bmatrix}
\boldsymbol{S}\boldsymbol{w}_{1} & \cdots & \boldsymbol{S}\boldsymbol{w}_{n_b}
\end{bmatrix}
\in \mathbb{R}^{N_{\widehat{u}}\times n_b}.
\]
Both $\boldsymbol{W}_{n_b}$ and $\boldsymbol{Z}_{n_b}$ are distributed in the same manner as the global trace vector, so no separate reduced mesh or coarse trace space is introduced.

For the current trace system \eqref{eq:schur-trace}, the reduced-basis initial guess is defined by the residual-minimization problem
\begin{equation}
\label{eq:rb-min}
\boldsymbol{c}_{n_b}
=
\arg\min_{\boldsymbol{c}\in\mathbb{R}^{n_b}}
\left\|
\boldsymbol{g}
-
\boldsymbol{Z}_{n_b}\boldsymbol{c}
\right\|_{2},
\qquad
\delta \widehat{\mathbf{U}}_{0}
=
\boldsymbol{W}_{n_b}\boldsymbol{c}_{n_b}.
\end{equation}
The resulting vector $\delta \widehat{\mathbf{U}}_{0}$ is the correction in $\operatorname{span}(\boldsymbol{W}_{n_b})$ that best approximates the current Newton step in the Euclidean residual norm. The basis is updated incrementally. After a GMRES solve has converged and the corresponding trace correction has been accepted, the newest basis vector $\boldsymbol{w}_{\mathrm{new}}$ is inserted into the reduced basis. If the maximum basis size $n_b^{\max}$ has already been reached, the oldest stored basis vector is discarded so that the basis remains a fixed window of recent corrections. Only one new projected vector is then required
$\boldsymbol{z}_{\mathrm{new}}
=
\boldsymbol{S}\boldsymbol{w}_{\mathrm{new}}$. Previously stored products $\boldsymbol{S}\boldsymbol{w}_{i}$ are retained, so the matrix $\boldsymbol{Z}_{n_b}$ is updated by appending or replacing a single column rather than being recomputed from scratch. The reduced-basis update therefore requires one application of the condensed trace operator per accepted linear solve, together with the bookkeeping needed to maintain the fixed-size window of recent columns in $\boldsymbol{W}_{n_b}$ and $\boldsymbol{Z}_{n_b}$.

\subsection{Distributed GPU Implementation}

The distributed implementation carefully exploits the decompositions induced by HDG to exploit modern accelerators on two principles. First, inter-partition communication is restricted to trace data and to the condensed face contributions needed to assemble or apply the trace operator. Second, element-interior states and element-local factorizations remain local and are dispatched in every case to appropriate GPU code, depending on the characteristics of the load. Since every aspect of our solver is designed to respect this decomposition, we naturally found or implemented kernels for each component. 

For the most part, we apply batched vendor linear algebra libraries for operations targeting dense element and face blocks. Assembly of the condensed operator, local recovery of the interior corrections, application of the additive Schwarz preconditioner, and construction of the reduced-basis initial guess are all dominated by dense linear algebra over many independent small to moderate-size blocks. These operations are well suited to accelerator hardware because they can be executed in batched form with relatively high arithmetic intensity, as also observed in GPU-oriented high-order DG and HDG solvers \cite{Nguyen2020gpu,CuongNguyen2022,Welter2026, King2013, klockner:DGGPU, King2013}.

For operations involving quadrature points within assembly where our implementation handles generic pointwise operations within the HDG weak formulation, Kokkos is used to express the custom parallel kernels in a performance-portable form \cite{Carter2014,Trott2022}. These kernels feed more traditional assembly kernels that look liked batched linear algebra. Although this strategy does require some amount of intermediate memory, to enable generality,  the implementation still manages with a small carefully managed scratchpad. We leave as future work carefully exploiting fine grained tradeoffs between fusion, data locality, and parallelism in this collection of kernels. For now, the implementation's management of memory still allow us to keep the principal solution, residual, and operator data resident in device memory during the stage solve. 

Our decomposition into Kokkos and batched linear algebra libraries matches the structure of HDG and is consistent with recent GPU realizations of implicit DG and HDG solvers for compressible-flow applications \cite{Nguyen2020gpu,Welter2026,CuongNguyen2022}: Kokkos handles the irregular but portable parallel work associated with assembly and trace-space bookkeeping, whereas the vendor libraries accelerate the dense local matrix operations that dominate the cost of condensation, preconditioning, and recovery.

Distributed communication is organized around the distinction between interior and interface work. During assembly and operator application, interface trace data are exchanged with neighboring ranks using nonblocking MPI, while independent work associated with strictly local elements proceeds concurrently. After the halo data have been received, the interface contributions are completed and accumulated into the local trace residual or trace matrix-vector product. The same face-based decomposition is used throughout the solver, so the condensed operator, the additive Schwarz preconditioner, and the reduced-basis approximation all act on a trace vector that is partitioned consistently with the mesh skeleton.

GMRES inherently requires global synchronization through inner products and residual norm evaluations, making the associated collective communication increasingly significant at large processor counts. Consequently, the efficiency of the proposed solver depends not only on accelerating the element-local dense linear algebra on GPUs, but also on reducing the number of Krylov iterations required to solve each condensed trace system. This objective is achieved through the combined use of additive Schwarz preconditioning and reduced-basis initialization, which improve the convergence of the Krylov solver while preserving the underlying HDG formulation. Together with static condensation, these components produce a solution strategy in which the dominant floating-point work remains element-local and dense, global communication is confined primarily to the distributed trace space.

\section{Subsonic Transitional Flows over the Eppler 387 Airfoil}

\subsection{Problem Description}

The first validation problem considers subsonic transitional flow over the Eppler 387 airfoil. This configuration is particularly sensitive to laminar separation, formation of a laminar separation bubble, transition of the separated shear layer, and subsequent turbulent reattachment. Accurate prediction of the aerodynamic loads therefore requires simultaneous resolution of the attached boundary layer, the separated shear layer, and the three-dimensional structures that emerge during transition. The Eppler 387 airfoil has been characterized experimentally over a range of low Reynolds numbers \cite{McGhee1988_2} and has also been used to assess high-order numerical methods for transitional aerodynamic flows \cite{Fernandez2017a}.

The flow is modeled as a calorically perfect gas with ratio of specific heats $\gamma=1.4$ and Prandtl number $Pr=0.72$. Four simulations are performed at $Re_c\in{1.0,,2.0,,3.0,,4.6}\times10^5$,
where $Re_c$ is based on the airfoil chord and freestream conditions. The freestream Mach number and angle of attack are fixed at $M_\infty=0.09$, $\alpha=6^\circ$.  The computational domain employs a body-fitted C-grid surrounding the airfoil, with the outer boundary placed sufficiently far from the surface to minimize contamination of the near-field solution. The downstream portion of the domain is extended to resolve the wake. The spanwise extent is $L_z/c=0.1$, and periodic boundary conditions are imposed on the two spanwise planes. A no-slip adiabatic wall condition is prescribed on the airfoil surface, while a compressible far-field condition is imposed on the outer boundary through the numerical flux.

The quantities of interest are the pressure coefficient $C_p$ and skin-friction coefficient $C_f$. The simulations are also examined in terms of separation, transition, and reattachment locations; mean velocity and pressure fields; Reynolds stresses and pressure fluctuations; boundary-layer thickness measures; and instantaneous vortical structures. Together, these diagnostics assess the ability of the proposed solver to simulate the transitional flow over the Eppler 387 airfoil over a wide range of Reynolds numbers.

\subsection{Simulation Setup}

\begin{figure}[t]
\centering
    \includegraphics[width=\linewidth]{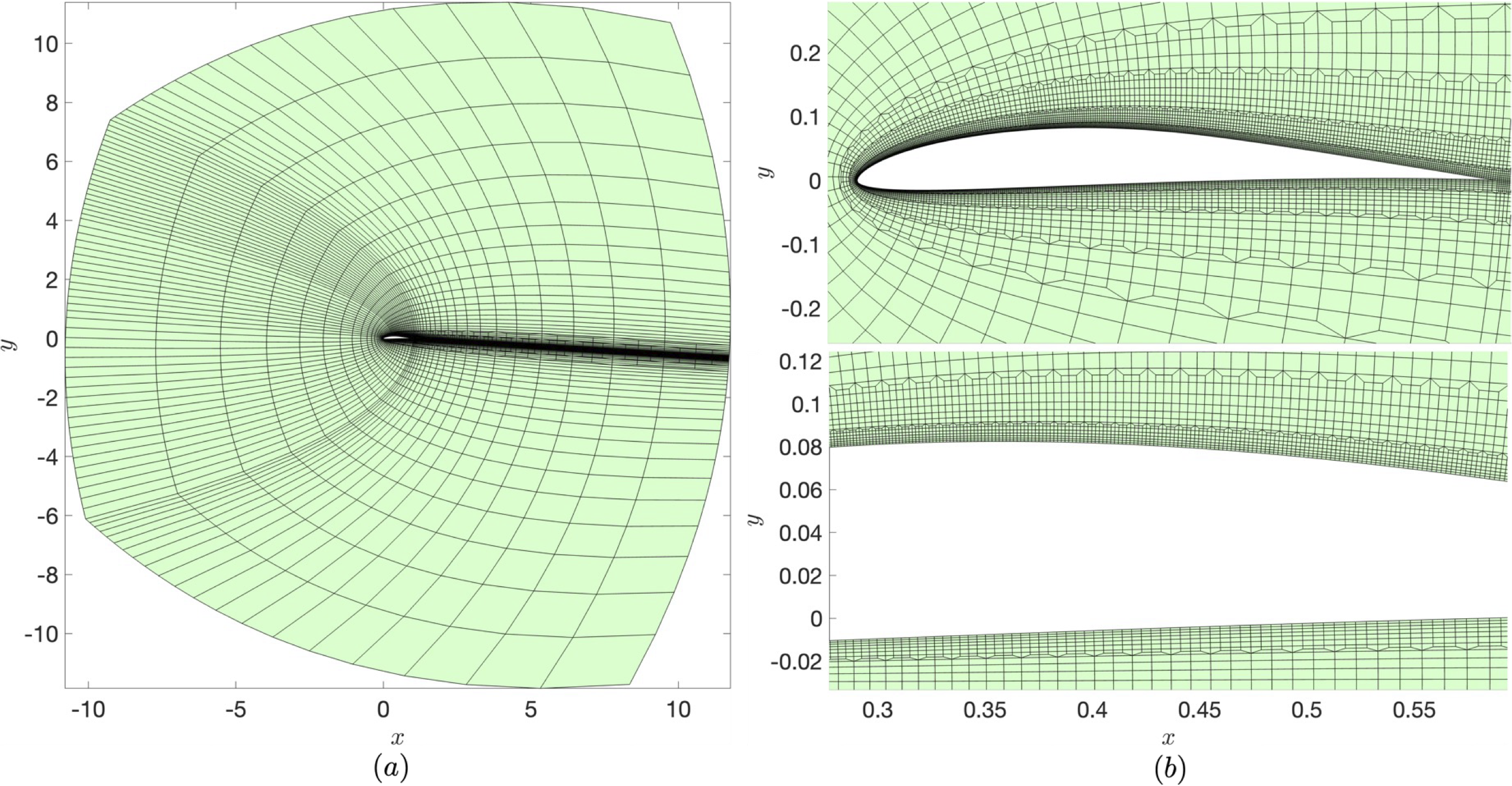}
\caption{The 2D computational mesh for the Eppler 387 airfoil. (a) Complete body-fitted curved quadrilateral C-grid used as the base mesh for the 3D simulations. (b) Enlarged views illustrate the progressive near-wall refinement, smooth mesh grading, and low-aspect-ratio elements used to  resolve the laminar boundary layer, laminar separation bubble, separated shear layer, and subsequent transition to turbulence. The 3D computational mesh is obtained by uniformly extruding this 2D mesh in the spanwise direction.}
\label{fig:eppler2dmesh}
\end{figure}

The three-dimensional computational mesh is obtained by uniformly extruding a body-fitted two-dimensional C-grid around the Eppler 387 airfoil in the spanwise direction. The underlying two-dimensional mesh consists of 20,680 curved quadrilateral elements, and its extrusion through 32 uniformly spaced spanwise layers over a spanwise extent of $L_z/c=0.1$ produces a three-dimensional mesh containing 661,760 curved hexahedral elements. Figure~\ref{fig:eppler2dmesh} shows the two-dimensional computational mesh together with an enlarged view of the near-airfoil discretization. The body-fitted C-grid topology provides smooth alignment with the airfoil geometry and an elongated wake block that enables accurate resolution of the downstream shear layer while keeping the far-field boundary sufficiently distant from the airfoil.

Polynomial degree $k=2$ is employed for both the geometric representation and the solution approximation. The mesh is progressively refined toward the airfoil surface through five successive refinement levels, producing the wall-normal resolution required to accurately resolve the viscous sublayer, the developing laminar boundary layer, and the separated shear layer while maintaining smooth element-size transitions and high-quality curved elements throughout the computational domain. Refinement is concentrated near the leading edge, along both airfoil surfaces, and in the near wake to capture the laminar separation bubble, shear-layer instabilities, and the subsequent transition to turbulence. This successive refinement strategy maintains relatively low element aspect ratios, below 2 on the suction surface and below 5 on the pressure surface, despite the extremely fine near-wall resolution required for DNS. Because transition occurs primarily on the suction surface, the mesh is intentionally refined more aggressively there, resulting in a substantially higher concentration of elements than on the pressure surface. This mesh refinement strategy provides the near-wall resolution required for direct numerical simulation while substantially reducing the total number of elements compared to uniform mesh refinement.

Each three-dimensional simulation is initialized from a converged two-dimensional precursor solution at the corresponding Reynolds number. The precursor solution is extruded uniformly in the spanwise direction, with the initial spanwise momentum set to zero. Time integration is performed using the three-stage, third-order DIRK method described in Section~2 with a fixed nondimensional time step of $\Delta t = 2.5\times10^{-3}$, and each case is advanced for 4,000 time steps. Spanwise and temporal averages are accumulated after the initial transient to obtain statistically converged flow statistics. All simulations were performed with Exasim \cite{Vila-Perez2021} on the Frontier supercomputer at the Oak Ridge Leadership Computing Facility using 92 compute nodes (736 AMD Instinct MI250X GPUs). Each Reynolds-number case required approximately two hours of wall-clock time.

\subsection{Transitional Flow Physics}
\begin{figure}[t]
\centering
\includegraphics[width=0.98\linewidth]{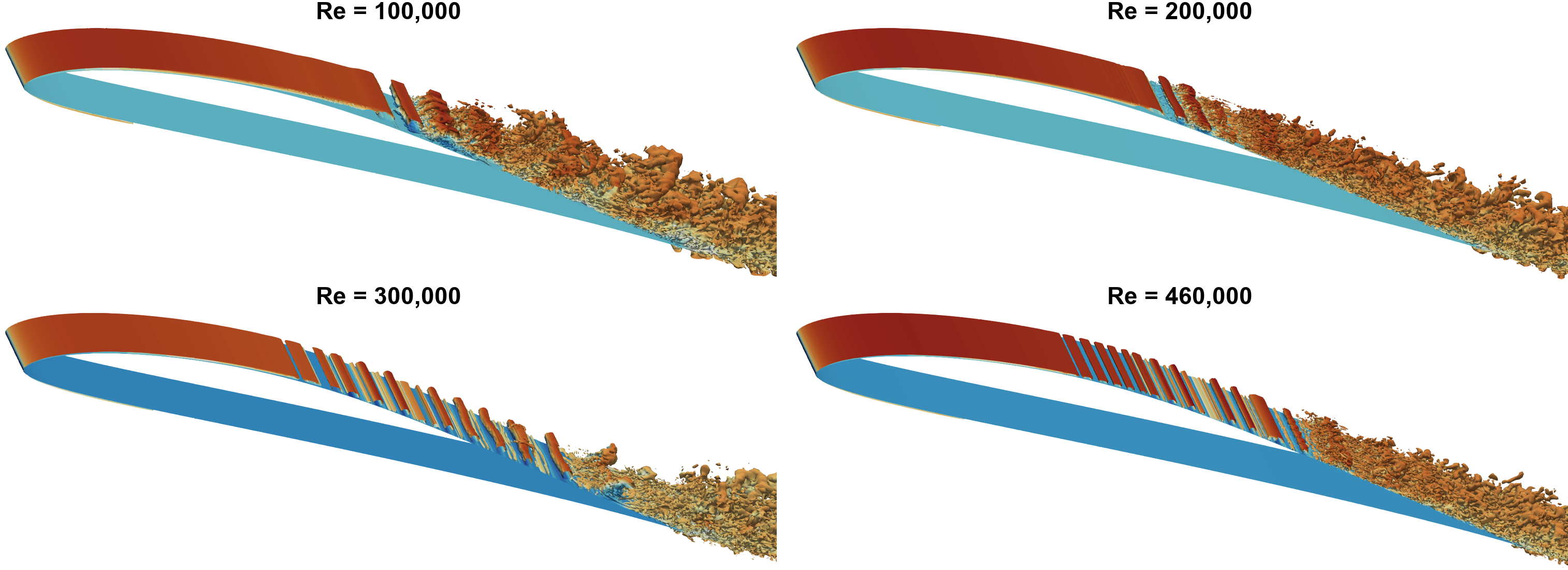}
\caption{Instantaneous vortical structures for the Eppler 387 simulations at $Re_c=1.0\times 10^{5}$, $2.0\times 10^{5}$, $3.0\times 10^{5}$, and $4.6\times 10^{5}$, identified by a positive $Q$-criterion isosurface. The panels visualize the emergence, streamwise development, and eventual breakdown of 3D structures associated with transition in the suction-side separated shear layer and wake.}
\label{fig:eppler-qcriterion}
\end{figure}

Figure~\ref{fig:eppler-qcriterion} provides an instantaneous view of the three-dimensional vortical structures that develop over the suction surface and in the near wake for the four Reynolds numbers. At the lowest Reynolds number, coherent three-dimensional structures first become visible only in the aft portion of the suction-side separated shear layer, after which they break down into a broad wake populated by disordered vortical motion. As the Reynolds number increases, the onset of these structures moves upstream and the band of organized vortices along the suction surface becomes longer and more clearly defined before breakdown. The highest-Reynolds-number case  exhibits a more compact transition process, in which the separated shear layer becomes three-dimensional earlier and the subsequent turbulent reattachment occurs over a shorter streamwise distance than in the low-Reynolds-number case.

\begin{figure}[ht]
\centering
\includegraphics[width=0.98\linewidth]{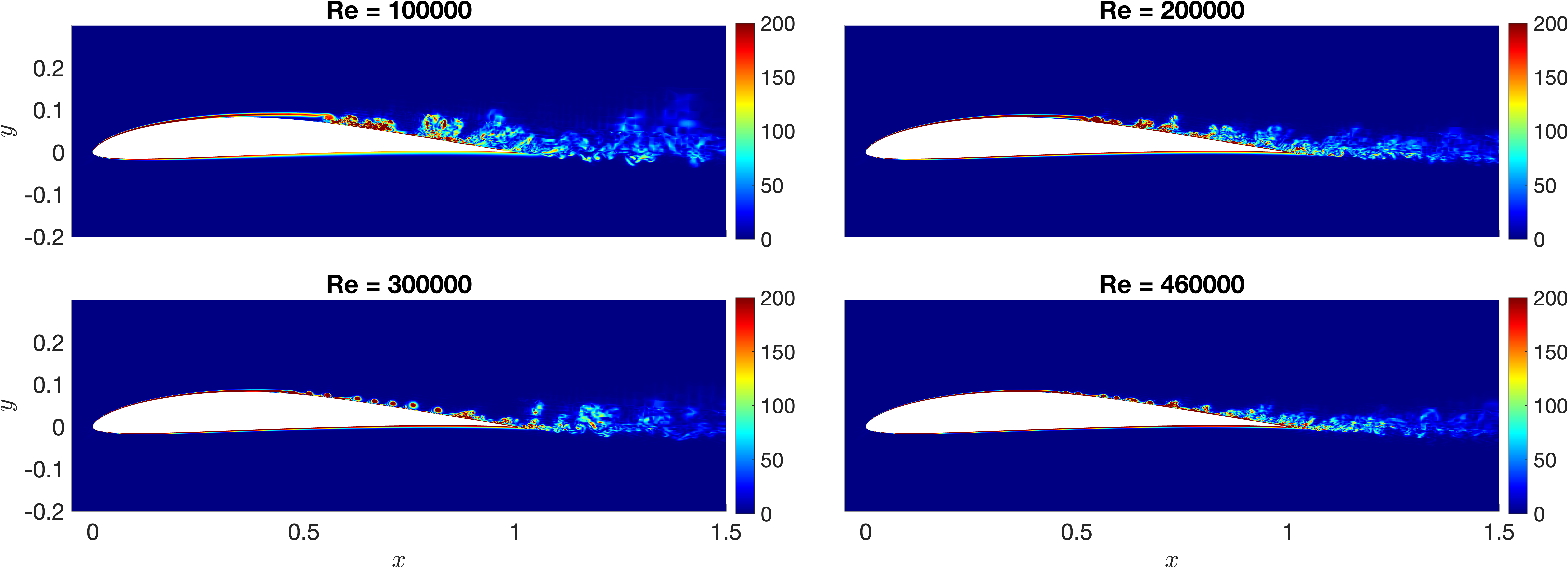}
\caption{Instantaneous vorticity-amplitude field for the four Eppler 387 simulations. The contours show where rotational disturbances are concentrated in the suction-side separated shear layer and in the near wake.}
\label{fig:eppler-vorticity-amplitude}
\end{figure}

The disturbance-amplitude field in Fig.~\ref{fig:eppler-vorticity-amplitude} reinforces this picture by showing where rotational disturbances are amplified as they convect through the separated shear layer. For $Re_c=1.0\times 10^{5}$, the strongest activity is concentrated farther downstream and is followed by a long region of disturbed wake motion. As the Reynolds number increases, the disturbance band shifts toward the suction-side separation region and becomes more closely attached to the airfoil before lifting into the wake, which is consistent with earlier instability growth and earlier turbulent breakdown. The spatial footprint of the amplified disturbances also contracts as the Reynolds number increases, indicating that the most dynamically active portion of the transitional shear layer occupies a shorter streamwise interval.

\begin{figure*}[ht]
\centering
\begin{subfigure}[t]{0.32\textwidth}
    \centering
    \includegraphics[width=\linewidth]{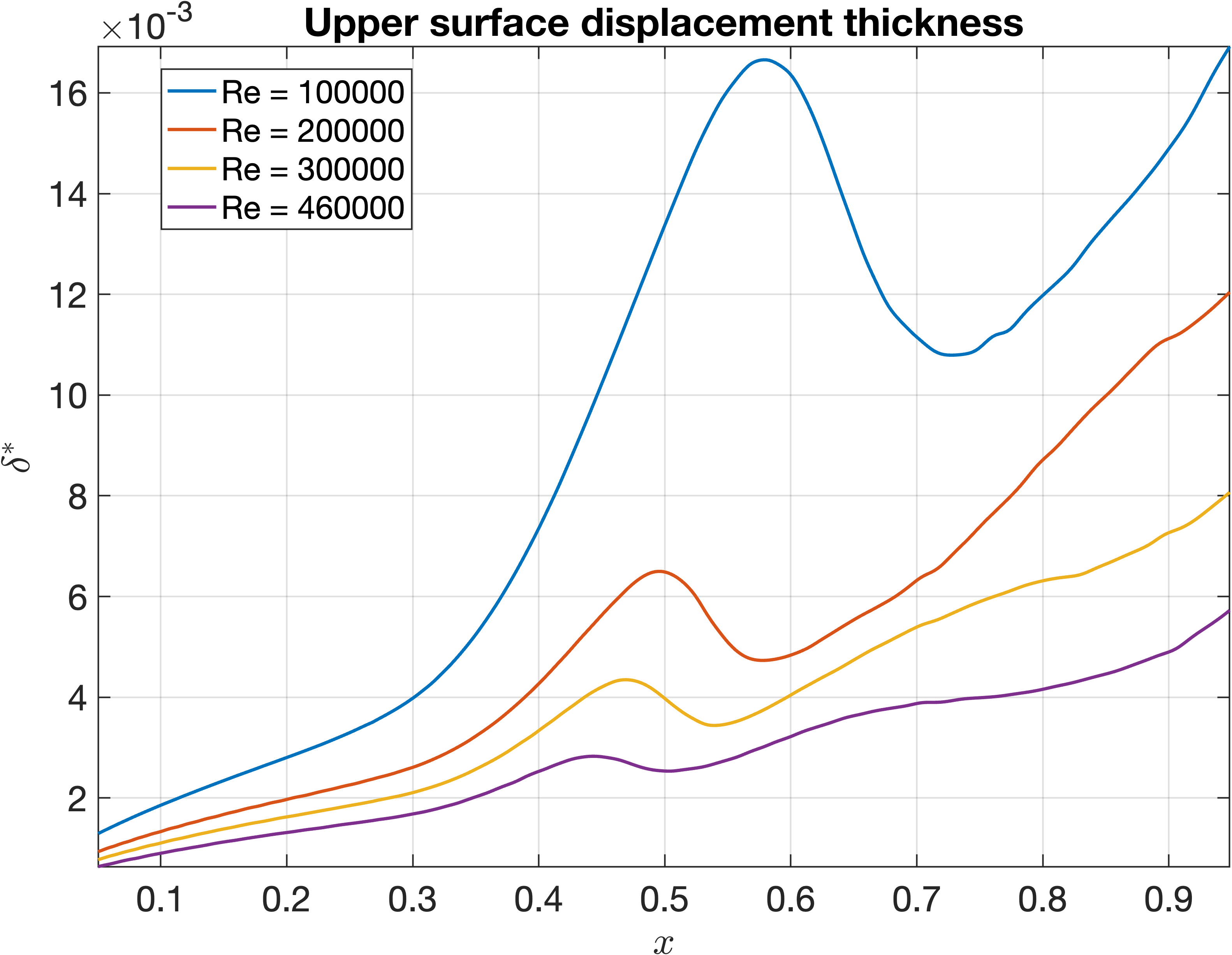}
    \caption{}
    \label{fig:eppler-upper-displacement-thickness}
\end{subfigure}
\hfill
\begin{subfigure}[t]{0.32\textwidth}
    \centering
    \includegraphics[width=\linewidth]{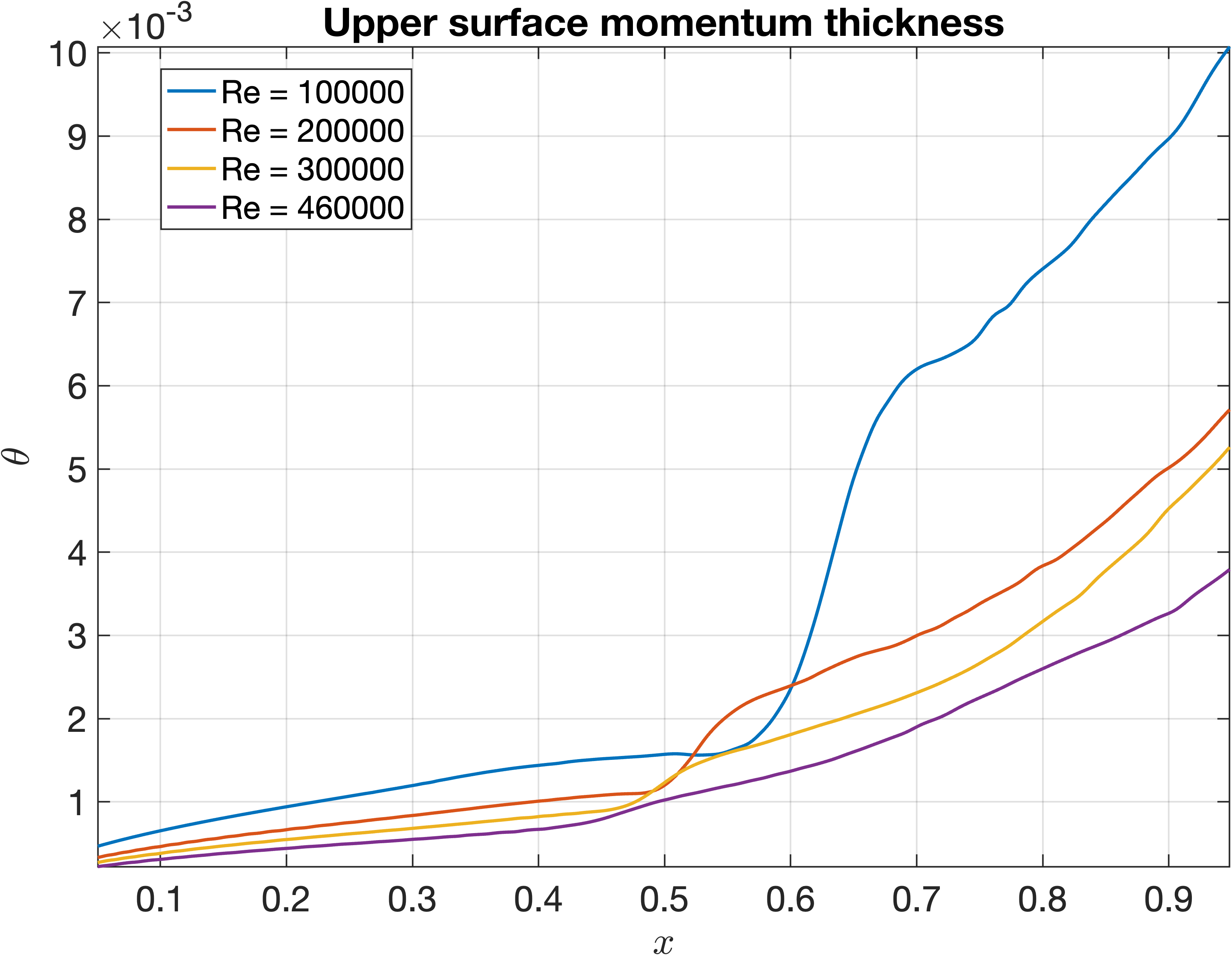}
    \caption{}
    \label{fig:eppler-upper-momentum-thickness}
\end{subfigure}
\hfill
\begin{subfigure}[t]{0.32\textwidth}
    \centering
    \includegraphics[width=\linewidth]{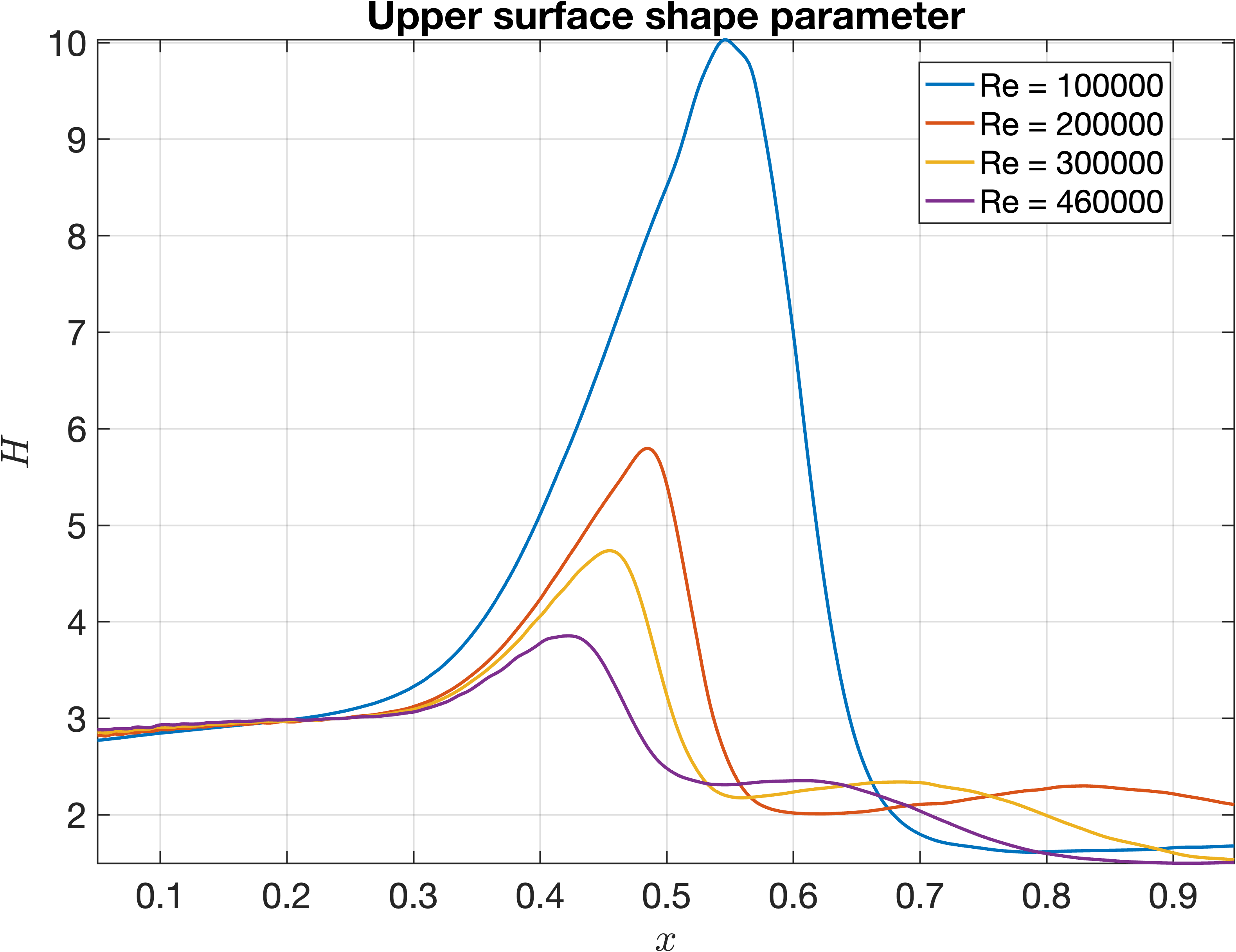}
    \caption{}
    \label{fig:eppler-upper-shape-factor}
\end{subfigure}
\caption{Upper-surface boundary-layer integral quantities for the four Eppler 387 simulations. Panel (a) shows the displacement thickness $\delta^\ast$, panel (b) shows the momentum thickness $\theta$, and panel (c) shows the shape factor $H$, each plotted as a function of chordwise coordinate on the suction surface. These quantities are widely used to characterize boundary-layer growth, laminar separation, transition, and turbulent reattachment in transitional airfoil simulations \cite{Fernandez2017a,Uranga2011b}.}
\label{fig:eppler-upper-boundary-layer-integrals}
\end{figure*}

The boundary-layer integral quantities in Fig.~\ref{fig:eppler-upper-boundary-layer-integrals} provide a complementary mean-flow description of the transition process. Displacement thickness, momentum thickness, and shape factor are widely used to characterize the growth of the boundary layer, laminar separation, transition, and turbulent reattachment, and they have previously been employed for transitional airfoil simulations by Fernandez et al.~\cite{Fernandez2017a} and Uranga et al.~\cite{Uranga2011b}. The suction-surface displacement thickness in Fig.~\ref{fig:eppler-upper-displacement-thickness} grows rapidly through the separated region and reaches its largest values in the lowest-Reynolds-number case, which is consistent with the more extended low-speed region visible in the instantaneous fields. The momentum thickness in Fig.~\ref{fig:eppler-upper-momentum-thickness} exhibits the same ordering, with the lowest-Reynolds-number case undergoing the largest downstream growth and the highest-Reynolds-number case remaining comparatively thin throughout the suction surface. The shape factor $H$ in Fig.~\ref{fig:eppler-upper-shape-factor} is particularly useful for identifying the separated and recovering portions of the boundary layer: all four cases show a pronounced increase through the laminar separation-bubble region followed by a sharp decrease once transition and reattachment occur. Increasing Reynolds number reduces both the peak value of $H$ and the streamwise extent over which elevated shape factor persists, indicating that the separated bubble becomes shorter and the recovering turbulent boundary layer is established earlier.

\begin{figure}[ht]
\centering
\includegraphics[width=0.99\linewidth]{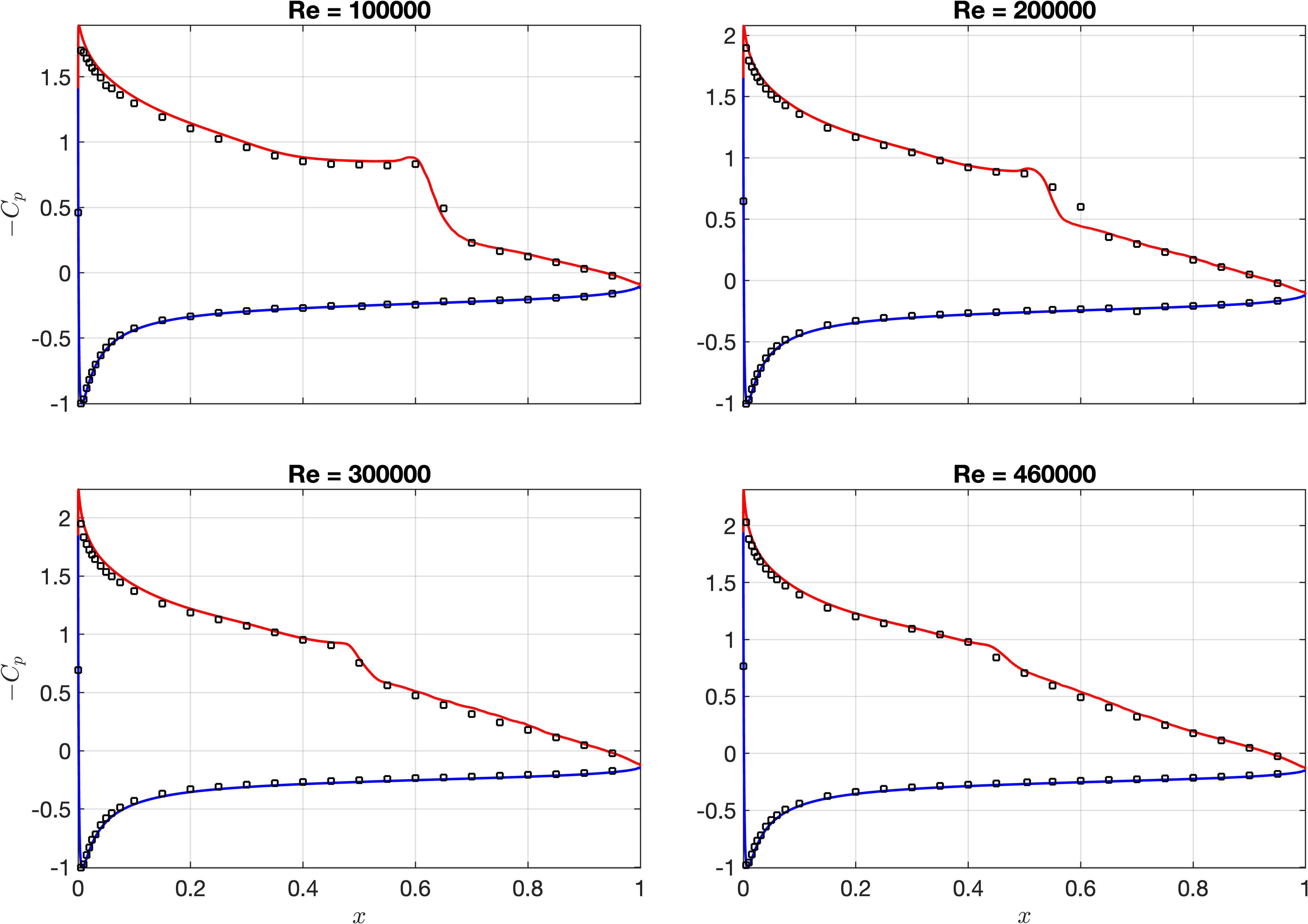}
\caption{Pressure coefficient $C_p$ along the airfoil surface for the four Eppler 387 cases at $Re_c=1.0\times 10^{5}$, $2.0\times 10^{5}$, $3.0\times 10^{5}$, and $4.6\times 10^{5}$. The solid curves are the present HDG predictions, and the circular symbols are the experimental measurements reported by McGhee et al.~\cite{McGhee1988_2}. The comparison is used to assess the Reynolds-number dependence of the suction-side pressure distribution and its recovery.}
\label{fig:eppler-Cp}
\end{figure}

\begin{figure}[ht]
\centering
\includegraphics[width=0.99\linewidth]{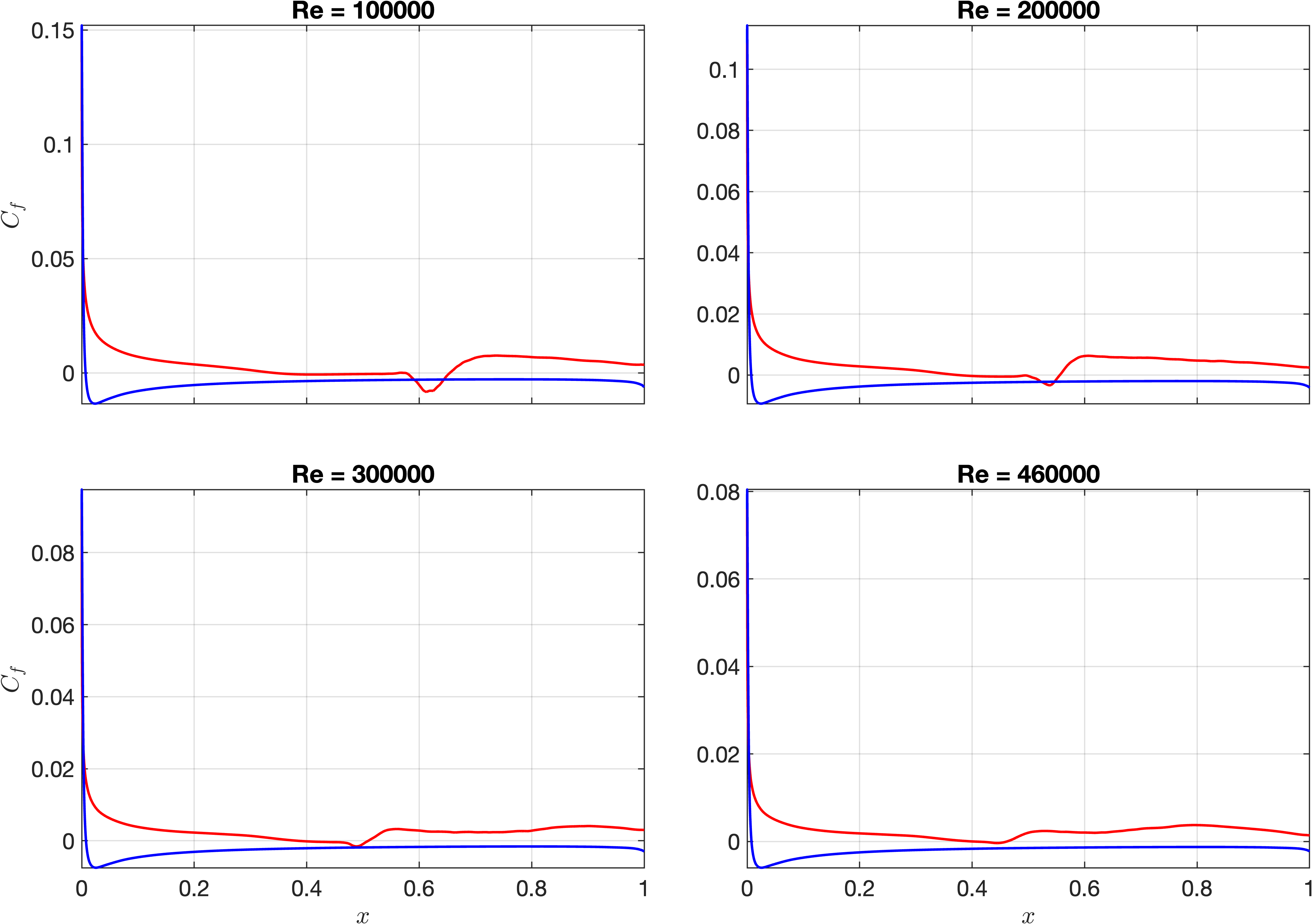}
\caption{Skin-friction coefficient $C_f$ along the airfoil surface for the four Eppler 387 simulations. The curves provide a surface-based indicator of separation, transition, and reattachment through the streamwise variation and sign changes of the wall shear.}
\label{fig:eppler-Cf}
\end{figure}

The surface pressure and skin-friction distributions in Figs.~\ref{fig:eppler-Cp} and \ref{fig:eppler-Cf} show how these internal flow changes are reflected in the aerodynamic response. In Fig.~\ref{fig:eppler-Cp}, the solid curves denote the present HDG predictions and the circles denote the experimental measurements of McGhee et al.~\cite{McGhee1988_2}. Across the four Reynolds numbers, the predicted pressure distributions reproduce the  experimental data accurately, including the leading-edge suction peak, the downstream pressure recovery, and the Reynolds-number dependence of the suction-side plateau. The skin-friction coefficient in Fig.~\ref{fig:eppler-Cf} shows the same trend more directly: the lowest-Reynolds-number case has the broadest region of near-zero or negative $C_f$, whereas the higher-Reynolds-number cases recover to positive $C_f$ over a shorter streamwise interval. The marked change in $C_f$ through the middle of the suction surface is therefore consistent with the transition and reattachment trends inferred from the integral thicknesses and the instantaneous vortical structures.

\begin{figure}[ht]
\centering
\includegraphics[width=0.98\linewidth]{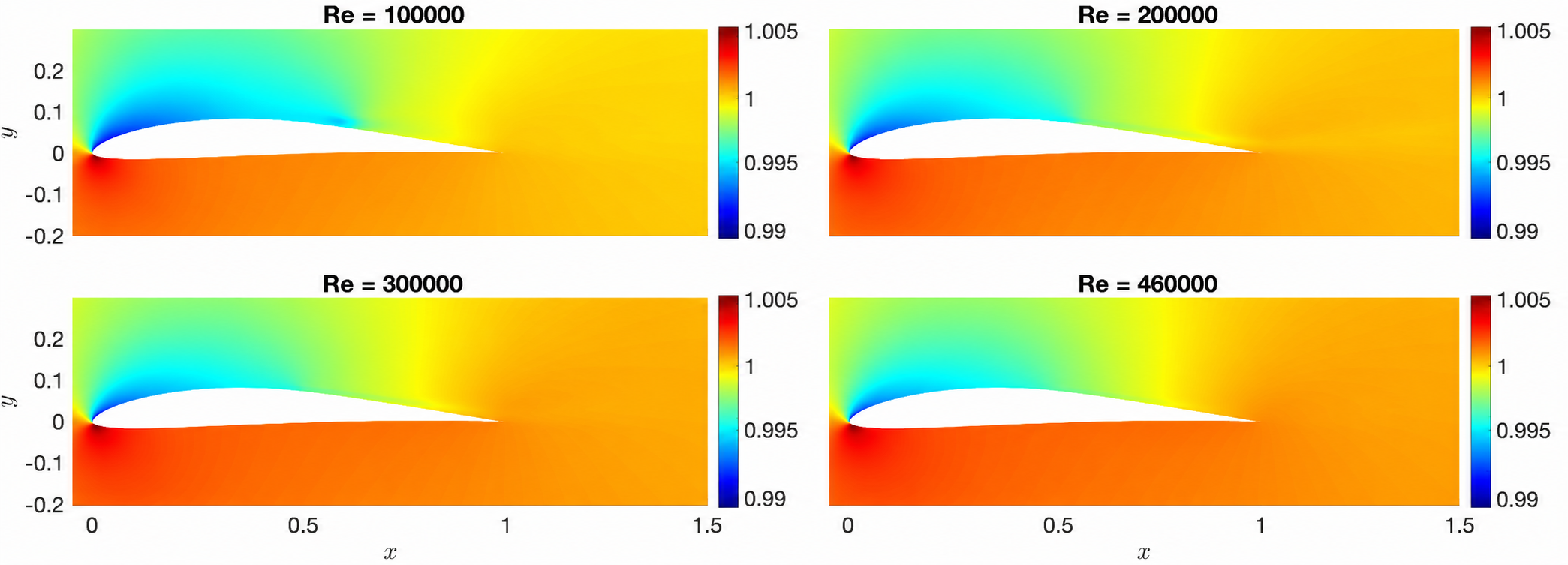}
\caption{Mean pressure normalized by the freestream pressure for the four Eppler 387 simulations.}
\label{fig:eppler-mean-pressure}
\end{figure}

\begin{figure}[ht]
\centering
\includegraphics[width=0.98\linewidth]{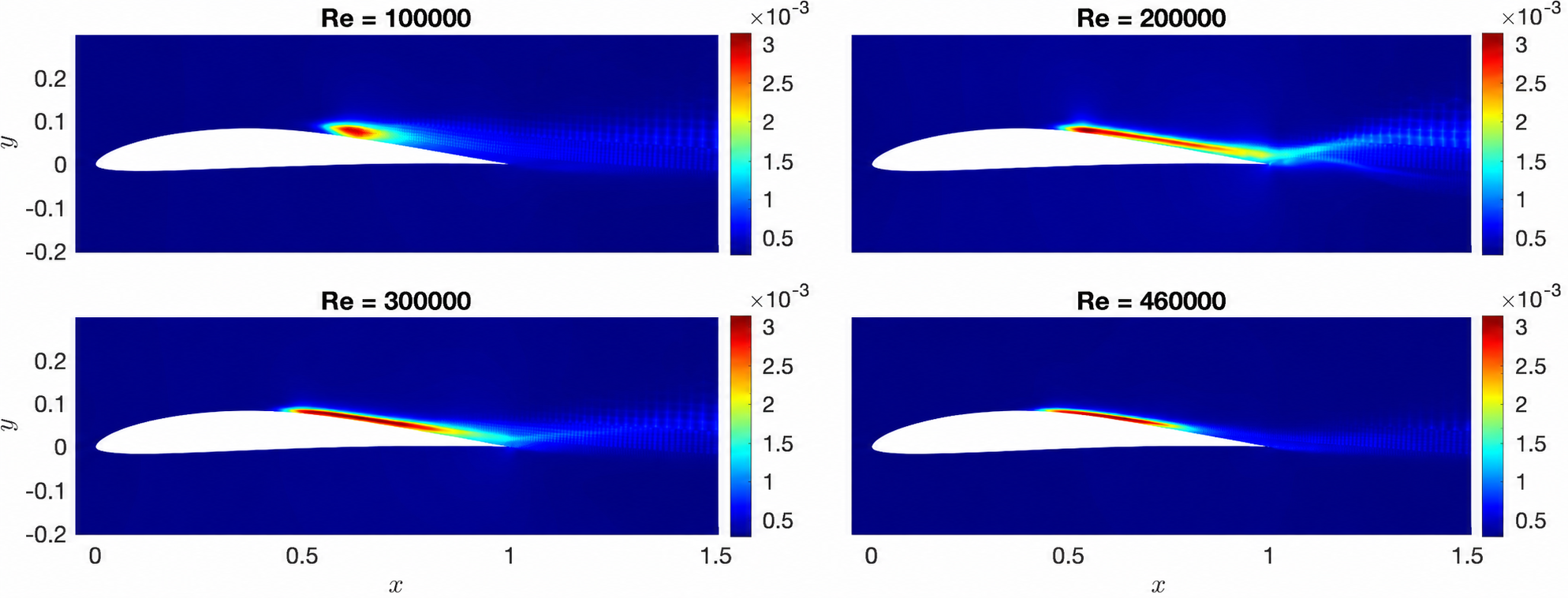}
\caption{Root-mean-square pressure fluctuation normalized by the freestream pressure for the four Eppler 387 simulations. The contours identify the suction-side and wake regions in which transitional unsteadiness is most pronounced.}
\label{fig:eppler-pressure-rms}
\end{figure}

The mean and fluctuating pressure fields in Figs.~\ref{fig:eppler-mean-pressure} and \ref{fig:eppler-pressure-rms} further localize the dynamically active region. The normalized mean pressure field remains smooth over most of the domain, but the suction-side pressure recovery and wake pressure footprint become more compact as the Reynolds number increases. By contrast, the pressure RMS field is strongly localized to the separated shear layer and the immediate wake, where transition-induced unsteadiness is concentrated. In all four cases the most active pressure-fluctuation region follows the separated shear layer from the aft suction surface into the wake, while its streamwise footprint moves upstream and shortens as the Reynolds number increases.

\begin{figure}[ht]
\centering
\includegraphics[width=0.98\linewidth]{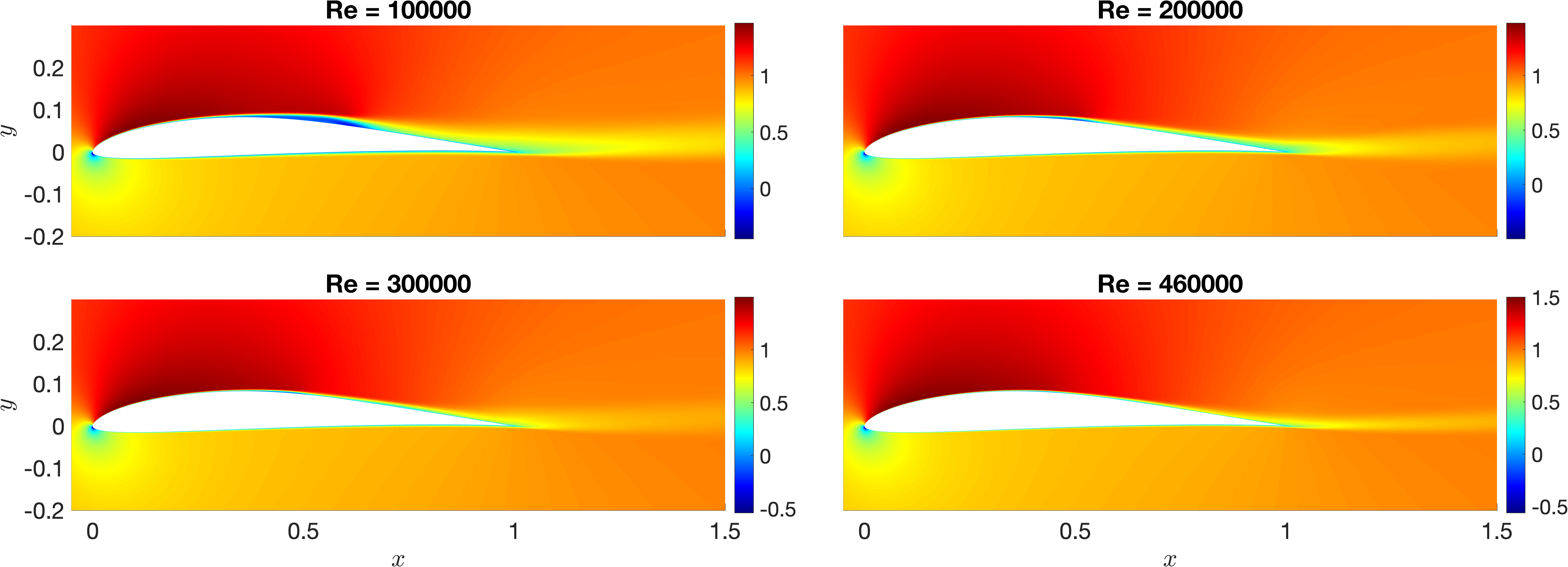}
\caption{Mean streamwise velocity field for the four Eppler 387 simulations. The contours show the development of the low-speed suction-side layer and its downstream recovery.}
\label{fig:eppler-mean-velocity}
\end{figure}

\begin{figure}[ht]
\centering
\includegraphics[width=0.98\linewidth]{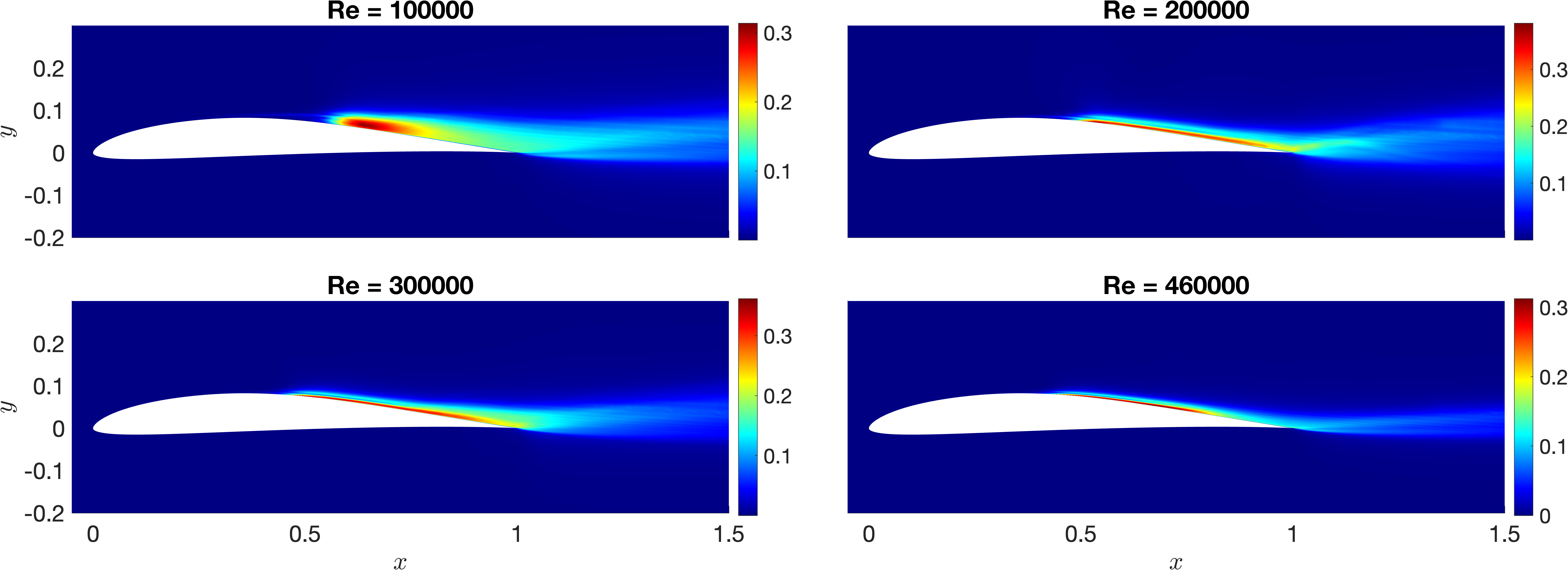}
\caption{Root-mean-square streamwise-velocity fluctuation for the four Eppler 387 simulations. The contours identify the growth, downstream convection, and decay of velocity fluctuations through the transitional shear layer and wake.}
\label{fig:eppler-velocity-rms}
\end{figure}

Figures~\ref{fig:eppler-mean-velocity} and \ref{fig:eppler-velocity-rms} show the corresponding mean and fluctuating velocity fields. The mean streamwise-velocity contours reveal a low-speed layer above the suction surface that thickens through the separated region and then contracts after reattachment, with the thickest and most extended deficit occurring at the lowest Reynolds number. The RMS field is confined to the same separated-shear-layer/wake corridor identified previously by the pressure fluctuations and vorticity amplitude. Its peak activity appears near the aft suction surface and then decays downstream, again with a progressively shorter streamwise extent as the Reynolds number increases.

\begin{figure}[t]
\centering
\includegraphics[width=0.98\linewidth]{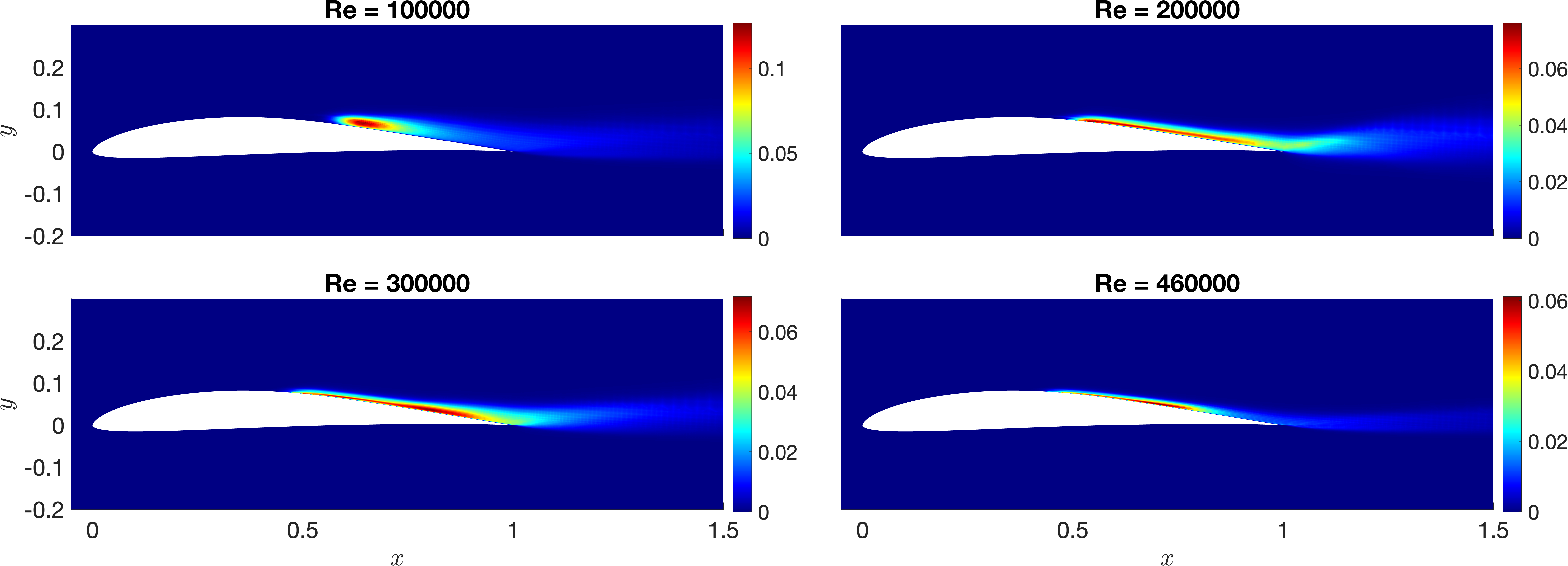}
\caption{Turbulent kinetic energy field for the four Eppler 387 simulations. The contours show where the resolved turbulent activity is concentrated in the separated shear layer and wake.}
\label{fig:eppler-tke}
\end{figure}

The turbulent kinetic energy field in Fig.~\ref{fig:eppler-tke} summarizes these trends in a single quantity. The energy is concentrated in the same separated-shear-layer band highlighted by the RMS fields and the $Q$-criterion visualizations, and it subsequently decays in the wake as the turbulence spreads downstream. The streamwise location of the energy-containing region shifts upstream with Reynolds number, while its spatial extent becomes shorter and more tightly confined to the suction-side transition zone. This behavior is consistent with the earlier onset of three-dimensionality, the shorter separation bubble, and the more rapid turbulent recovery inferred from the preceding figures.

\section{Supersonic Taylor--Green Vortex}

\subsection{Problem Description}

The supersonic Taylor--Green vortex provides a canonical periodic configuration for assessing whether a compressible flow solver can represent the transition from smooth large-scale vortical motion to fully three-dimensional turbulence while simultaneously capturing shock-associated compressibility effects. In contrast to wall-bounded transition problems, the geometry and boundary conditions are deliberately simple, so the case isolates the interaction between numerical dissipation, shock capturing, and the multiscale dynamics of compressible turbulence \cite{Gassner2013,Fernandez2017}. The present study follows the supersonic extension of the classical Taylor--Green vortex introduced for the assessment of low-dissipation shock-capturing strategies and subsequently adopted for cross-comparisons of high-order numerical methodologies for this flow \cite{Lusher2020,Chapelier2024}.

The problem is posed on the triply periodic cube $\Omega = (0 , 2\pi L)^3$, where $L = 1$ is the reference length scale of the initial vortices. The initial condition is specified analytically in terms of the primitive variables. Using the notation $(x,y,z)$ for the Cartesian coordinates and $(v_1,v_2,v_3)$ for the velocity components, the initial velocity field is
\begin{equation}
\label{eq:tgv-initial-velocity}
\begin{aligned}
v_1(x,y,z,0) &=  \sin\!\left(\frac{x}{L}\right)\cos\!\left(\frac{y}{L}\right)\cos\!\left(\frac{z}{L}\right), \\
v_2(x,y,z,0) &= - \cos\!\left(\frac{x}{L}\right)\sin\!\left(\frac{y}{L}\right)\cos\!\left(\frac{z}{L}\right), \\
v_3(x,y,z,0) &= 0,
\end{aligned}
\end{equation}
and the temperature is initially uniform,
\begin{equation}
\label{eq:tgv-initial-temperature}
T(x,y,z,0)=1.
\end{equation}
The pressure field is chosen so that the initial state is compatible with the compressible equations,
\begin{equation}
\label{eq:tgv-initial-pressure}
p(x,y,z,0)=
\frac{1}{\gamma M_0^2}
+\frac{1}{16}
\left[\cos\!\left(\frac{2x}{L}\right)+\cos\!\left(\frac{2y}{L}\right)\right]
\left[\cos\!\left(\frac{2z}{L}\right)+2\right].
\end{equation}
and the initial density follows from the ideal-gas relation. The resulting flow is therefore isothermal at $t=0$, but it is not spatially uniform in density or pressure. Starting from this smooth large-scale vortex field, the solution develops progressively finer vortical structures, followed by strong compressive regions, shock formation, and sustained shock--turbulence interaction.

The simulations considered in this paper adopt the supersonic benchmark parameters
\begin{equation}
\label{eq:tgv-parameters}
Re = 1600,
\qquad
M_0 = 1.25,
\end{equation}
where the Reynolds number is defined from the reference density $\rho_0$, velocity scale $U_0$, length scale $L$, and reference viscosity $\mu_0$, and the Mach number is defined from $U_0$ and the reference acoustic speed associated with the initial thermodynamic state. The dynamic viscosity is taken to be a function of temperature through Sutherland's law, and the Prandtl number is fixed at $Pr=0.71$. This regime is materially different from the low-Mach Taylor--Green vortex often used to study under-resolved turbulence. At $M_0=1.25$, compressibility and dilatational dissipation are no longer weak corrections to an essentially incompressible decay process. Instead, the evolving vortex field generates shock waves that interact with one another and with the developing turbulent structures, making the case a stringent test of both the spatial discretization and the localized stabilization strategy. Following \cite{Lusher2020}, the temperature-dependent viscosity is written in nondimensional form as
\begin{equation}
\label{eq:tgv-sutherland}
\mu(T)=\frac{1.4042\,T^{3/2}}{T+0.40417},
\end{equation}
which is used together with the constant Prandtl number $Pr=0.71$ in the viscous and heat-conduction terms.

The principal quantities examined later are the volume-integrated kinetic energy, the solenoidal dissipation, and the dilatational dissipation, together with instantaneous flow visualizations that reveal the coupled evolution of vortical structures and shocks. Following \cite{Lusher2020}, the kinetic energy is defined by
\begin{equation}
\label{eq:tgv-kinetic-energy}
E_k
=
\frac{1}{2|\Omega|}
\int_\Omega \rho v_j v_j\, d\Omega,
\end{equation}
and the total viscous dissipation rate is decomposed as
\begin{equation}
\label{eq:tgv-total-dissipation}
\varepsilon_T=\varepsilon_S+\varepsilon_D
=
\frac{1}{Re |\Omega| }\int_\Omega \mu\,\omega_i\omega_i\,d\Omega
+
\frac{4}{3Re |\Omega|}\int_\Omega \mu\left(\frac{\partial v_j}{\partial x_j}\right)^2 d\Omega,
\end{equation}
where $\omega_i=\epsilon_{ijk}\,\partial v_k/\partial x_j$ denotes the vorticity. These diagnostics provide complementary information. The kinetic energy primarily reflects the large-scale decay of the flow, the solenoidal dissipation is sensitive to the development of the vortical turbulence cascade, and the dilatational dissipation measures the intensity of compressive activity associated with shock formation and shock--turbulence interaction. Because the initial condition is analytical, the domain is periodic, and well-resolved comparison data are available from previous studies, the supersonic Taylor--Green vortex offers a controlled setting in which to evaluate whether the present implicit HDG solver can recover the correct transition to compressible turbulence without excessive numerical damping of either the vortical or shock-dominated features.

\subsection{Simulation Setup}

The supersonic Taylor--Green vortex is discretized on the periodic cube $\Omega=(0,2\pi)^3$ using uniform tensor-product hexahedral meshes. Opposite faces of the cube are paired periodically in all three coordinate directions, and the HDG approximation uses polynomial degree $k=2$ on every mesh.  Shock capturing is provided by the scalar artificial-viscosity regularization introduced in Section~2. For a mesh with $N$ elements in each coordinate direction, the characteristic sensor length is taken as $h_m=\frac{2\pi}{N}$. The artificial-viscosity coefficient is set to $2.0\times 10^{-3}$, the HDG stabilization parameter is $\tau=5.0$, and the lagged artificial-viscosity treatment described earlier is used together with two smoothing iterations for the auxiliary artificial-viscosity field. A mesh sequence with $N^3$ elements is considered for
\[
N\in\{32,48,64,80,96,112,128\},
\]
so that the study includes the resolutions $32^3$, $48^3$, $64^3$, $80^3$, $96^3$, $112^3$, and $128^3$. The corresponding simulations were performed on 4, 8, 16, 24, 36, 45, and 57 nodes (4 AMD MI250X GPUs per node), respectively, of the Frontier supercomputer at the Oak Ridge Leadership Computing Facility. Each simulation required approximately one hour of wall-clock time. 

Time integration is performed with the three-stage, third-order DIRK scheme used throughout the paper. All computations use the same fixed nondimensional time step, $\Delta t=10^{-2}$, and are advanced to the final time $t=20$, corresponding to $2000$ time steps. The kinetic energy and the solenoidal and dilatational contributions to the viscous dissipation, defined in \eqref{eq:tgv-kinetic-energy}--\eqref{eq:tgv-total-dissipation}, are evaluated as volume-integrated diagnostic quantities at every time step. These histories form the principal basis for the grid-resolution comparisons reported below, while instantaneous flow visualizations are used to examine the development of shocklets and three-dimensional vortical structures.

\subsection{Compressible Turbulence}
\begin{figure*}[htbp]
\centering
\begin{subfigure}[t]{0.49\textwidth}
    \centering
    \includegraphics[width=\linewidth]{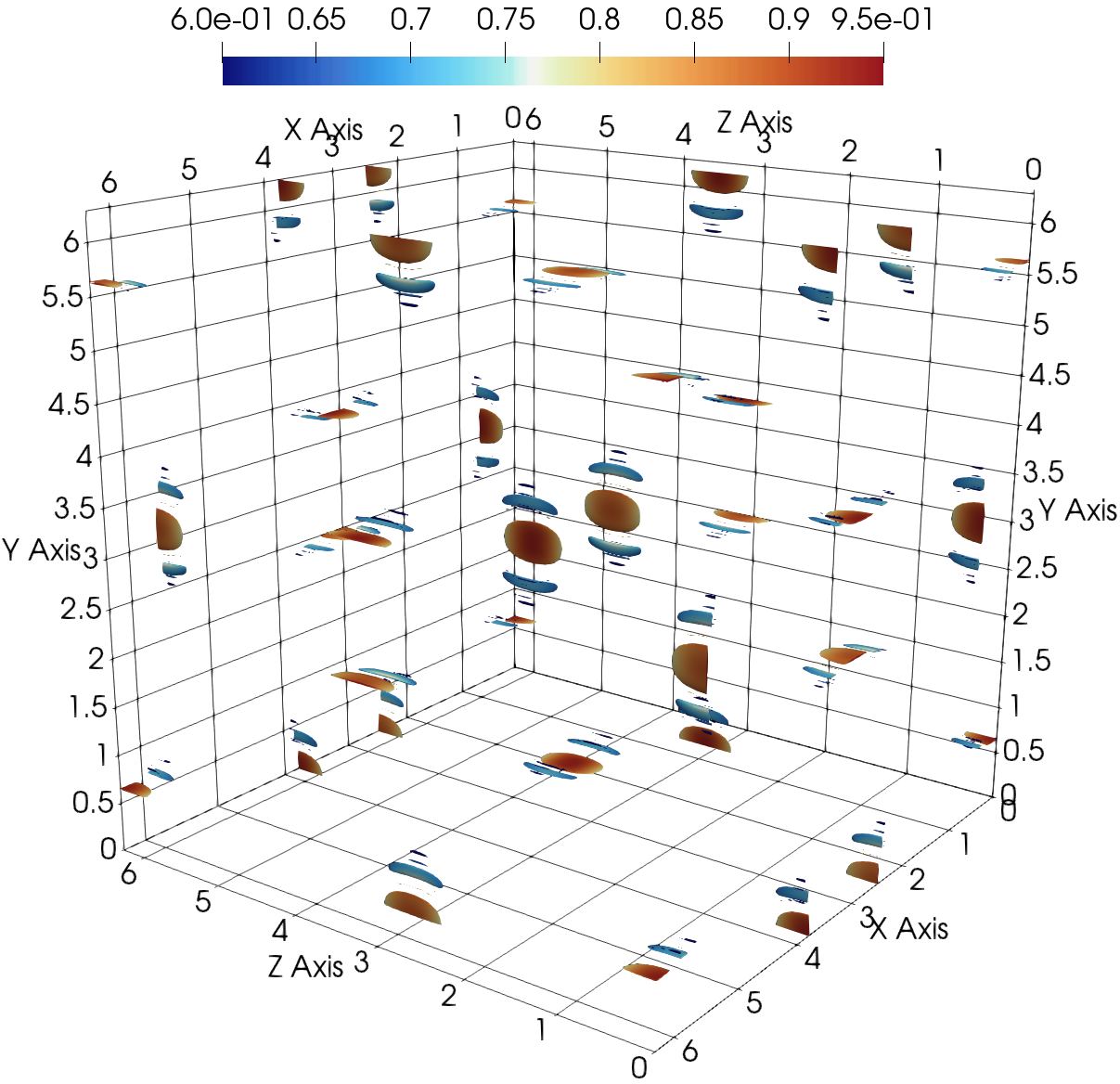}
    \caption{}
\end{subfigure}
\hfill
\begin{subfigure}[t]{0.49\textwidth}
    \centering
    \includegraphics[width=\linewidth]{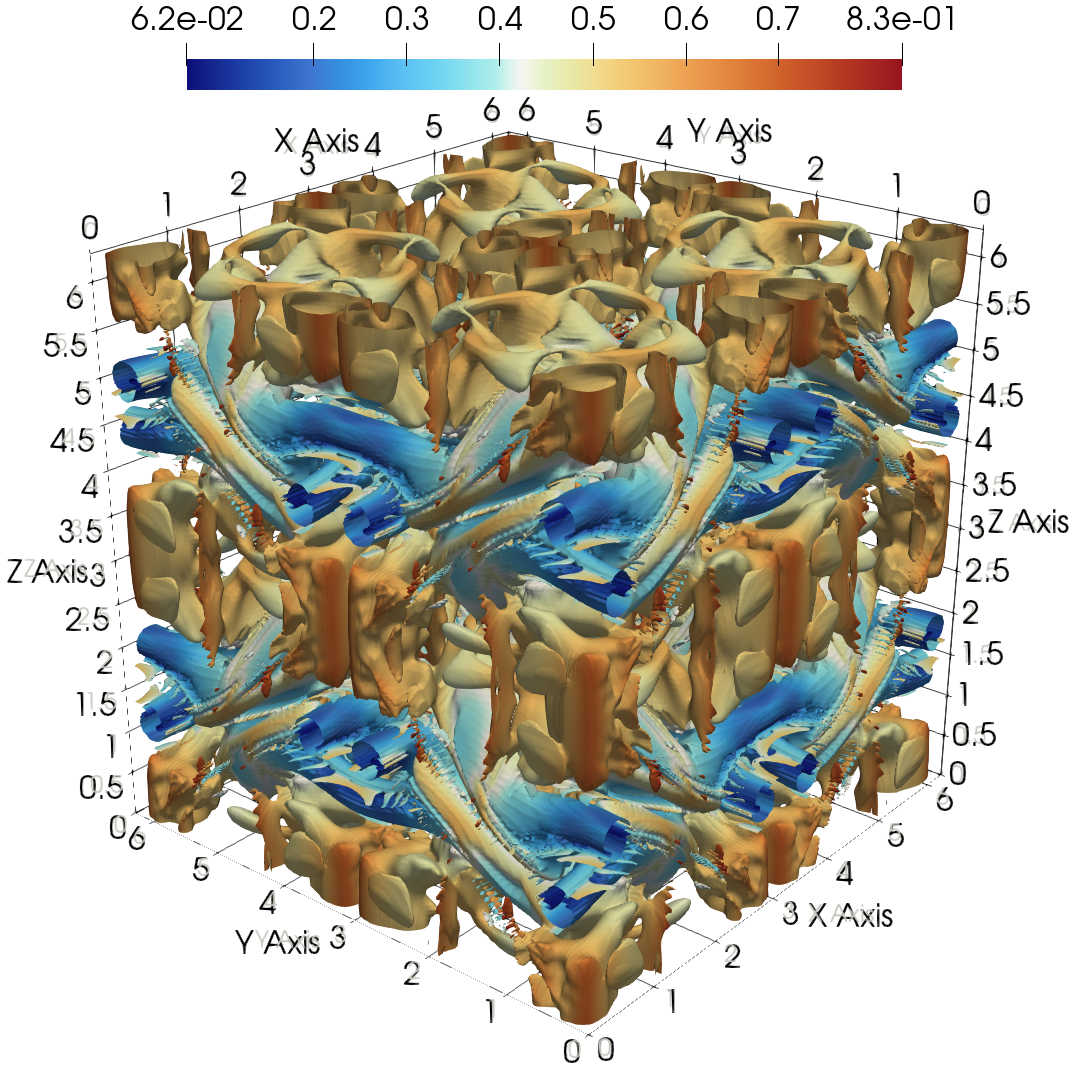}
    \caption{}
\end{subfigure}

\medskip

\begin{subfigure}[t]{0.49\textwidth}
    \centering
    \includegraphics[width=\linewidth]{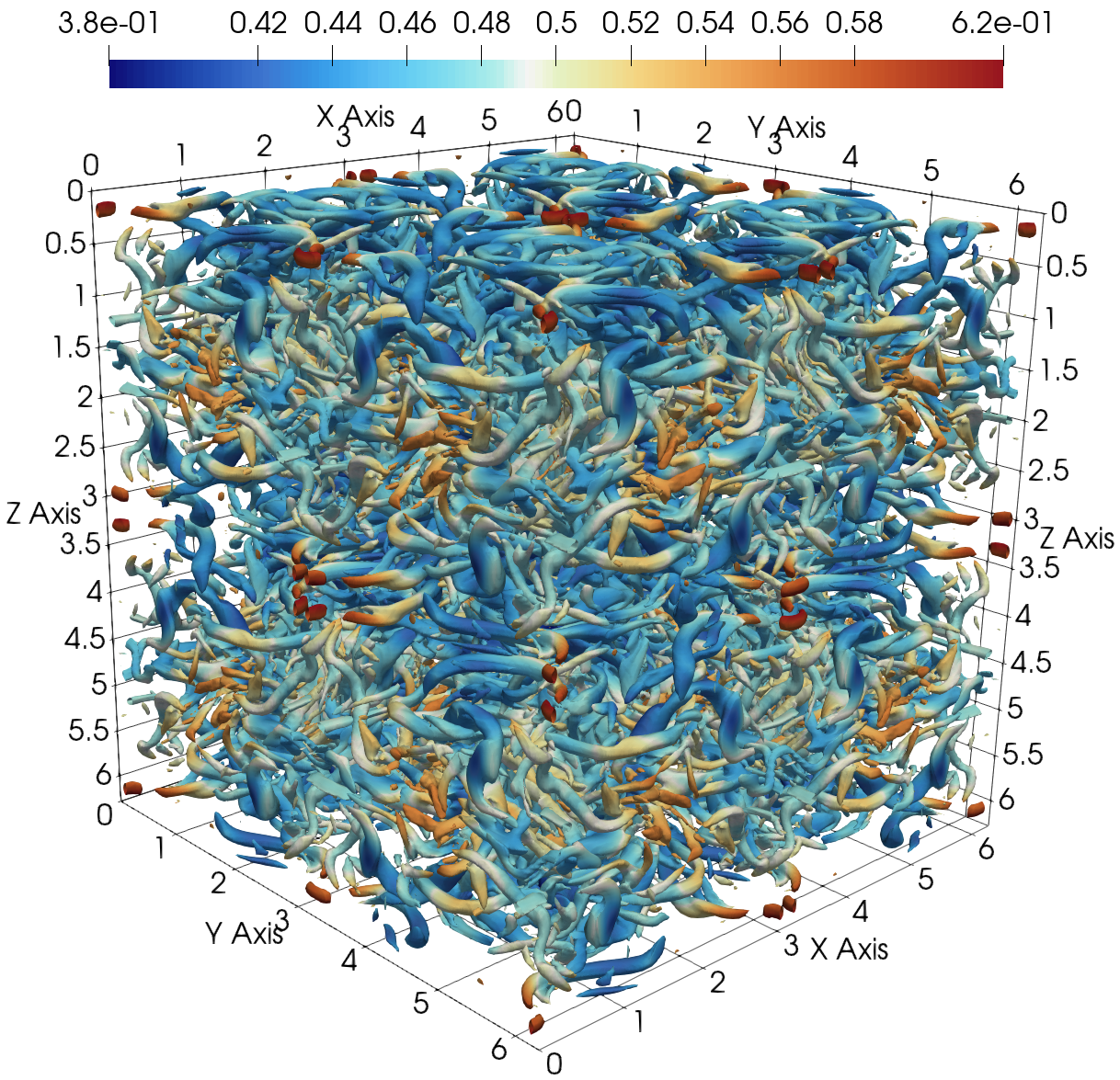}
    \caption{}
\end{subfigure}
\hfill
\begin{subfigure}[t]{0.49\textwidth}
    \centering
    \includegraphics[width=\linewidth]{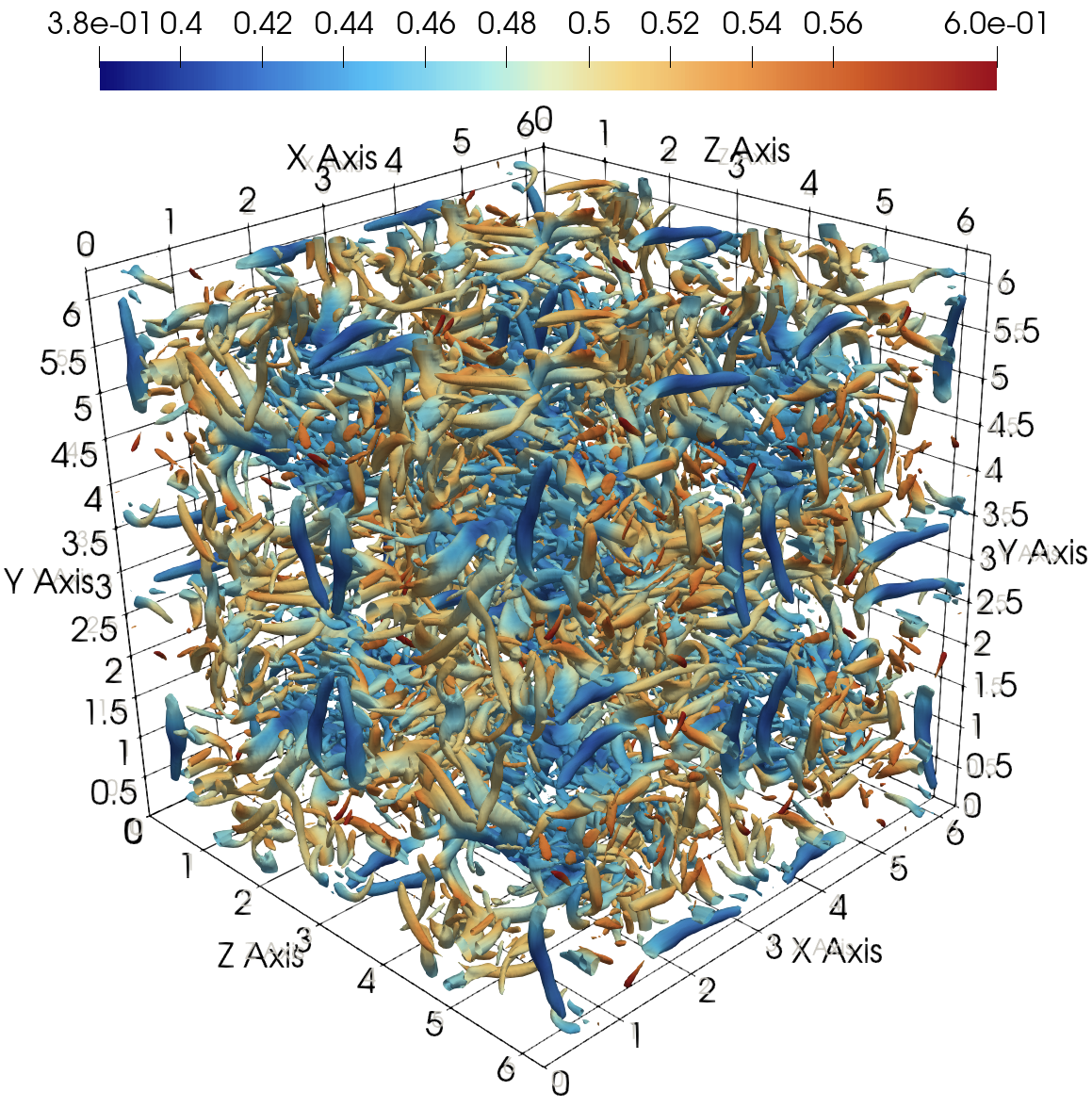}
    \caption{}
\end{subfigure}
\caption{Q-criterion iso-surfaces for the supersonic Taylor--Green vortex obtained on the $128^3$-element mesh at representative times: (a) $t=2$, (b) $t=6$, (c) $t=16$, and (d) $t=20$. The iso-surfaces are colored by pressure to highlight the coupling between vortical structures and the compressible field.}
\label{fig:tgv-qcriterion}
\end{figure*}

\begin{figure*}[htbp]
\centering
\begin{subfigure}[t]{0.49\textwidth}
    \centering
    \includegraphics[width=\linewidth]{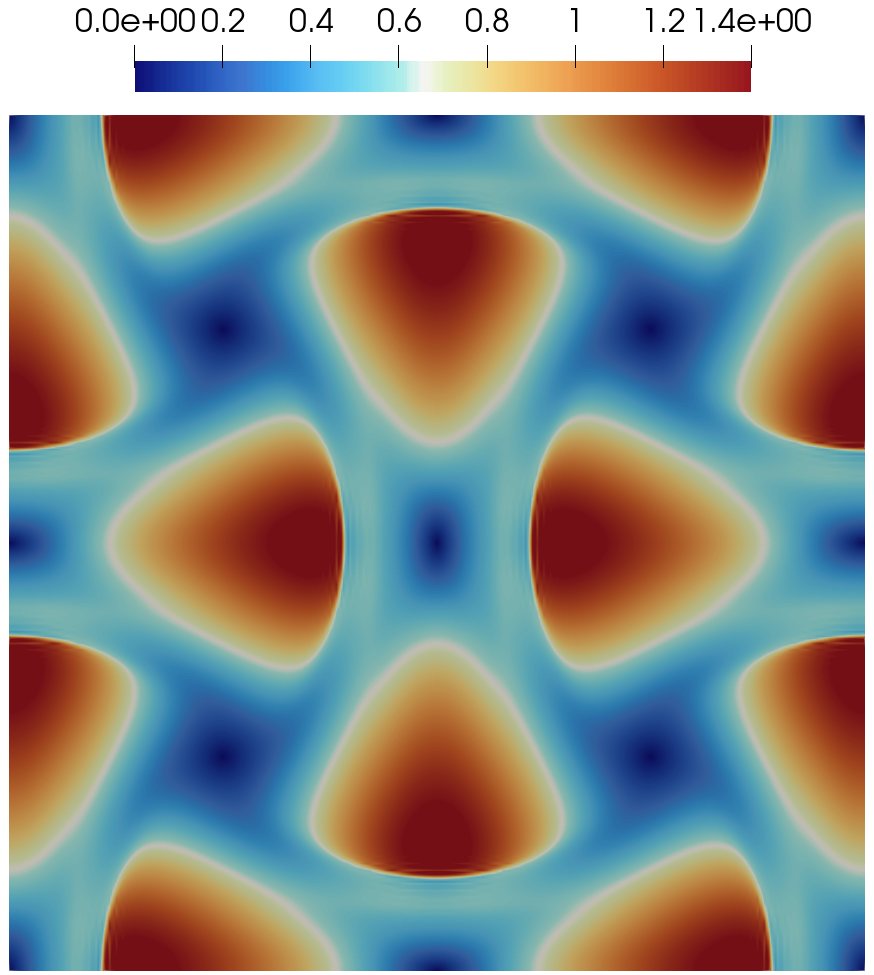}
    \caption{}
\end{subfigure}
\hfill
\begin{subfigure}[t]{0.49\textwidth}
    \centering
    \includegraphics[width=\linewidth]{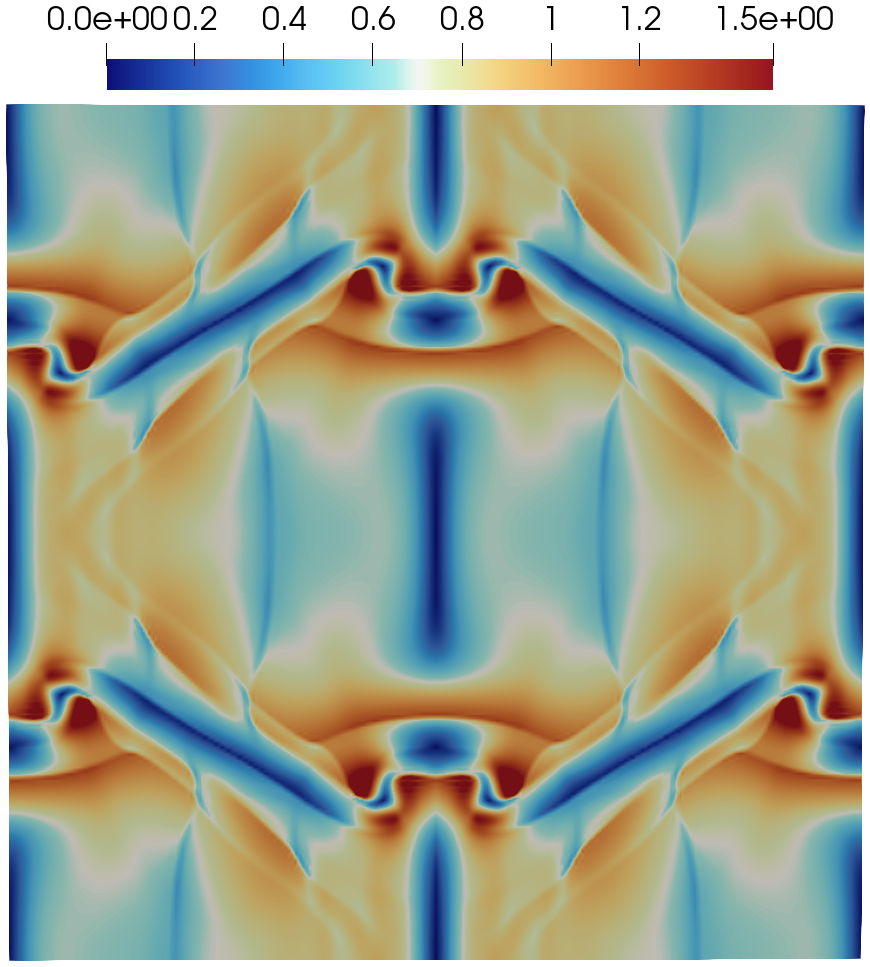}
    \caption{}
\end{subfigure}

\medskip

\begin{subfigure}[t]{0.49\textwidth}
    \centering
    \includegraphics[width=\linewidth]{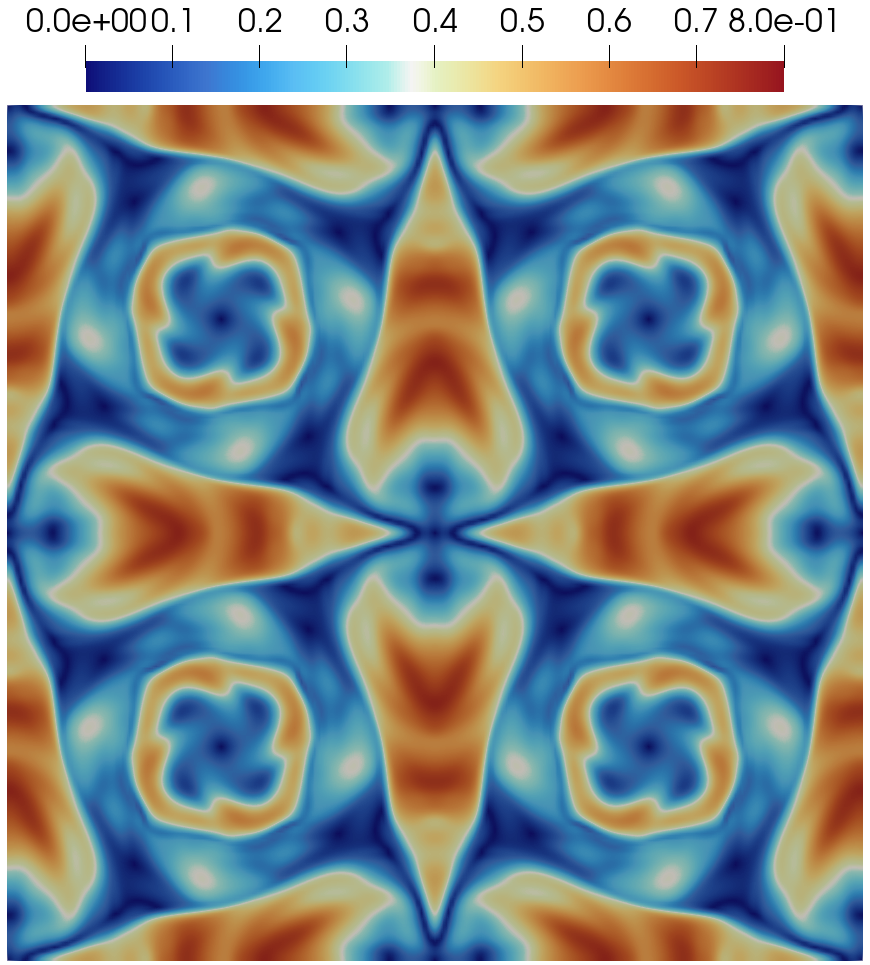}
    \caption{}
\end{subfigure}
\hfill
\begin{subfigure}[t]{0.49\textwidth}
    \centering
    \includegraphics[width=\linewidth]{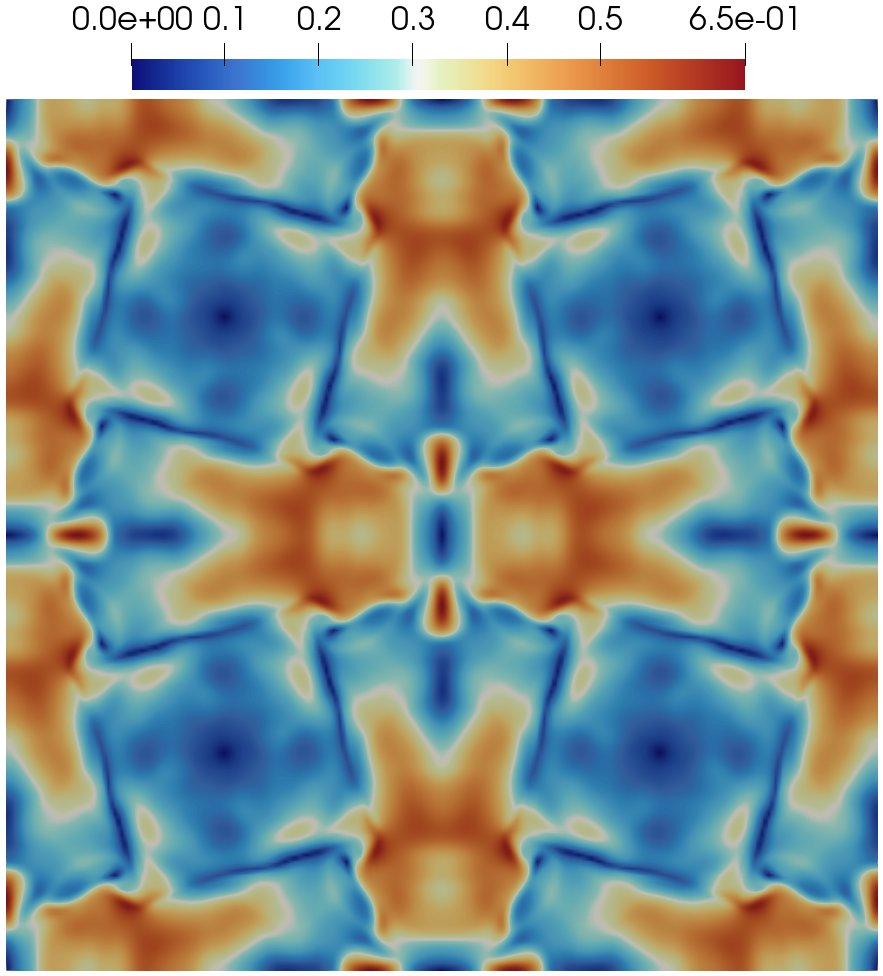}
    \caption{}
\end{subfigure}
\caption{Local Mach-number field in the $(x,y)$-plane for the supersonic Taylor--Green vortex obtained on the $128^3$-element mesh at (a) $t=2$, (b) $t=6$, (c) $t=16$, and (d) $t=20$. The slices show the progressive deformation of the large-scale compressible field, the emergence of steep compressive structures, and the subsequent weakening of the local Mach-number fluctuations during turbulent decay.}
\label{fig:tgv-mach-slices}
\end{figure*}

\begin{figure*}[htbp]
\centering
\begin{subfigure}[t]{0.49\textwidth}
    \centering
    \includegraphics[width=\linewidth]{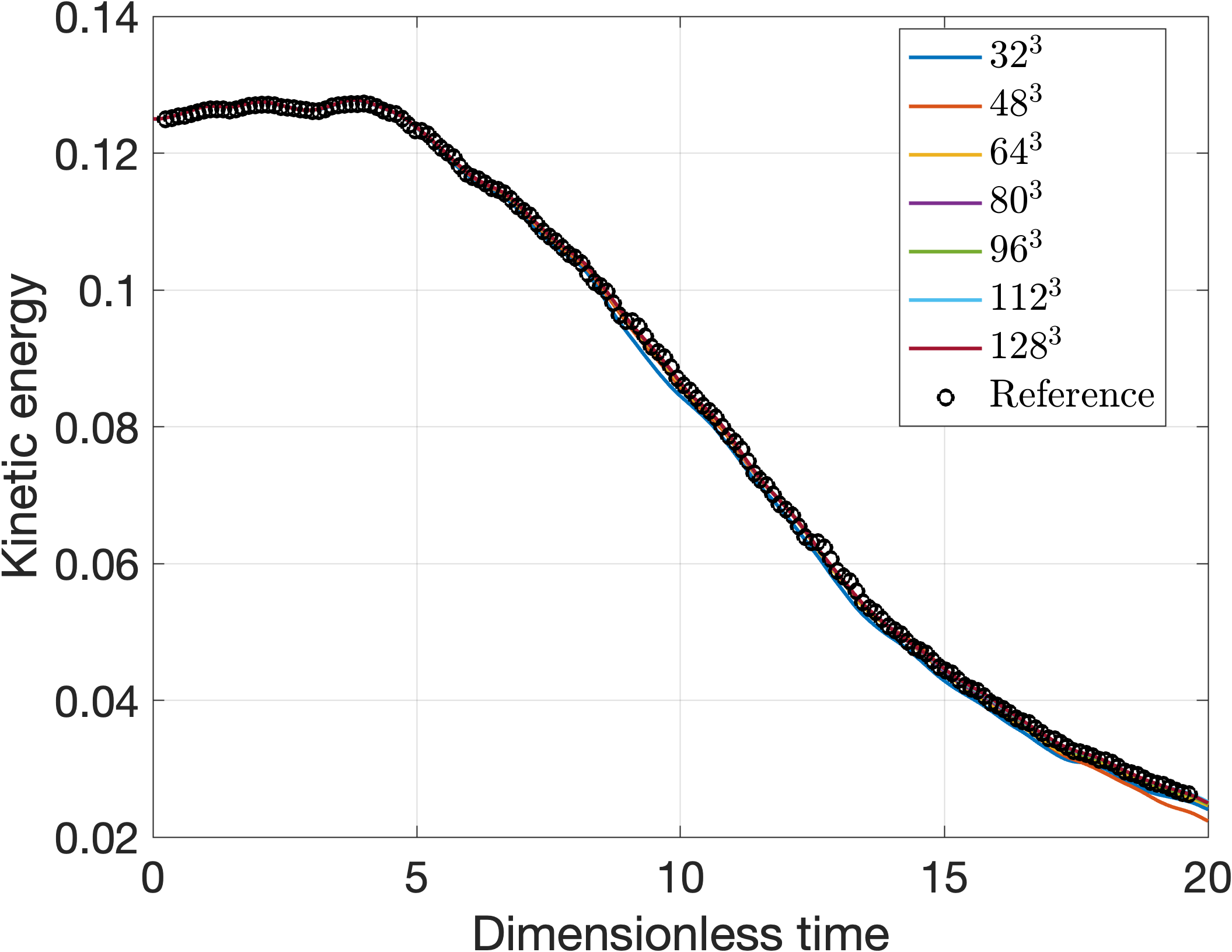}
    \caption{}
\end{subfigure}
\hfill
\begin{subfigure}[t]{0.49\textwidth}
    \centering
    \includegraphics[width=\linewidth]{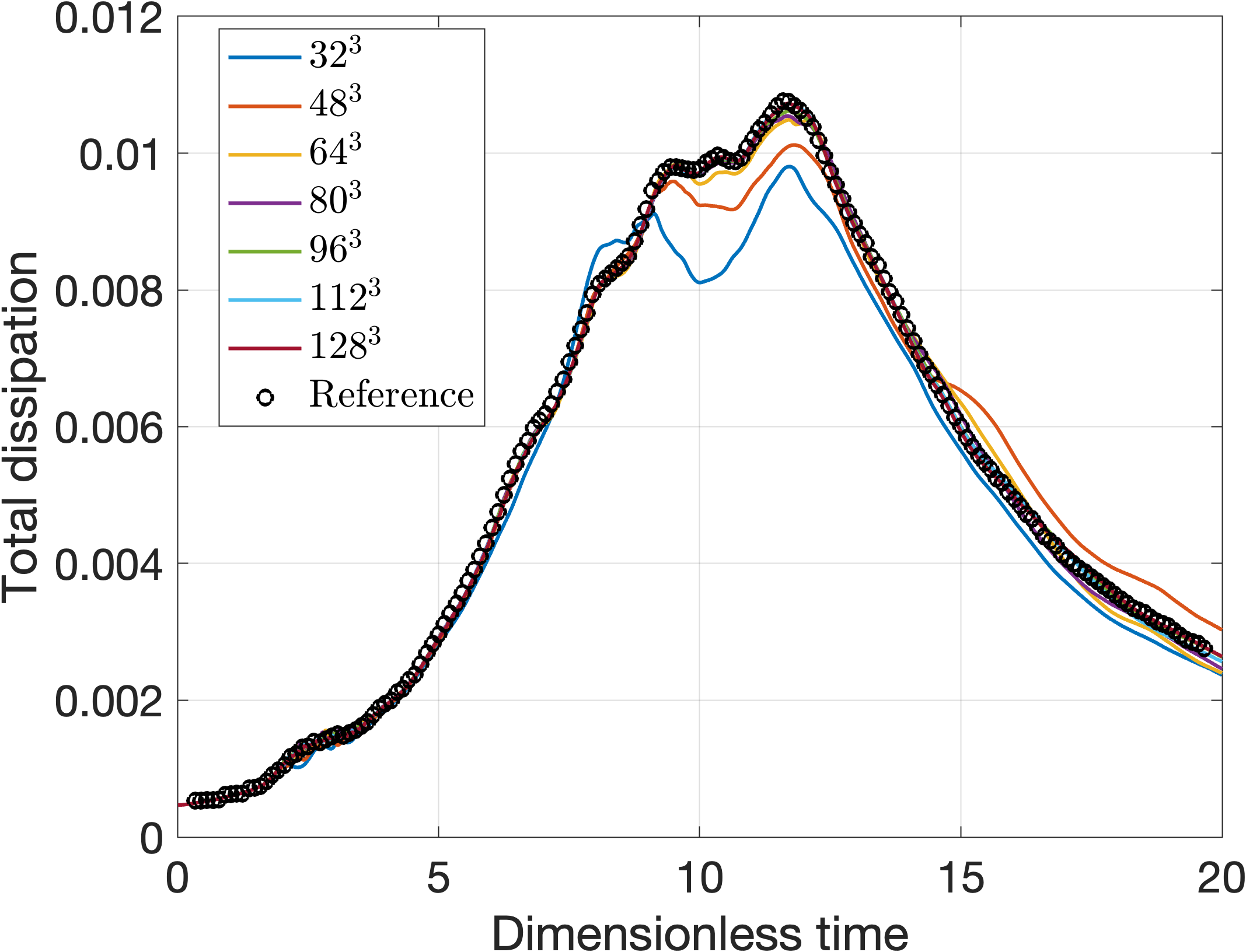}
    \caption{}
\end{subfigure}
\caption{Time histories of the integral diagnostics for the supersonic Taylor--Green vortex obtained on a sequence of successively refined meshes: (a) kinetic energy and (b) total viscous dissipation. The mesh sequence is used to assess spatial convergence, while the open circles denote reference DNS data from \cite{Lusher2020,Chapelier2024}.}
\label{fig:tgv-energy-total-dissipation}
\end{figure*}

\begin{figure*}[htbp]
\centering
\begin{subfigure}[t]{0.49\textwidth}
    \centering
    \includegraphics[width=\linewidth]{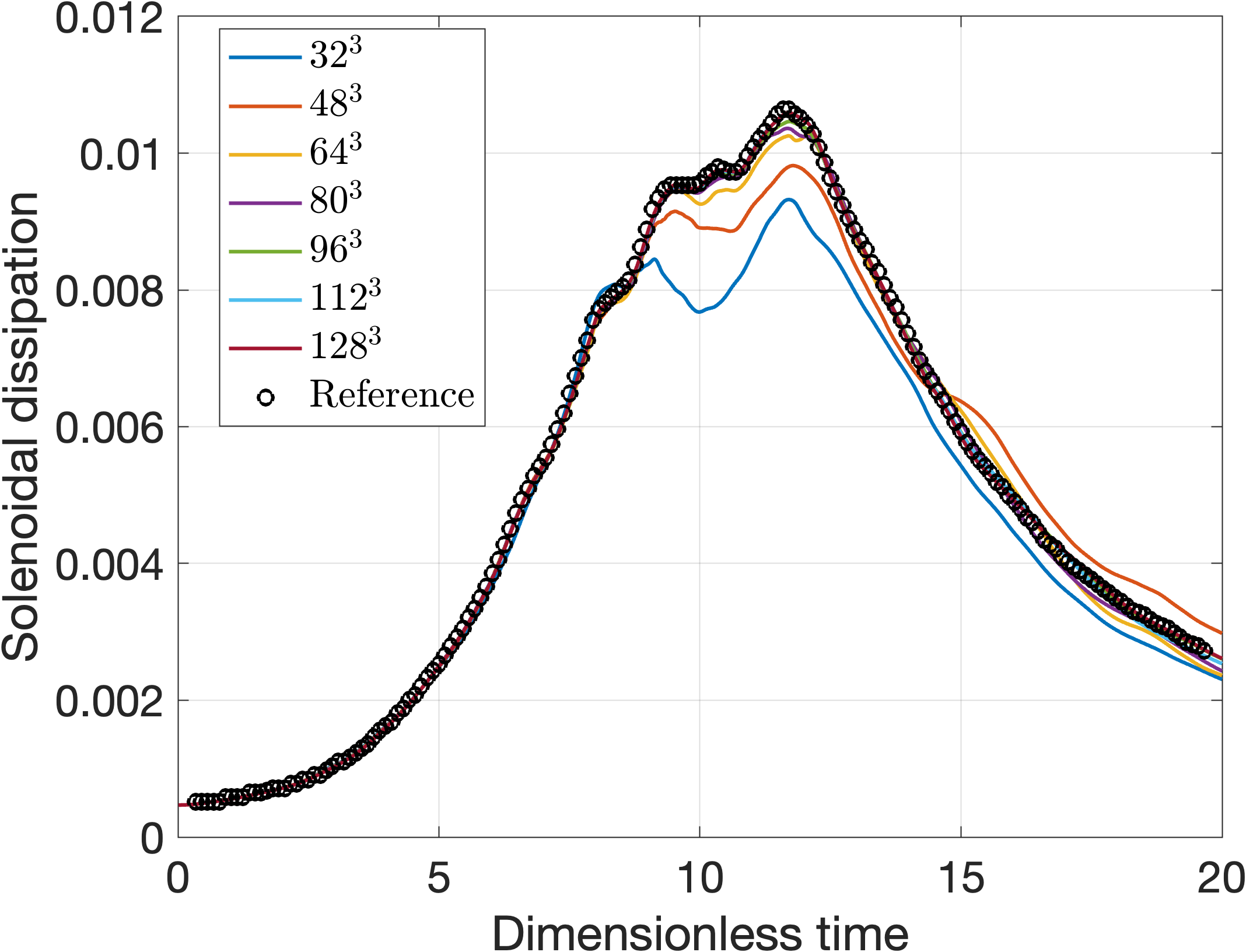}
    \caption{}
\end{subfigure}
\hfill
\begin{subfigure}[t]{0.46\textwidth}
    \centering
    \includegraphics[width=\linewidth]{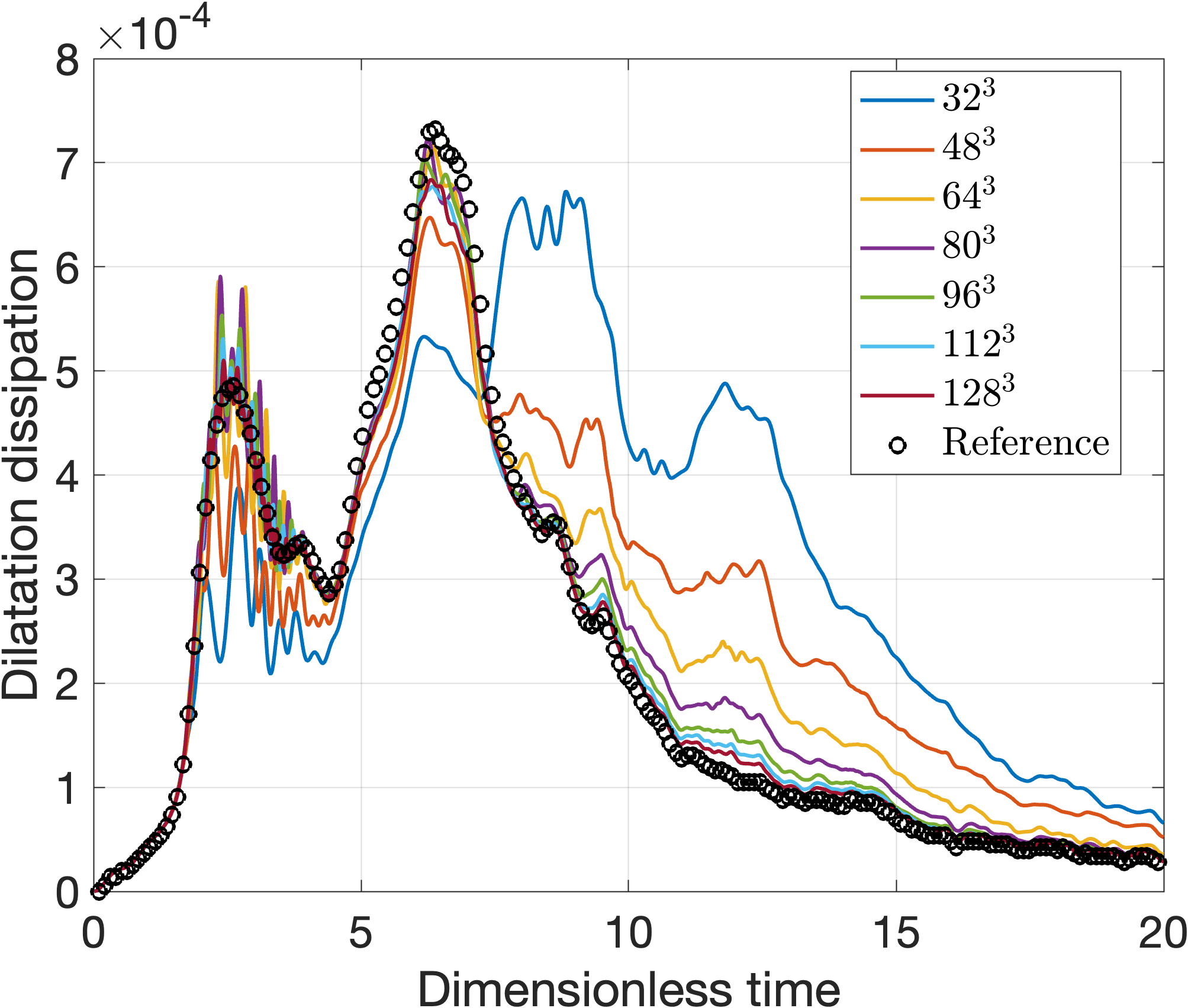}
    \caption{}
\end{subfigure}
\caption{Decomposition of the viscous dissipation rate for the supersonic Taylor--Green vortex obtained on a sequence of successively refined meshes: (a) solenoidal dissipation and (b) dilatational dissipation. The mesh sequence is used to assess spatial convergence, while the open circles denote reference DNS data from \cite{Lusher2020,Chapelier2024}.}
\label{fig:tgv-dissipation-components}
\end{figure*}

The supersonic Taylor--Green vortex provides a compact setting in which the transition from an initially smooth, highly organized vortex field to decaying compressible turbulence can be examined through complementary instantaneous and integral diagnostics. In the present computations, the instantaneous Q-criterion iso-surfaces in Fig.~\ref{fig:tgv-qcriterion} and the local Mach-number slices in Fig.~\ref{fig:tgv-mach-slices} are shown from the finest $128^3$-element simulation in order to visualize the vortical and compressive structure of the flow at representative times. The temporal histories in Figs.~\ref{fig:tgv-energy-total-dissipation} and \ref{fig:tgv-dissipation-components}, by contrast, are reported on the full mesh sequence $32^3$--$128^3$ defined in the preceding subsection and are used to assess mesh convergence against the reference data of \cite{Lusher2020,Chapelier2024}. Together, these diagnostics quantify the breakdown of the large-scale vortex system and the subsequent turbulent decay.

On the finest $128^3$-element mesh, the Q-criterion iso-surfaces show that the vortical field at the early time $t=2$ remains strongly constrained by the symmetries of the initial condition. The structures are sparse and highly organized, and the pressure coloring already indicates a non-negligible coupling between the vortex field and the compressible response. By $t=6$, these coherent structures have stretched and rolled into a substantially more intricate three-dimensional network. The iso-surfaces now trace elongated vortex tubes and sheets spanning a large fraction of the domain, which is consistent with the rapid growth of velocity gradients and the onset of strong nonlinear interaction among the primary vortices. At $t=16$, the original large-scale organization has been largely destroyed and replaced by a dense population of smaller-scale vortical structures filling the interior of the cube. By the final time $t=20$, the flow remains fully three-dimensional, but the structures are less coherent and less densely packed than at $t=16$, reflecting the progressive decay of the turbulent field. Across this sequence, the narrowing range of the pressure coloring is consistent with the attenuation of the strongest compressive events as the flow moves away from the most active transition stage.

The Mach-number slices provide the complementary view of the compressive dynamics. At $t=2$, the field is smooth and dominated by the large-scale imprint of the initial condition, with broad regions of locally supersonic and subsonic motion separated by relatively mild gradients. By $t=6$, this smooth structure has been replaced by sharply deformed bands and localized regions of intense compression and expansion. These steep features coincide with the stage at which the Q-criterion visualization shows strong stretching and interaction of the primary vortices, indicating that the evolving vortex system is actively generating compressive motions and shocklet-like structures. At later times the Mach-number field becomes increasingly multiscale, but its amplitude decreases. In particular, the color scale drops from values above unity at $t=2$ and $t=6$ to maxima below unity at $t=16$ and $t=20$, showing that the strongest supersonic pockets are transient and that the decay of turbulence is accompanied by a substantial weakening of the compressive field. Taken together, Figs.~\ref{fig:tgv-qcriterion} and \ref{fig:tgv-mach-slices} show that the most intense compressibility effects occur during the breakdown of the original vortex system rather than during the late-time decay.

The integral quantities quantify this evolution through the mesh sequence $32^3$, $48^3$, $64^3$, $80^3$, $96^3$, $112^3$, and $128^3$, and thereby provide a direct convergence study against the benchmark data reported by \cite{Lusher2020,Chapelier2024}. Figure~\ref{fig:tgv-energy-total-dissipation}(a) shows that the kinetic-energy histories converge monotonically under refinement: the coarse meshes begin to lose energy too rapidly once the flow becomes strongly multiscale, whereas the $96^3$, $112^3$, and $128^3$ solutions lie much closer to one another and to the reference curves of \cite{Lusher2020,Chapelier2024} over the full time interval. The energy is nearly constant during the earliest stage of the simulation, with a slight increase from approximately $0.125$ to about $0.127$ before the large-scale vortex system begins its sustained decay, and it subsequently falls to approximately $2.6\times 10^{-2}$ by $t=20$. The total dissipation in Fig.~\ref{fig:tgv-energy-total-dissipation}(b) exhibits a stronger mesh dependence than the kinetic energy. As the vortical field breaks down, the refined solutions approach a common broad maximum slightly above $10^{-2}$ around $t\approx 11$--$12$, in close agreement with the reference data of \cite{Lusher2020,Chapelier2024}, whereas the $32^3$ and $48^3$ meshes underestimate the peak and produce a broader late-time tail. The largest discrepancies therefore occur during the period of maximum dissipation, which coincides with the most intricate vortical structures in Fig.~\ref{fig:tgv-qcriterion} and is consistent with the transfer of energy toward smaller scales and the intensification of viscous dissipation during the breakdown process.

The dissipation decomposition in Fig.~\ref{fig:tgv-dissipation-components}, again shown on the full mesh sequence and compared with the benchmark data of \cite{Lusher2020,Chapelier2024}, clarifies the role of compressibility and the resolution required to capture it. The solenoidal contribution in Fig.~\ref{fig:tgv-dissipation-components}(a) dominates the total dissipation over most of the simulation and closely tracks the growth and decay of the vortical cascade. Its peak occurs at essentially the same time as the peak in total dissipation, confirming that the main dissipation maximum is controlled primarily by the amplification of vortical gradients. Under mesh refinement, the solenoidal-dissipation curves approach one another steadily, and the $96^3$--$128^3$ solutions reproduce both the timing and the magnitude of the reference maximum much more closely than the coarser meshes. The dilatational contribution in Fig.~\ref{fig:tgv-dissipation-components}(b) is roughly an order of magnitude smaller, with peak values of approximately $7\times 10^{-4}$, but it is far from negligible. It exhibits pronounced early-time activity, including a distinct maximum near $t\approx 6$--$7$, which is consistent with the appearance of the steep compressive structures in Fig.~\ref{fig:tgv-mach-slices}(b). The dilatational curves also show stronger resolution sensitivity than the solenoidal contribution, especially after the early peak: coarse meshes overpredict the dissipation level over a substantial portion of the decay interval, whereas the $96^3$, $112^3$, and $128^3$ solutions remain much closer to the reference behavior of \cite{Lusher2020,Chapelier2024}. This difference is physically consistent with the sharper compressive features that contribute to dilatational dissipation and that therefore require finer spatial resolution than the broader vortical structures dominating the solenoidal term.

The Q-criterion and Mach-number fields from the finest simulation demonstrate the transition from organized large-scale vortices to a complex three-dimensional turbulent state containing transient shocklet-like structures, while the mesh-sequence study shows that the associated kinetic-energy and dissipation histories converge toward the benchmark results of \cite{Lusher2020,Chapelier2024}. The convergence is more rapid for kinetic energy than for the dissipation measures, and the dilatational component remains the most resolution-sensitive quantity, which indicates that capturing the coupled vortical and compressive dynamics of the supersonic Taylor--Green vortex requires finer spatial resolution than would be inferred from energy decay alone. Within that context, the HDG solver recovers the essential temporal and structural features of compressible turbulent decay in this benchmark problem.


\section{Hypersonic Boundary-Layer Transition}

\subsection{Problem Description}

\begin{figure}[t]
\centering
\includegraphics[alt={sample image},width=\textwidth]{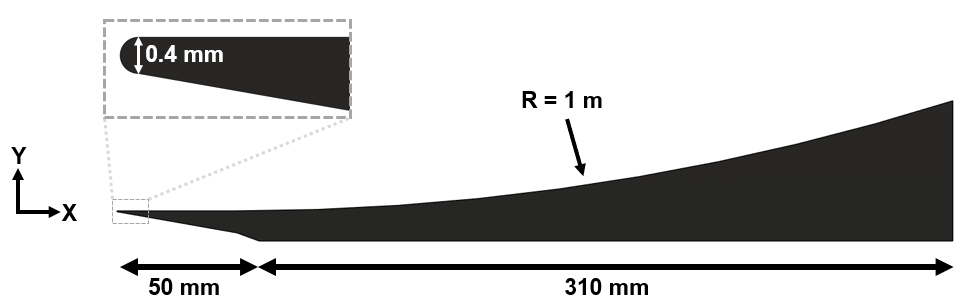}
\caption{Streamwise definition of the Mach-7 hypersonic boundary-layer-transition geometry. The figure shows the rounded leading edge, the $50\,\mathrm{mm}$ forebody, the $310\,\mathrm{mm}$ streamwise-concave flared segment with curvature radius $R=1\,\mathrm{m}$, and the $0.4\,\mathrm{mm}$ nose-tip diameter corresponding to a nose radius of $0.2\,\mathrm{mm}$.}
\label{Streamwisegeo}
\end{figure}

Hypersonic boundary-layer transition provides a stringent  benchmark because the relevant dynamics involve strong compressibility, very thin near-wall layers, steep wall-normal gradients, and the coexistence of multiple instability mechanisms over widely separated spatial and temporal scales. In high-Mach-number boundary layers, second-mode, or Mack-mode, disturbances are a canonical acoustic instability mechanism \cite{Kendall1975,Zhang2013}. Curved hypersonic walls introduce an additional route to transition: streamwise concave curvature produces centrifugal destabilization, leading to the growth of streamwise-oriented G{\"o}rtler vortices, streaks, secondary instability, and eventual breakdown \cite{Balakumar2015,Chynoweth2019,Hader2019}. The present Mach-7 configuration is selected to expose these instability mechanisms. The geometry consists of a smooth wall with fixed streamwise concavity and controlled spanwise convex curvature. The streamwise profile, shown in Fig.~\ref{Streamwisegeo}, has a total length of $360\,\mathrm{mm}$. A $50\,\mathrm{mm}$ forebody follows the rounded leading edge and is followed by a $310\,\mathrm{mm}$ flared segment with streamwise radius of curvature $R=1\,\mathrm{m}$. The nose radius is $2\times10^{-4}\,\mathrm{m}$, equivalently $0.2\,\mathrm{mm}$, corresponding to the $0.4\,\mathrm{mm}$ nose-tip diameter shown in the inset. This streamwise concavity provides the centrifugal forcing associated with the primary G{\"o}rtler instability. The spanwise extrusion uses a lateral computational extent of $30\,\mathrm{mm}$, with the homogeneous lateral direction treated periodically. The freestream and wall conditions are summarized in Table~\ref{tab:simulation_parameters}. The flow is computed at $M_\infty=7$ and $\alpha=0^\circ$ with $\gamma=1.4$ and $Pr=0.71$. The wall is maintained isothermal at $T_{\mathrm{wall}}=296\,\mathrm{K}$, while freestream temperature is $T_{\mathrm{ref}}=51.86\,\mathrm{K}$.

\begin{table}[htbp]
  \centering
  \caption{Flow, thermal, and nondimensional reference parameters for the Mach-7 hypersonic boundary-layer-transition simulations. Dimensional quantities are given in SI units; the nondimensional freestream state corresponds to the solver normalization used in Section~6.4.}
  \label{tab:simulation_parameters}
  \begin{tabular}{@{}p{0.25\textwidth}p{0.49\textwidth}p{0.2\textwidth}@{}}
    \toprule
    \textbf{Category} & \textbf{Parameter} & \textbf{Value} \\
    \midrule
    Fluid properties
      & Specific heat ratio ($\gamma$) & 1.4 \\
      & Prandtl number ($Pr$) & 0.71 \\
    \midrule
    Flow conditions
      & Freestream Mach number ($M_\infty$) & 7.0 \\
      & Reynolds per meter ($Re$) & $8.088\times 10^6$ \\
      & Angle of attack ($\alpha$) & $0^\circ$ \\
    \midrule
    Thermal state
      & Freestream temperature ($T_{\mathrm{ref}}$) & $51.86\,\mathrm{K}$ \\
      & Isothermal wall temperature ($T_{\mathrm{wall}}$) & $296.0\,\mathrm{K}$ \\
    \midrule
    Dimensionless
      & Freestream density ($\rho_\infty$) & 1.0 \\
      & Freestream pressure ($p_\infty$) & 0.01446 \\
      & Freestream total energy ($E_\infty$) & 0.53614 \\
    \bottomrule
  \end{tabular}
\end{table}

Relative to the previous subsonic airfoil and supersonic Taylor--Green benchmarks, the present configuration adds thin hypersonic boundary layers, strong compressibility, curved three-dimensional geometry, curvature-driven G{\"o}rtler structures, high-frequency acoustic disturbances, and nonlinear secondary breakdown. The problem therefore tests whether the proposed HDG solver can capture instability growth and transition in a strongly compressible wall-bounded flow without resorting to transition or turbulence models.

\subsection{Simulation Setup}

\begin{figure}[t]
    \centering
    \includegraphics[width=0.98\linewidth]{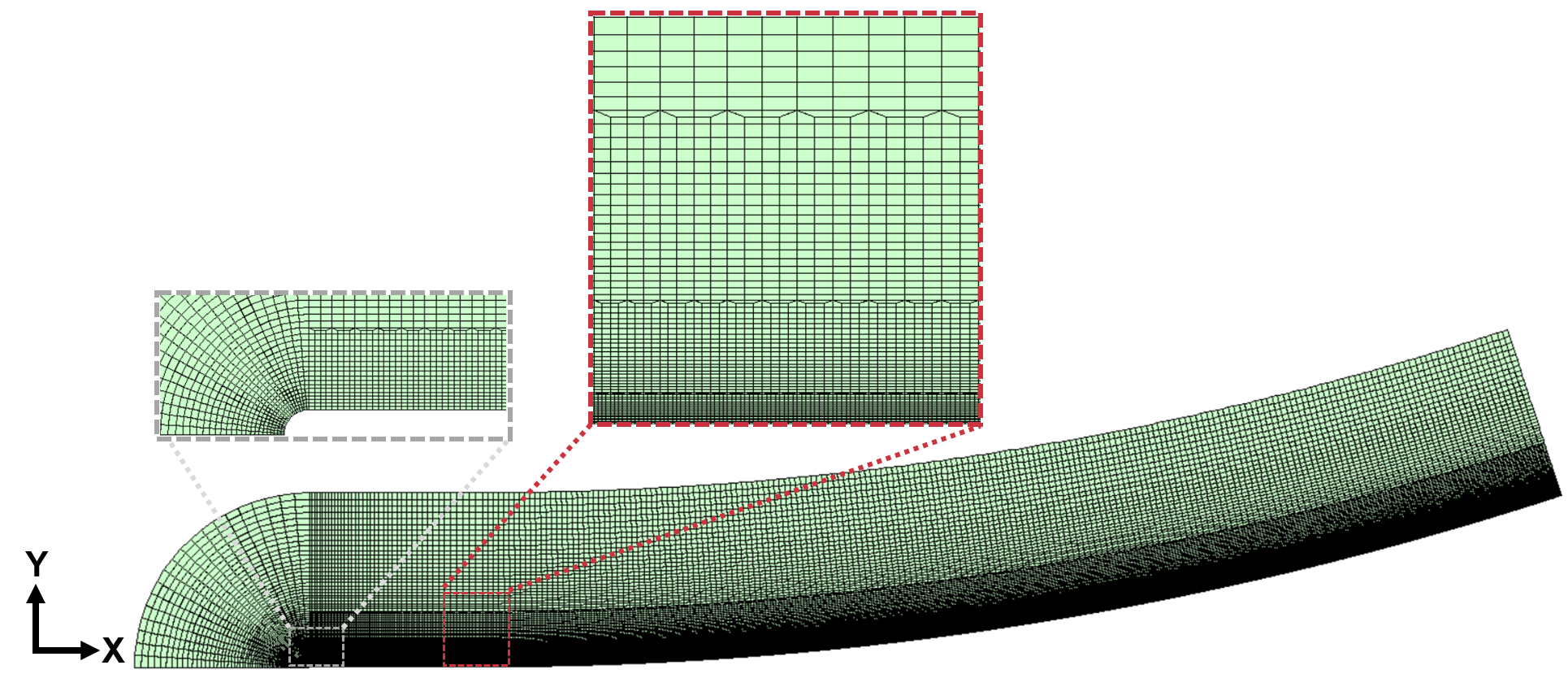}
    \caption{Representative side view of the high-order mesh used for the Mach-7 hypersonic boundary-layer-transition simulations. The mesh is clustered near the rounded leading edge and along the isothermal wall to resolve the hypersonic boundary layer and the streamwise-concave flare region, while the outer elements are smoothly stretched toward the freestream boundary.}
    \label{fig:meshside}
\end{figure}

The hypersonic calculations use the same implicit HDG methodology described in Sections~2 and 3, but on a much more fine wall-bounded mesh than the preceding benchmark problems. The computational grid is a body-fitted, high-order discretization of the flared geometry described above. It is strongly clustered near the leading edge and along the isothermal wall to resolve the thin hypersonic boundary layer, the near-wall thermal layer, and the early growth of streamwise G{\"o}rtler streaks. The outer part of the mesh expands smoothly away from the wall, keeping the freestream boundary well separated from the boundary-layer region while avoiding unnecessary far-field refinement. A representative streamwise side view of the mesh is shown in Fig.~\ref{fig:meshside}. The three-dimensional production meshes contain more than $12$ million quadratic elements; the same polynomial degree, $k=2$, is used for the curved geometric mapping and for the HDG solution approximation. Time advancement is performed with the three-stage, third-order DIRK scheme used throughout the paper. The time step is fixed at $\Delta t = 2\times 10^{-4}$, and each case is advanced for $9200$ time steps, corresponding to a final nondimensional time of $t=1.84$.

The unsteady simulations are initialized from a converged 2D laminar precursor solution on the same streamwise geometry. The precursor state is mapped onto the 3D mesh and extended in the homogeneous lateral direction, with the spanwise velocity initially set by the mean precursor flow. This initialization removes the startup transient associated with forming the basic hypersonic boundary layer and allows the subsequent calculation to focus on the growth of the imposed disturbances, centrifugal streaks, and nonlinear breakdown. The artificial-viscosity field for shock and compression regularization is initialized consistently with the precursor solution and then evolved through the stage-lagged treatment described in Section~2.

Boundary data are imposed weakly through the HDG numerical fluxes. The inflow boundary prescribes the Mach-7 mean state together with the synthetic velocity perturbations described in the next subsection. The wall is no-slip and isothermal at the temperature given in Table~\ref{tab:simulation_parameters}. A symmetry condition is used on the upstream portion of the boundary preceding the physical wall, the downstream boundary is treated as a supersonic outflow, and the outer boundary uses the freestream condition. The homogeneous lateral direction is periodic. The simulations reported below use the same mesh, wall condition, freestream state, and time integration while varying the prescribed inflow turbulence intensity over $I=0.5\%$, $0.2\%$, and $0.1\%$; this isolates the effect of disturbance amplitude on the transition process.

All hypersonic simulations were performed with the GPU implementation of the HDG solver. A representative production run required approximately $6{,}600$ node-hours and used $1{,}536$ AMD Instinct MI250X GPUs on Frontier at the Oak Ridge Leadership Computing Facility; corresponding runs were also carried out on Tuolumne using $1{,}024$ AMD Instinct MI300A GPUs. These calculations constitute the largest wall-bounded DNS demonstrations in the present study and provide a stringent test of the distributed trace-space solver on highly anisotropic hypersonic meshes.

\subsection{Synthetic Inflow Disturbances}

Synthetic velocity disturbances are imposed at the inlet by adding a random Fourier field to the mean inflow profile. Following the stochastic construction of Kraichnan~\cite{Kraichnan1970}, B\'{e}chara et al.~\cite{Bechara1994}, and Bailly and Juv\'{e}~\cite{Bailly1999}, the perturbation velocity is written as
\begin{equation}
\mathbf{u}'(\mathbf{x},t)
=
2\sum_{n=1}^{N}
\hat{u}_n \boldsymbol{\sigma}_n
\cos\!\left[
\mathbf{k}_n \cdot
\left(\mathbf{x}-\mathbf{u}_c t\right)
+\psi_n+\omega_n t
\right],
\label{eq:synthetic_turbulence}
\end{equation}
where $\hat{u}_n$, $\boldsymbol{\sigma}_n$, $\mathbf{k}_n$, $\psi_n$, and $\omega_n$ are the amplitude, direction vector, wave vector, phase, and angular frequency of mode $n$, respectively, and $\mathbf{u}_c=\mathbf{u}_{\infty}$. The imposed inlet velocities are therefore given by
\begin{equation}
\label{eq:inlet_velocity}
\mathbf{u}_{\mathrm{in}}(\mathbf{x},t)
=
{\mathbf{u}}_{\infty}
+
\mathbf{u}'(\mathbf{x},t),
\end{equation}
where ${\mathbf{u}}_{\infty}$ is the freestream velocity.

The wave-vector orientations are sampled uniformly on the sphere, and the modal amplitudes are set by the target spectrum,
\begin{equation}
\hat{u}_n
=
\sqrt{E(k_n)\,\Delta k_n},
\label{eq:mode_amplitude}
\end{equation}
where $\Delta k_n$ is the local spacing in wave-number space. The wave numbers are distributed logarithmically from
\begin{equation}
k_1=\frac{k_e}{5}
\qquad\text{and}\qquad
k_N=k_{\eta},
\end{equation}
with intermediate values defined by
\begin{equation}
k_n=rk_{n-1},
\qquad
r=
\left(\frac{k_N}{k_1}\right)^{1/(N-1)}.
\label{eq:wavenumber_spacing}
\end{equation}
The random variables $\boldsymbol{\sigma}_n$, $\mathbf{k}_n$, $\psi_n$, and $\omega_n$ follow the stochastic procedure of Bailly and Juv\'{e}~\cite{Bailly1999}. The prescribed von K\'{a}rm\'{a}n--Pao spectrum is
\begin{equation}
E(k)
=
\alpha
\frac{I^2 u_{\infty}^2}{k_e}
\frac{(k/k_e)^4}
     {\left[1+(k/k_e)^2\right]^{17/6}}
\exp\!\left[-2\left(\frac{k}{k_{\eta}}\right)^2\right],
\label{eq:von_karman_pao}
\end{equation}
where $I$ is the prescribed turbulence intensity and $u_{\infty}$ is the freestream velocity magnitude. The dissipative and energy-containing wave numbers are
\begin{equation}
\label{eq:kolmogorov_wavenumber}
k_{\eta}
=
\epsilon^{1/4}\nu^{-3/4},
\qquad k_e
=
\frac{9\pi\alpha}{55\Lambda}.
\end{equation}
where $\nu$ is the kinematic viscosity, $\epsilon$ is the turbulent kinetic energy dissipation rate, and $\Lambda$ is the integral length scale. The coefficient $\alpha=1.453$ normalizes the spectrum so that
\begin{equation}
\overline{k}
=
\int_{0}^{\infty} E(k)\,\mathrm{d}k
=
\frac{3}{2}\left(u_{\infty}I\right)^2.
\label{eq:turbulent_kinetic_energy}
\end{equation}
The integral length scale satisfies
\begin{equation}
\Lambda
=
\frac{\pi}{2\left(u_{\infty}I\right)^2}
\int_{0}^{\infty}
\frac{E(k)}{k}\,\mathrm{d}k.
\label{eq:integral_length_scale}
\end{equation}
In the present calculations, the integral length scale is tied to the near-wall resolution $\Lambda=10h_{\min}$,
 where $h_{\min}$ is the minimum wall-normal grid spacing, and the dissipation rate is estimated from
\begin{equation}
\epsilon
=
c_{\mu}^{3/4}
\frac{\overline{k}^{3/2}}{\Lambda},
\qquad
c_{\mu}=0.09.
\label{eq:dissipation_rate}
\end{equation}
The free-stream intensity $I$ controls the amplitude of the free-stream disturbance. We study the effects of this parameter on hypersonic boundary layer instabilities and the onset of laminar-to-turbulent transition.

\subsection{Boundary Layer Instability and Transition}

This problem involves the simultaneous evolution of freestream disturbances  and curvature-induced vortical structures. As discussed above, the streamwise-concave wall promotes centrifugal instability and the amplification of streamwise G{\"o}rtler vortices, while the cold hypersonic boundary layer can also support high-frequency acoustic disturbances, including Mack-mode activity \cite{Kendall1975,Zhang2013,Balakumar2015,Chynoweth2019,Hader2019}. The results presented in this subsection examine how this instability system responds to reductions in the imposed synthetic freestream turbulence intensity from $0.5\%$ to $0.2\%$ and $0.1\%$. The transition process is characterized using complementary evidence from instantaneous and time-averaged pressure and density fields, wall heat flux, skin friction, and spanwise/time-averaged wall coefficients.

\begin{figure}[t]
    \centering    
    \begin{subfigure}[b]{0.48\linewidth}
        \centering
        \includegraphics[width=\linewidth]{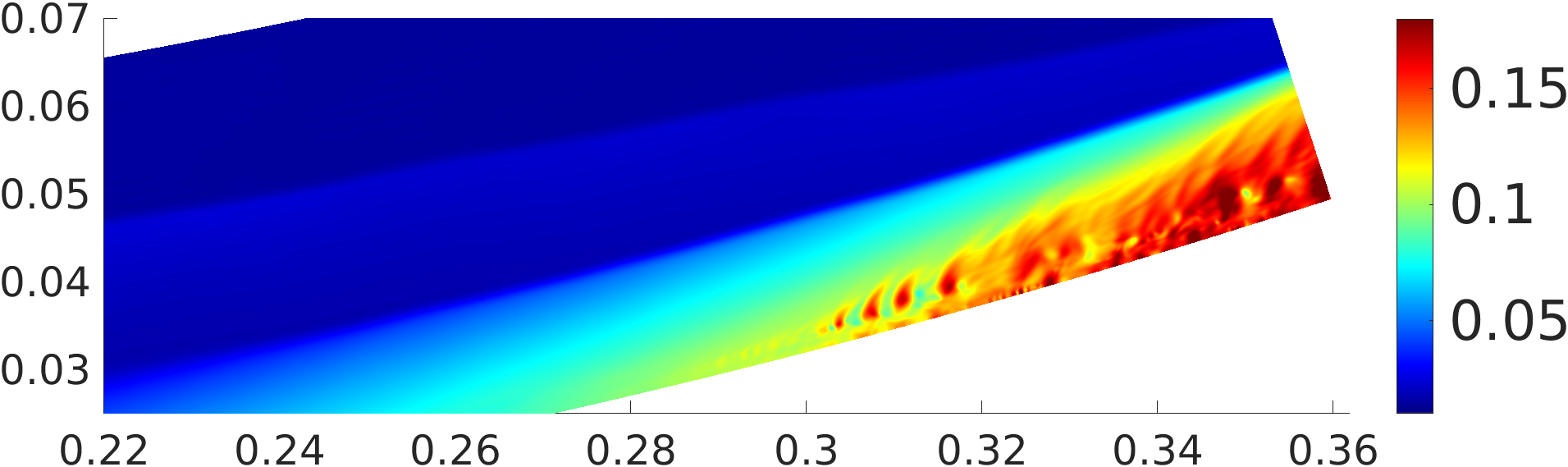}
        \caption{0.5\% TI}
    \end{subfigure}
    \hfill
    \begin{subfigure}[b]{0.48\linewidth}
        \centering
        \includegraphics[width=\linewidth]{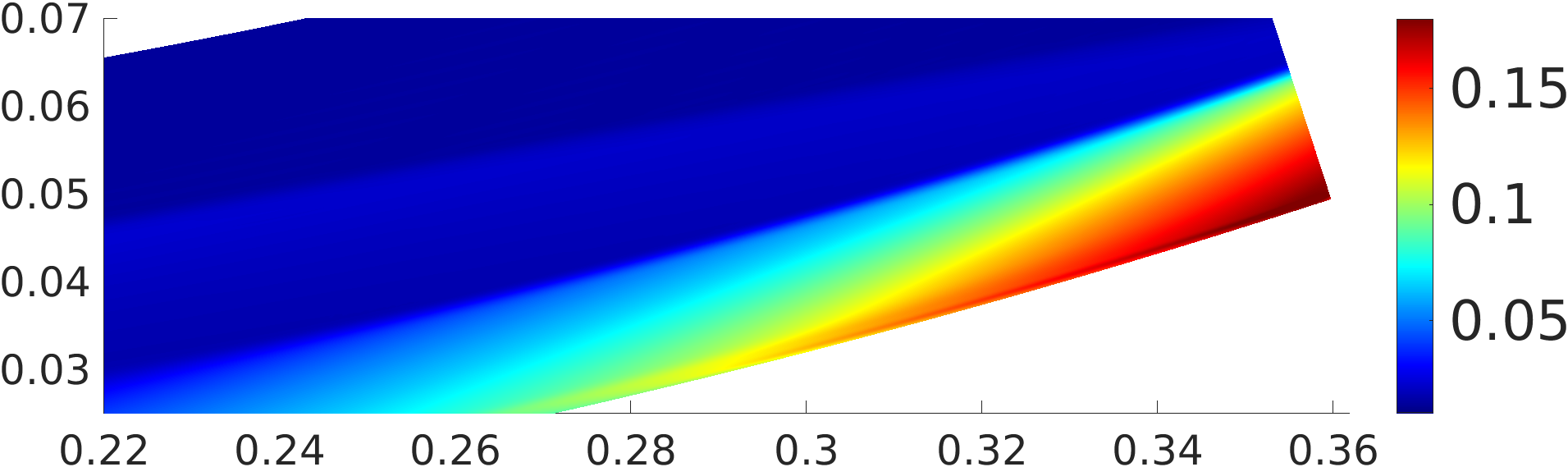}
        \caption{0.5\% TI}
    \end{subfigure}

    \begin{subfigure}[b]{0.48\linewidth}
        \centering
        \includegraphics[width=\linewidth]{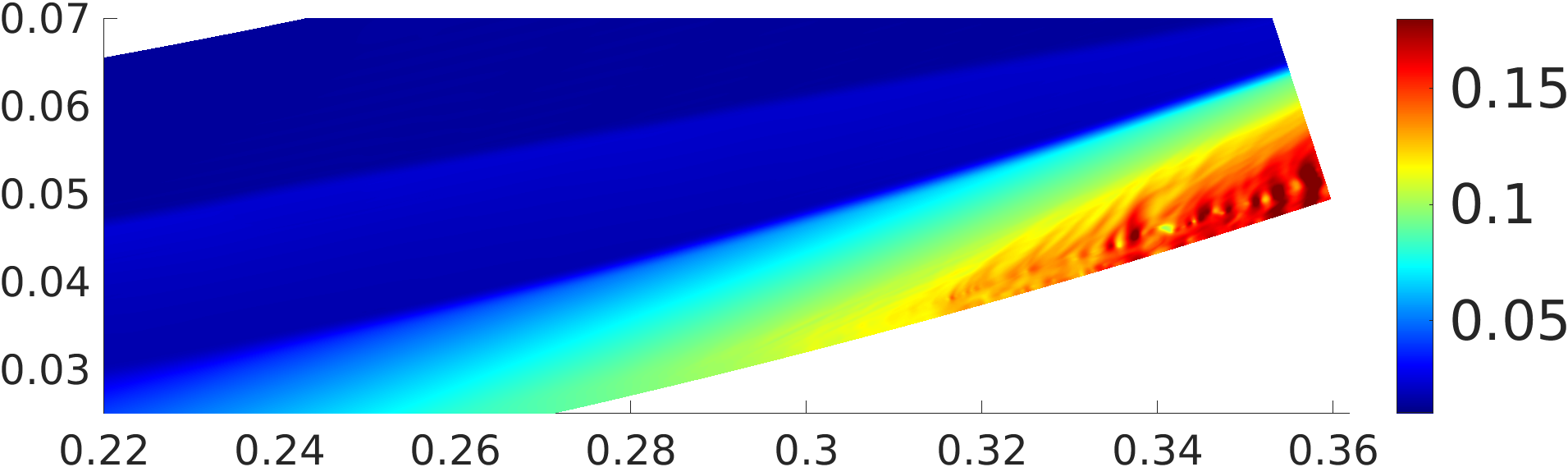}
        \caption{0.2\% TI}
    \end{subfigure}
    \hfill    
    \begin{subfigure}[b]{0.48\linewidth}
        \centering
        \includegraphics[width=\linewidth]{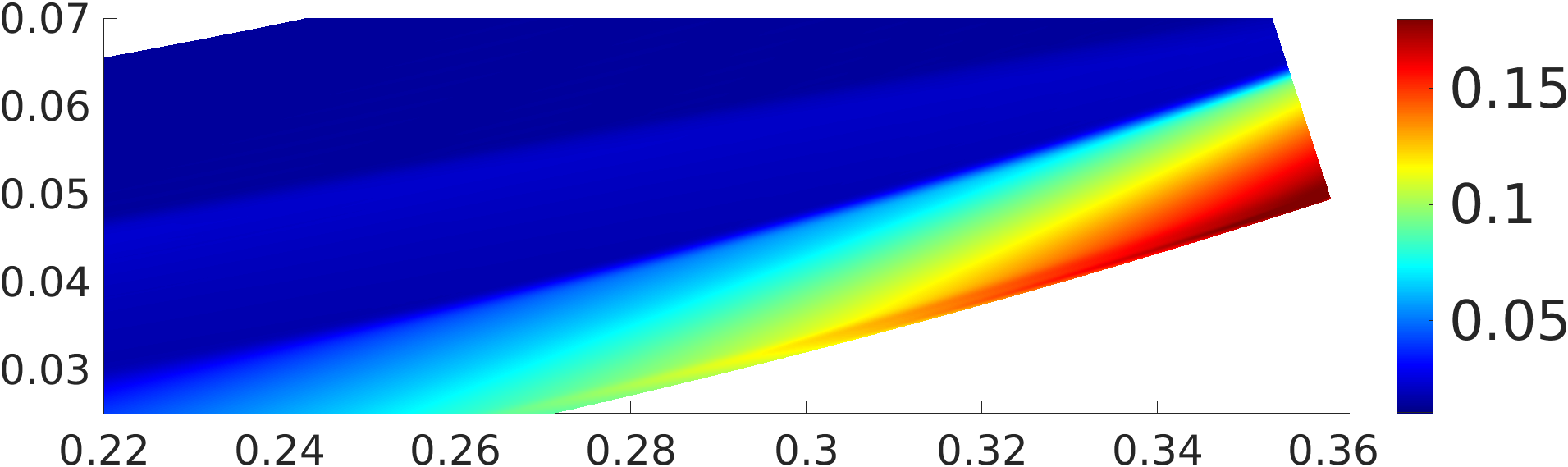}
        \caption{0.2\% TI}
    \end{subfigure}

    \begin{subfigure}[b]{0.48\linewidth}
        \centering
        \includegraphics[width=\linewidth]{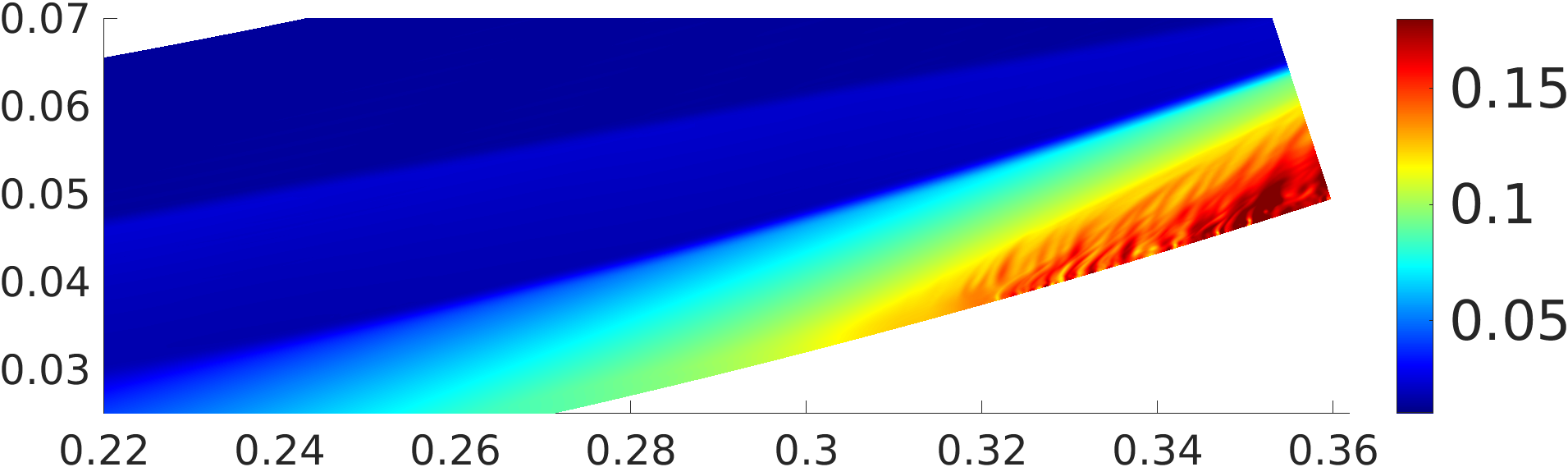}
        \caption{0.1\% TI}
    \end{subfigure}
    \hfill
    \begin{subfigure}[b]{0.48\linewidth}
        \centering
        \includegraphics[width=\linewidth]{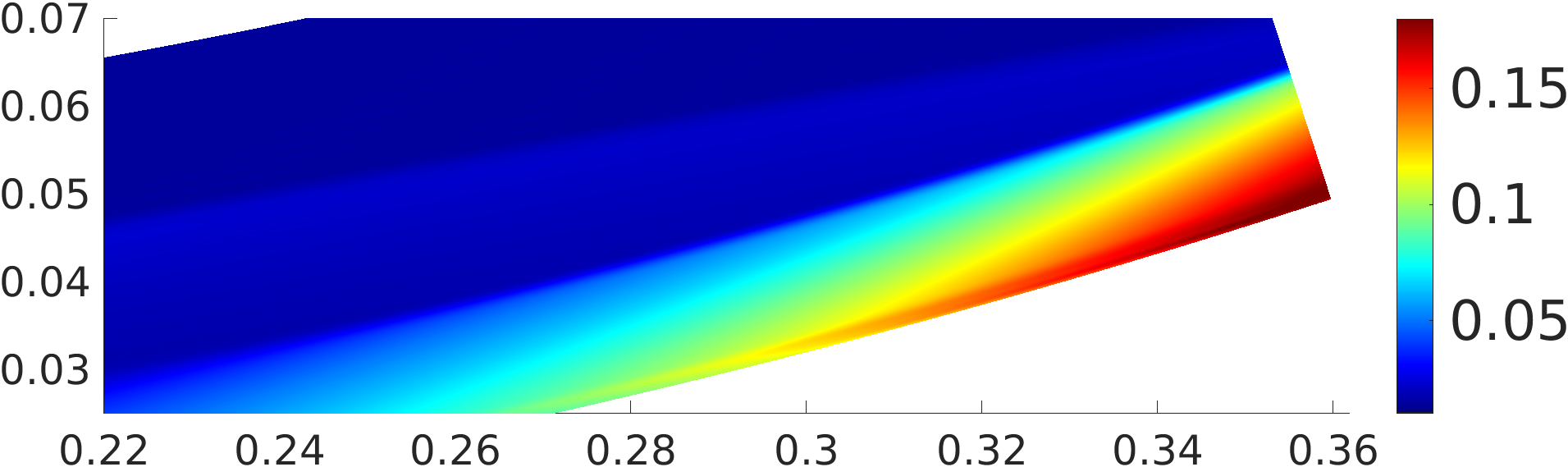}
        \caption{0.1\% TI}
    \end{subfigure}

    \caption{Dimensionless pressure at a representative side-view slice for
    freestream turbulence intensities of $0.5\%$, $0.2\%$, and $0.1\%$.
    Instantaneous fields are shown in the left column and the corresponding
    time-averaged fields in the right column.}
    \label{pressureside}
\end{figure}

\begin{figure}[h]
    \centering
    \begin{subfigure}[b]{0.48\linewidth}
        \centering
        \includegraphics[width=\linewidth]{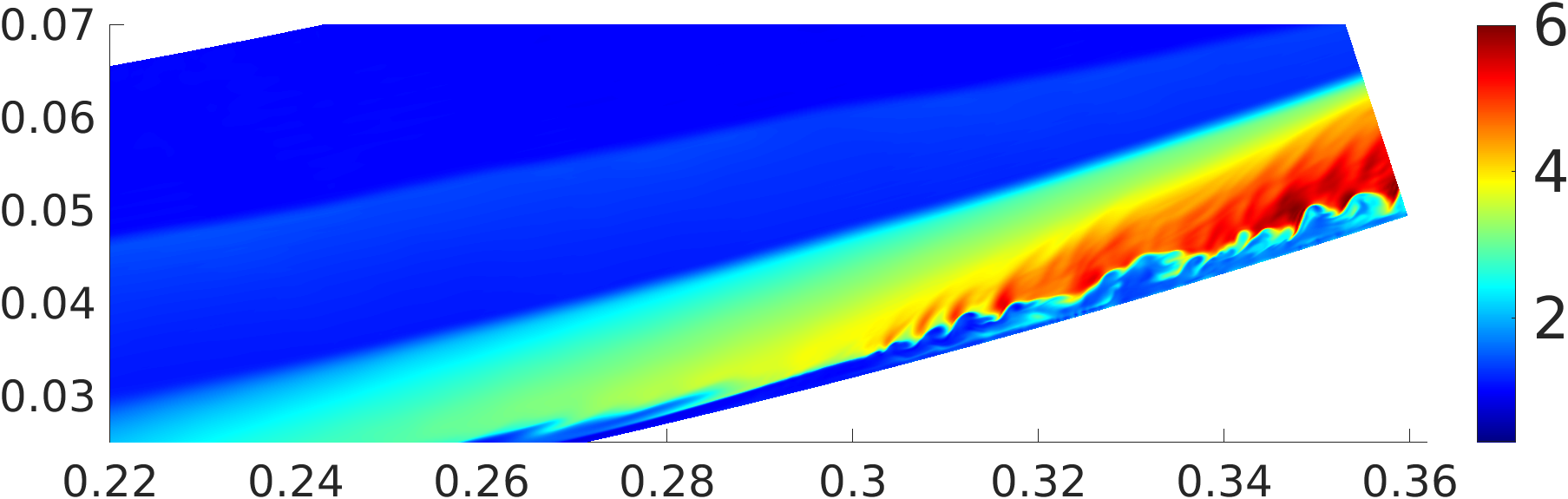}
        \caption{0.5\% TI}
    \end{subfigure}
    \hfill
    \begin{subfigure}[b]{0.48\linewidth}
        \centering
        \includegraphics[width=\linewidth]{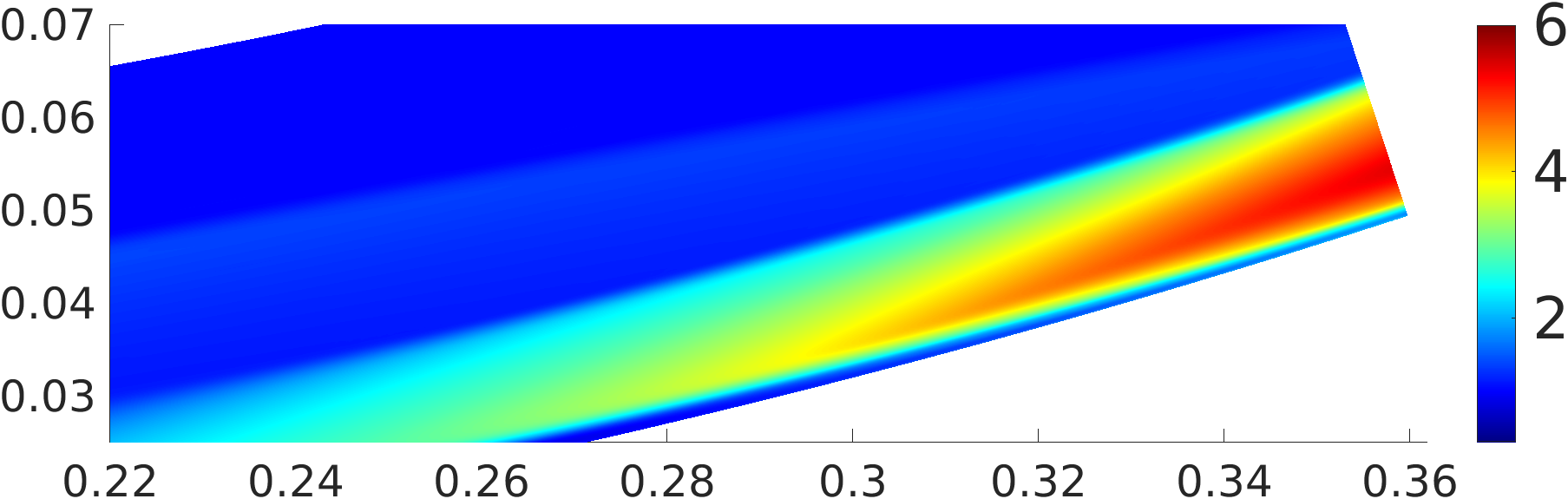}
        \caption{0.5\% TI}
    \end{subfigure}

    \begin{subfigure}[b]{0.48\linewidth}
        \centering
        \includegraphics[width=\linewidth]{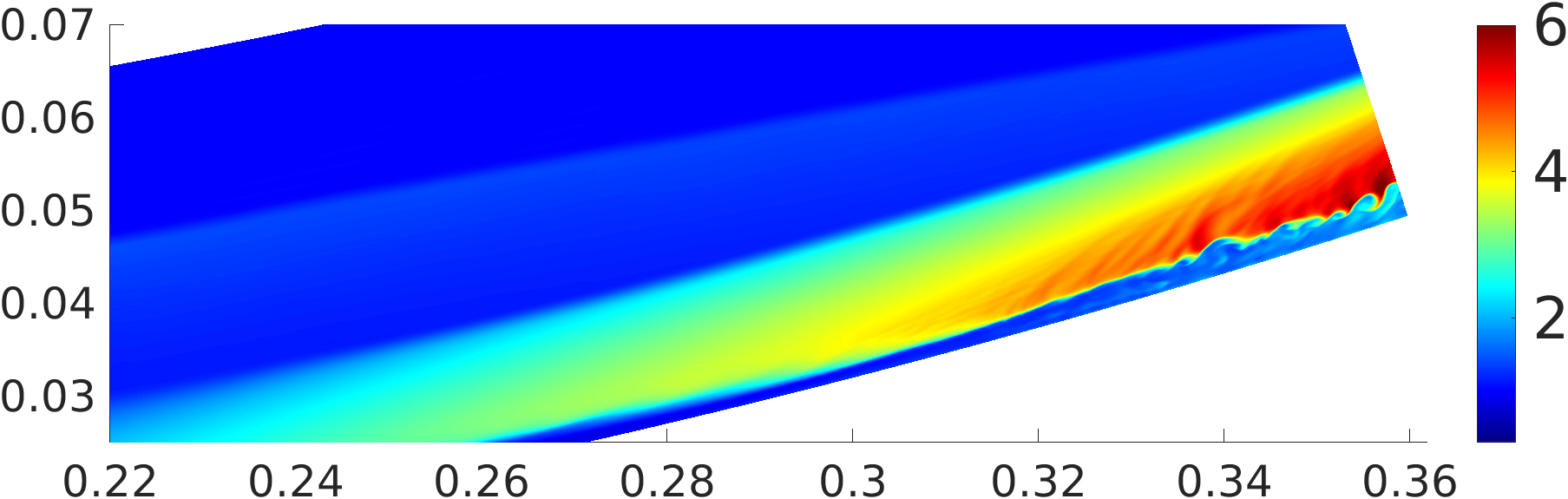}
        \caption{0.2\% TI}
    \end{subfigure}
    \hfill
    \begin{subfigure}[b]{0.48\linewidth}
        \centering
        \includegraphics[width=\linewidth]{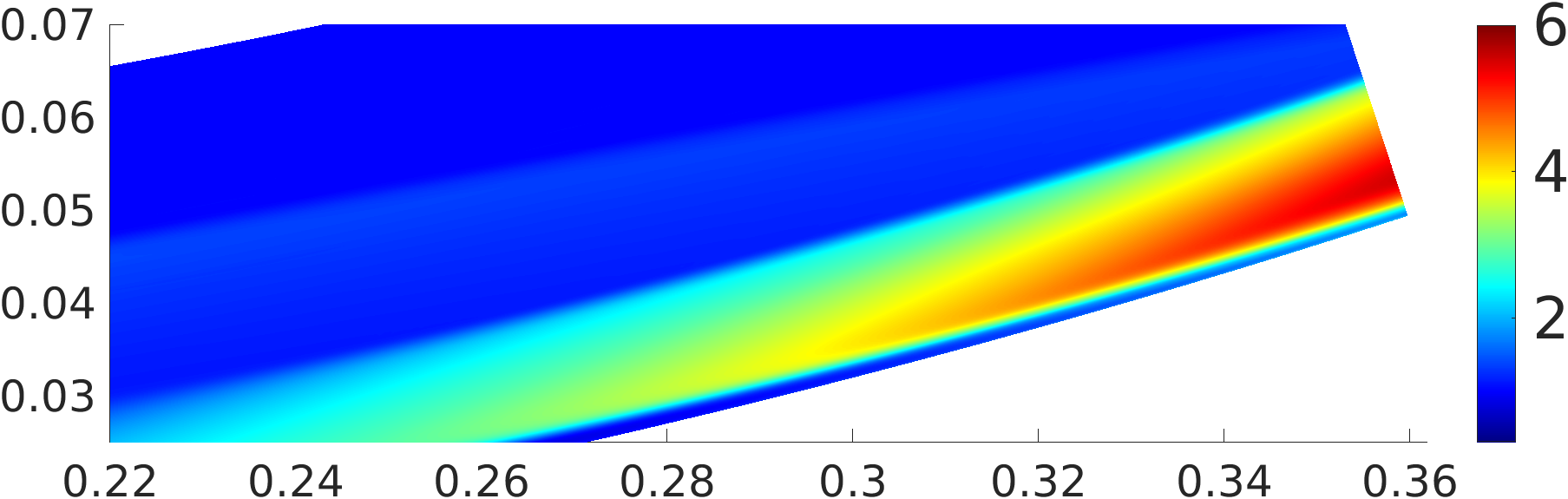}
        \caption{0.2\% TI}
    \end{subfigure}

    \begin{subfigure}[b]{0.48\linewidth}
        \centering
        \includegraphics[width=\linewidth]{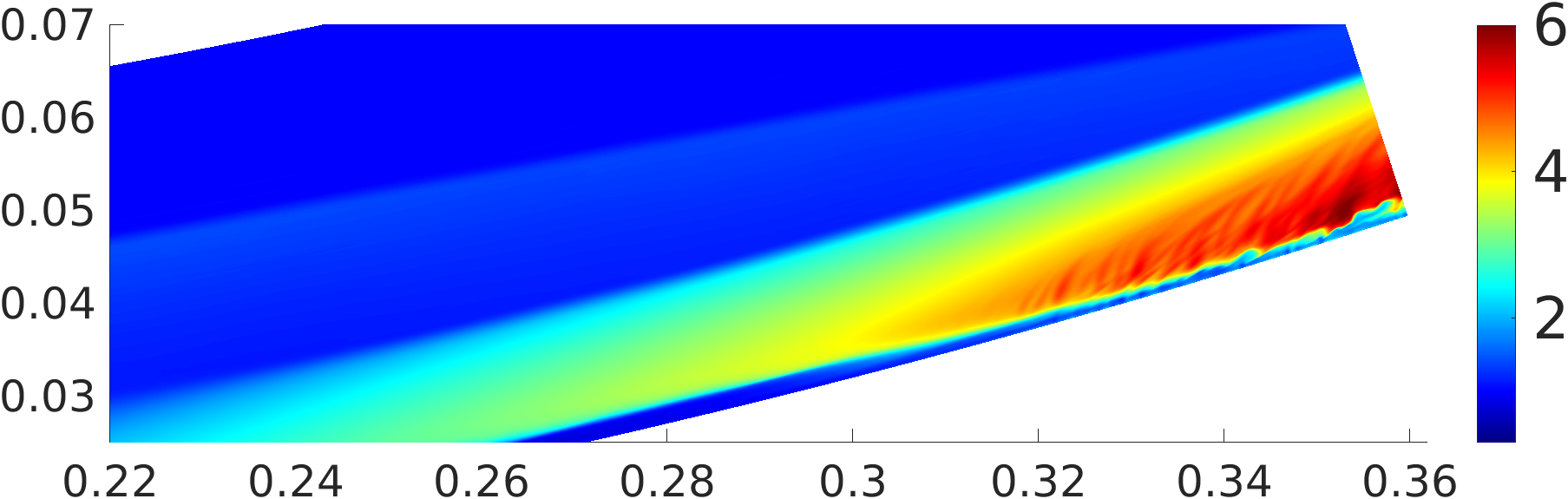}
        \caption{0.1\% TI}
    \end{subfigure}
    \hfill    
    \begin{subfigure}[b]{0.48\linewidth}
        \centering
        \includegraphics[width=\linewidth]{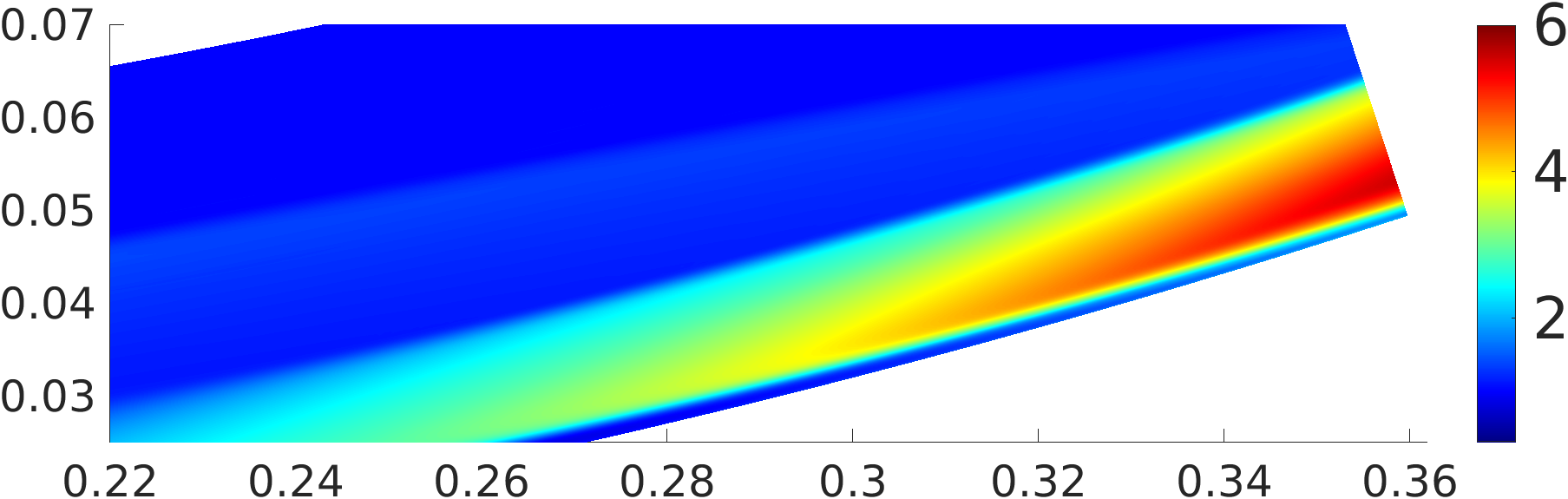}
        \caption{0.1\% TI}
    \end{subfigure}
    \caption{Dimensionless density at the representative side-view slice used
    in Fig.~\ref{pressureside} for freestream turbulence intensities of
    $0.5\%$, $0.2\%$, and $0.1\%$. Instantaneous fields are shown in the left
    column and the corresponding time-averaged fields in the right column.}
    \label{densityside}
\end{figure}

The instantaneous pressure fields in Fig.~\ref{pressureside} illustrate the downstream amplification of the disturbance field. The upstream boundary layer remains comparatively smooth, whereas farther downstream the pressure field develops increasingly pronounced wave-like variations followed by a more irregular near-wall structure associated with nonlinear transition. The corresponding density fields in Fig.~\ref{densityside} exhibit the same progression more clearly: an initially smooth compressed layer develops finite-amplitude distortions that increase in spatial complexity with downstream distance before giving way to a strongly fluctuating transitional region. These instantaneous fields are consistent with the simultaneous presence of compressible waves and deformation of the boundary layer by the vortical streaks. Comparison among the three forcing levels indicates that appreciable small-scale fluctuations appear farther upstream for $0.5\%$ TI and are progressively displaced downstream as the turbulence intensity is reduced. Temporal averaging removes much of the oscillatory and intermittent content, leaving primarily the mean compression and boundary-layer development. 

\begin{figure}[h]
    \centering
    \begin{subfigure}[b]{\textwidth}
        \centering
        \includegraphics[width=\textwidth]{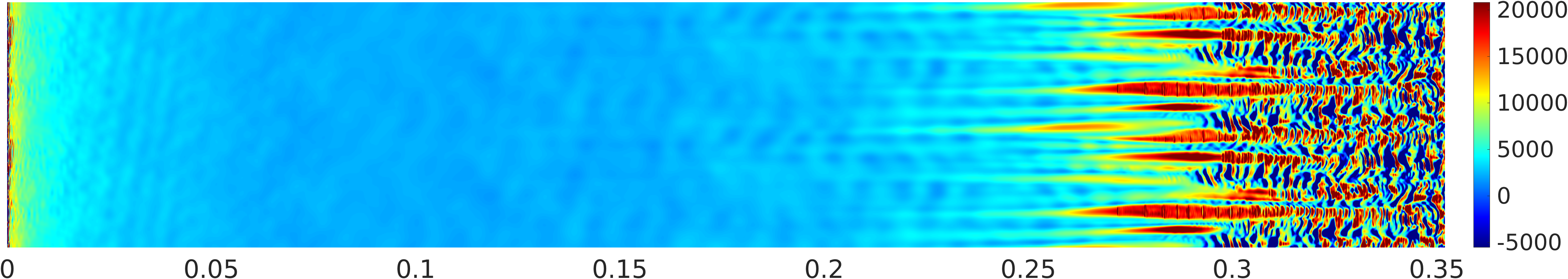}
        \caption{0.5\% TI}
    \end{subfigure}

    \begin{subfigure}[b]{\textwidth}
        \centering
        \includegraphics[width=\textwidth]{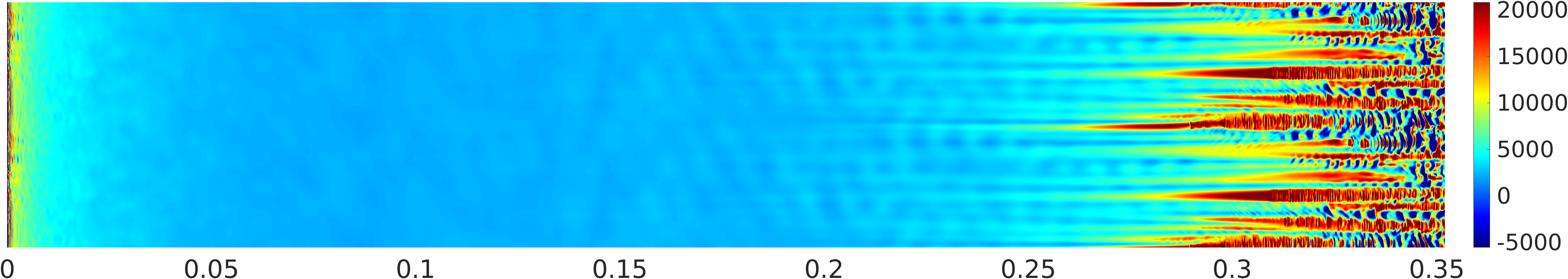}
        \caption{0.2\% TI}
    \end{subfigure}    

    \begin{subfigure}[b]{\textwidth}
        \centering
        \includegraphics[width=\textwidth]{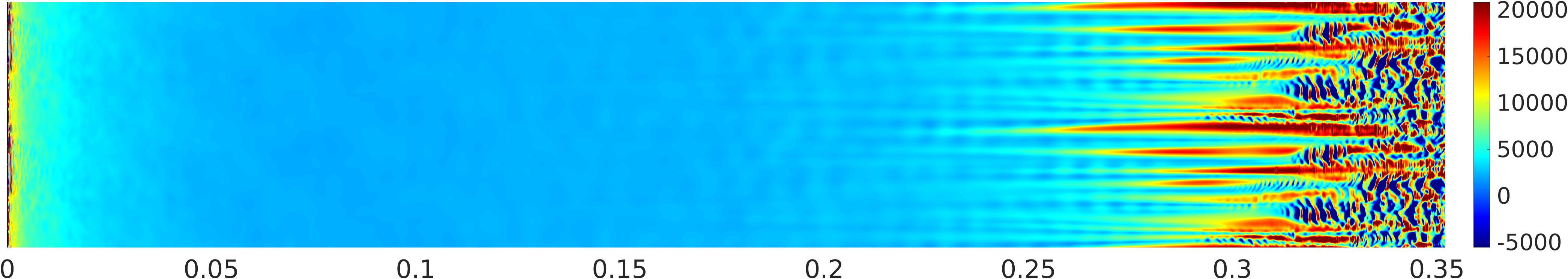}
        \caption{0.1\% TI}
    \end{subfigure}

    \caption{Instantaneous wall heat flux, in $\mathrm{W\,m^{-2}}$, for
    freestream turbulence intensities of $0.5\%$, $0.2\%$, and $0.1\%$ at
    time step 4600. Elongated streamwise bands provide the wall signature of
    the G{\"o}rtler-vortex/streak system, while their downstream loss of
    coherence accompanies nonlinear breakdown.}
    \label{fig:heatfluxTI}
\end{figure}

The wall heat-flux distributions provide a particularly clear signature of the centrifugal-instability dynamics. In the instantaneous fields of Fig.~\ref{fig:heatfluxTI}, weak upstream perturbations are followed by elongated streamwise bands of enhanced and reduced wall heat flux. These structures are consistent with the surface imprint of the G{\"o}rtler-vortex/streak system. Counter-rotating streamwise vortices redistribute momentum and thermal energy across the boundary layer: downwash transports relatively high-momentum, high-enthalpy fluid toward the wall and locally enhances heat transfer, whereas neighboring upwash produces regions of reduced wall heating. With increasing downstream distance, the initially coherent streaks become increasingly modulated in the spanwise direction, lose regularity, and eventually fragment into an irregular wall pattern characteristic of nonlinear breakdown. This sequence occurs for all three turbulence intensities, but at systematically different streamwise locations. The $0.5\%$ case exhibits the earliest loss of streak coherence, whereas the $0.2\%$ and $0.1\%$ cases retain organized streamwise structures farther downstream. The weak mottled disturbances near the upstream portion of the wall are also more pronounced at higher TI, consistent with the larger amplitude of the imposed freestream forcing.

\begin{figure}[h]
    \centering
    \begin{subfigure}[b]{\textwidth}
        \centering
        \includegraphics[width=\textwidth]{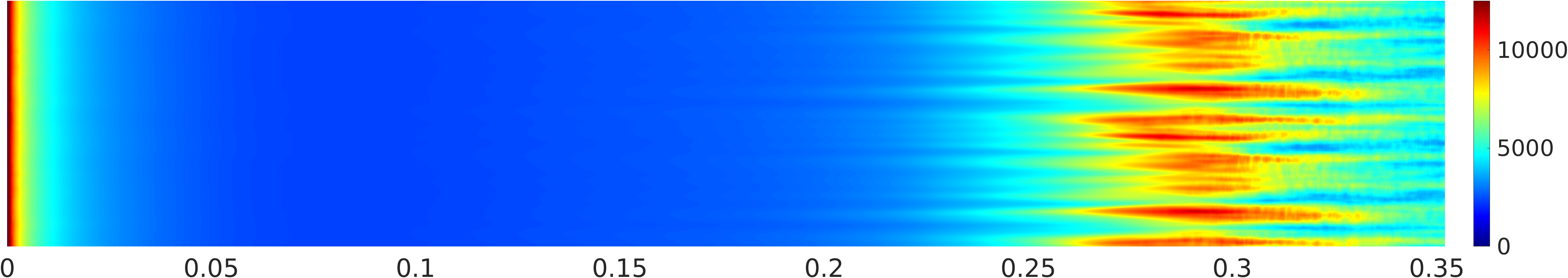}
        \caption{0.5\% TI}
    \end{subfigure}    

    \begin{subfigure}[b]{\textwidth}
        \centering
        \includegraphics[width=\textwidth]{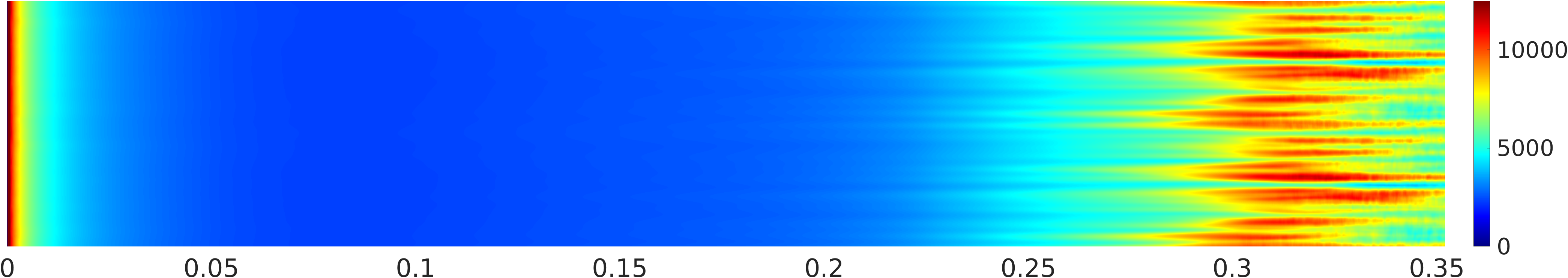}
        \caption{0.2\% TI}
    \end{subfigure}    

    \begin{subfigure}[b]{\textwidth}
        \centering
        \includegraphics[width=\textwidth]{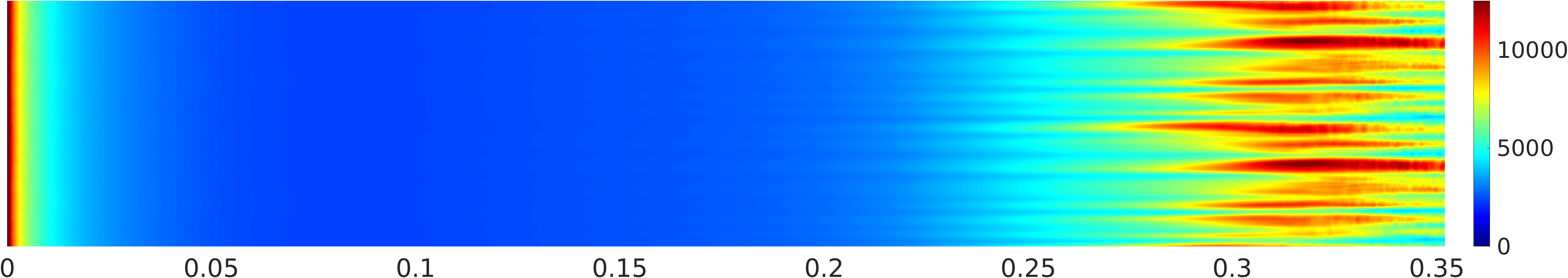}
        \caption{0.1\% TI}
    \end{subfigure}

    \caption{Time-averaged wall heat flux, in $\mathrm{W\,m^{-2}}$, for
    freestream turbulence intensities of $0.5\%$, $0.2\%$, and $0.1\%$.
    Persistent streamwise bands reveal statistically coherent thermal streaks,
    while their downstream displacement with decreasing TI is consistent with
    delayed transition.}
    \label{fig:heatfluxta}
\end{figure}

\begin{figure}[h]
    \centering
    \begin{subfigure}[b]{\textwidth}
        \centering
        \includegraphics[width=\textwidth]{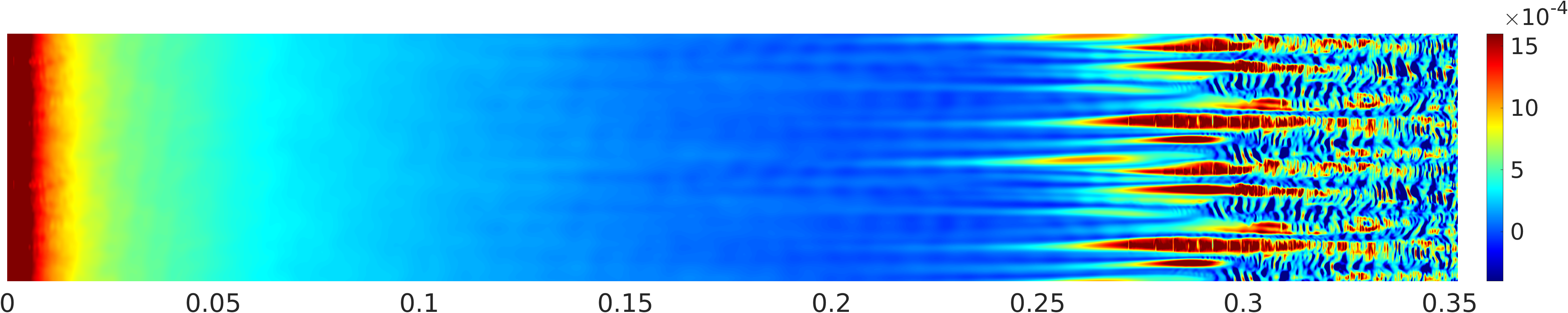}
        \caption{0.5\% TI}
    \end{subfigure}    

    \begin{subfigure}[b]{\textwidth}
        \centering
        \includegraphics[width=\textwidth]{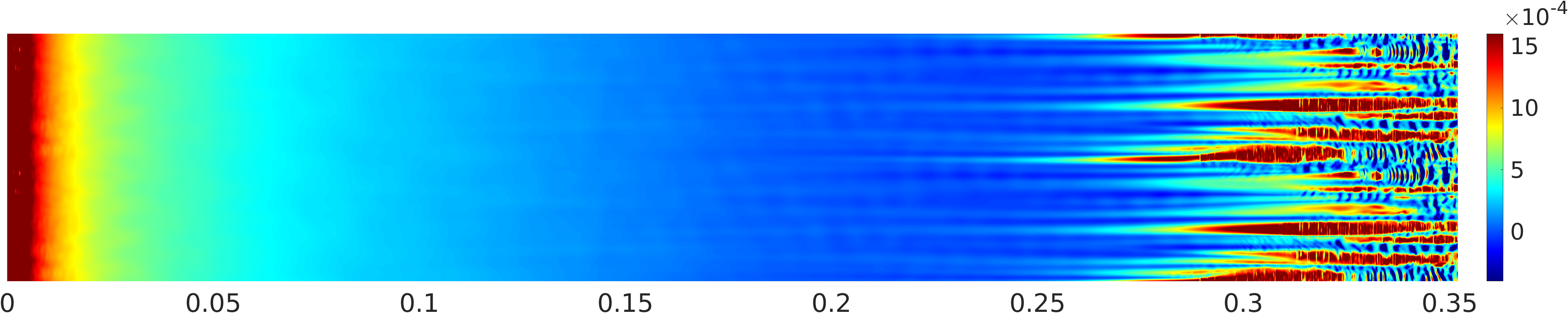}
        \caption{0.2\% TI}
    \end{subfigure}    

    \begin{subfigure}[b]{\textwidth}
        \centering
        \includegraphics[width=\textwidth]{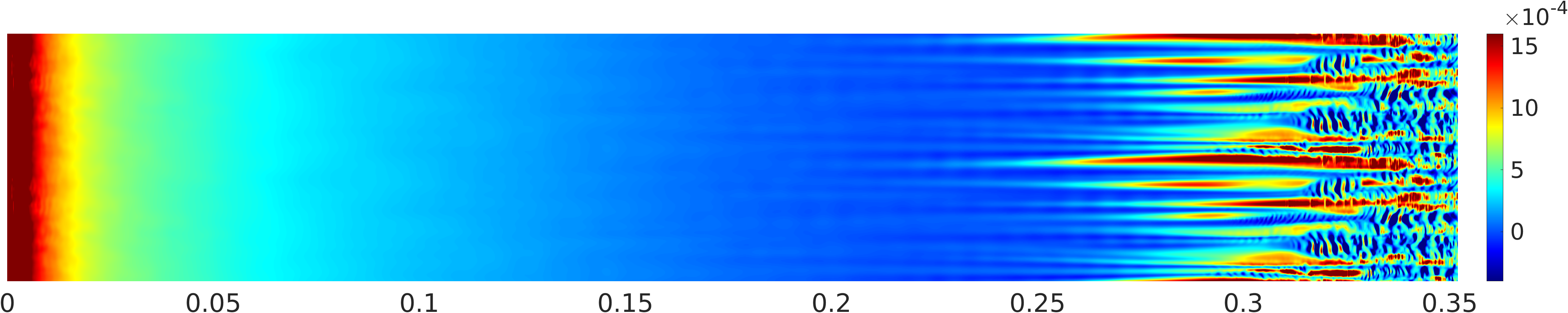}
        \caption{0.1\% TI}
    \end{subfigure}
    \caption{Instantaneous dimensionless skin friction for freestream
    turbulence intensities of $0.5\%$, $0.2\%$, and $0.1\%$ at time step
    4600. The streamwise wall-shear streaks provide an independent signature
    of the near-wall momentum redistribution associated with the
    G{\"o}rtler-vortex system.}
    \label{fig:skinfrictionTI}
\end{figure}

\begin{figure}[h]
    \centering
    \begin{subfigure}[b]{\textwidth}
        \centering
        \includegraphics[width=\textwidth]{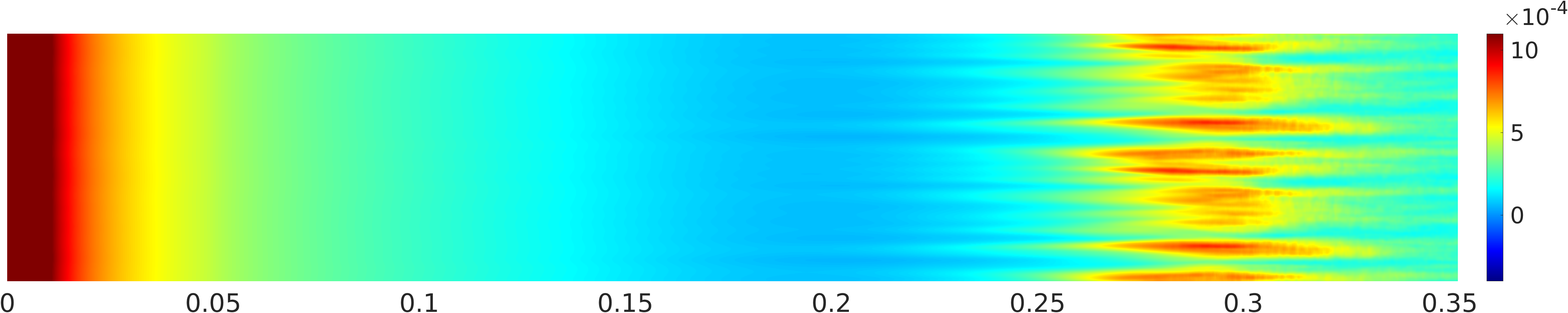}
        \caption{0.5\% TI}
    \end{subfigure}

    \begin{subfigure}[b]{\textwidth}
        \centering
        \includegraphics[width=\textwidth]{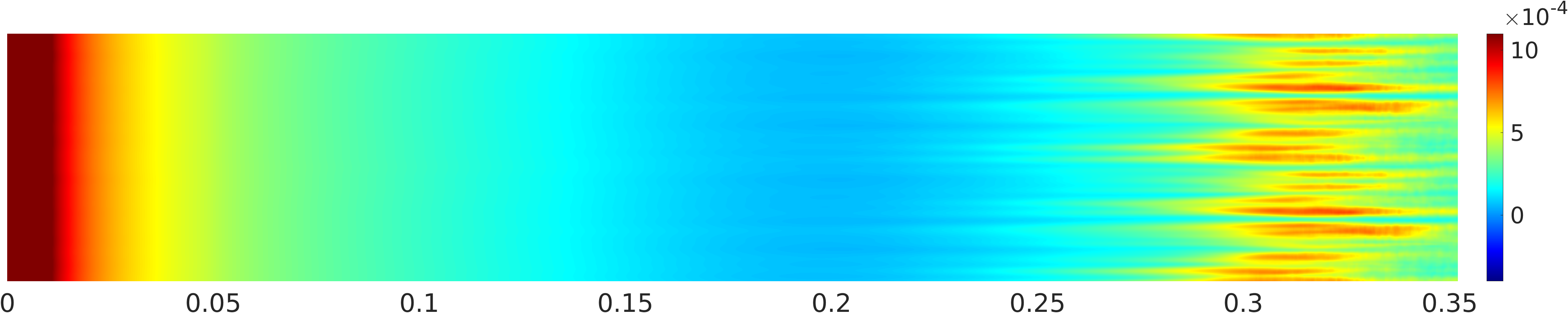}
        \caption{0.2\% TI}
    \end{subfigure}

    \begin{subfigure}[b]{\textwidth}
        \centering
        \includegraphics[width=\textwidth]{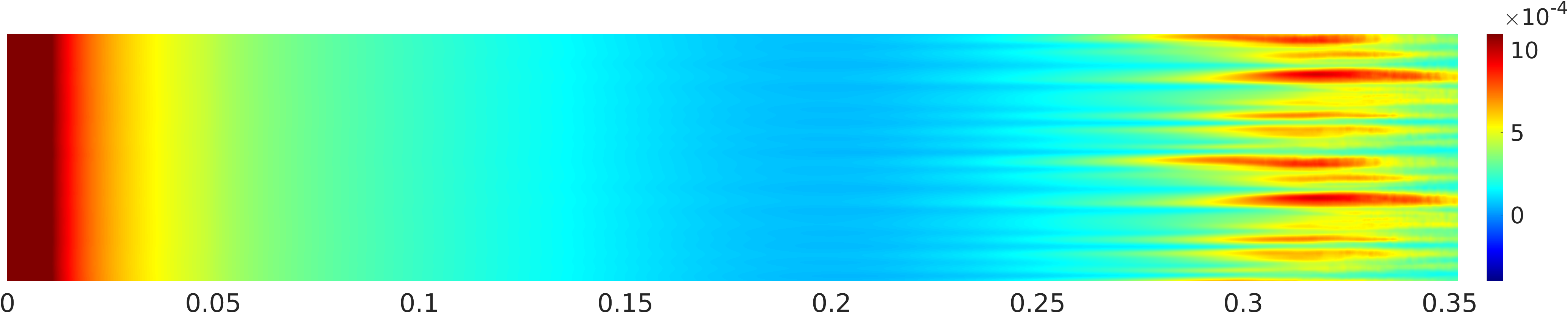}
        \caption{0.1\% TI}
    \end{subfigure}

    \caption{Time-averaged dimensionless skin friction for freestream
    turbulence intensities of $0.5\%$, $0.2\%$, and $0.1\%$. The averaged
    fields retain the persistent streak pattern and show a downstream
    displacement of the transition region as the turbulence intensity is
    reduced.}
    \label{fig:skinfrictionta}
\end{figure}

Temporal averaging distinguishes statistically persistent structures from instantaneous intermittency. The time-averaged heat-flux distributions in Fig.~\ref{fig:heatfluxta} retain streamwise bands over a substantial portion of the surface, demonstrating that the dominant streak pattern is not purely transient but persists statistically over the averaging interval. As the freestream turbulence intensity decreases, the region of strong streak amplification and the subsequent high-heat-flux transitional zone move downstream. In particular, the lower-TI cases preserve coherent thermal streaks over a greater streamwise extent before the wall field becomes spatially disordered. Thus, weaker freestream forcing delays nonlinear breakdown while allowing the curvature-driven streak system to evolve coherently over a longer downstream distance.

The skin-friction distributions in Figs.~\ref{fig:skinfrictionTI} and \ref{fig:skinfrictionta} provide an independent measure of the same near-wall dynamics. The instantaneous skin-friction fields contain elongated streamwise bands and downstream fragmentation that closely parallel the structures observed in the wall heat-flux distributions. This correspondence indicates that the thermal streaks are associated with the underlying redistribution of near-wall momentum rather than with an isolated thermal response. The time-averaged skin-friction fields preserve the dominant streak pattern and again show a downstream displacement of the transitional region as the forcing amplitude is reduced. The consistent behavior of heat flux and skin friction therefore provides complementary evidence that increasing freestream turbulence intensity accelerates the transition process while leaving the overall curvature-driven transition pathway qualitatively unchanged.

\begin{figure}[h]
    \centering
    \begin{subfigure}[b]{0.49\textwidth}
        \centering
        \includegraphics[width=\linewidth]{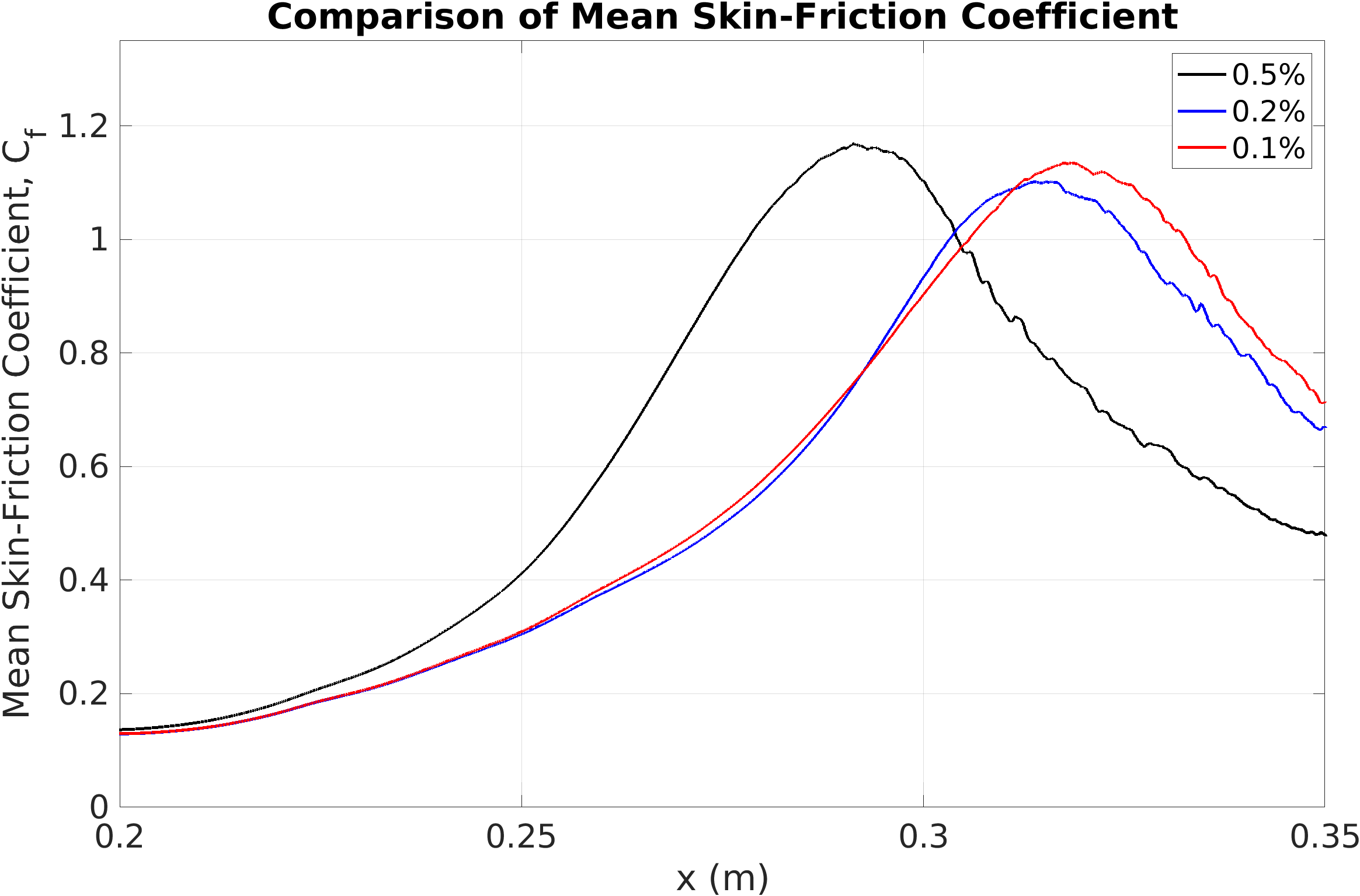}
        \caption{Skin-friction coefficient $C_f$}
        \label{fig:cf}
    \end{subfigure}    
    \hfill
    \begin{subfigure}[b]{0.49\textwidth}
        \centering
        \includegraphics[width=\linewidth]{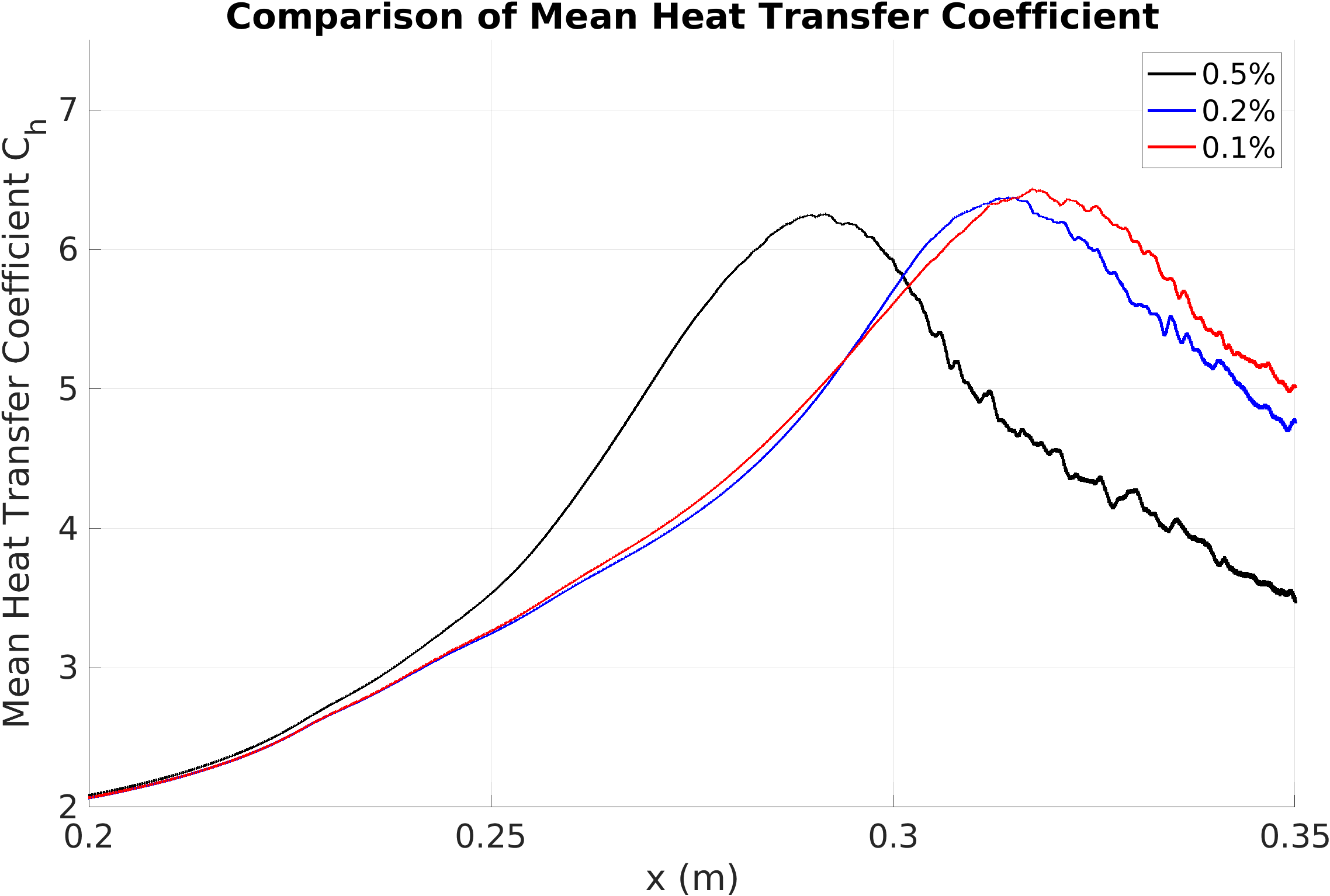}
        \caption{Heat-transfer coefficient $C_h$}
        \label{fig:ch}
    \end{subfigure}

    \caption{Spanwise/time-averaged wall coefficients for freestream
    turbulence intensities of $0.5\%$, $0.2\%$, and $0.1\%$: (a)
    skin-friction coefficient $C_f$ and (b) heat-transfer coefficient $C_h$.
    The downstream displacement of the principal rise and peak with decreasing
    TI is consistent with delayed transition.}
    \label{fig:averaged_coefficients}
\end{figure}

These trends are summarized quantitatively by the spanwise/time-averaged wall coefficients in Fig.~\ref{fig:averaged_coefficients}. Both $C_f$ and $C_h$ vary smoothly in the upstream region before exhibiting a more rapid increase as the streamwise streaks amplify and the boundary layer approaches breakdown. The $0.5\%$ TI case departs from the lower-intensity curves farther upstream and reaches its strongest response earlier, consistent with the upstream shift of transition observed in the instantaneous and time-averaged wall fields. The $0.2\%$ and $0.1\%$ cases remain close over a larger portion of the upstream domain before rising and reaching their maxima farther downstream. In particular, the heat-transfer response of the $0.1\%$ case extends farther downstream than that of the $0.2\%$ case. This behavior is consistent with the wall maps: reducing the forcing amplitude postpones nonlinear breakdown, while the longer coherent development of the G{\"o}rtler streaks can still produce substantial local wall heating and shear before the flow becomes fully disordered. The subsequent reduction of the averaged coefficients downstream of their maxima reflects the evolution beyond the localized transitional peak within the domain.

Together, interior and wall results demonstrate a consistent transition path. The upstream boundary layer remains comparatively smooth apart from the imposed disturbance field. Over the streamwise-concave surface, centrifugal instability amplifies streamwise G{\"o}rtler vortices and produces persistent streaks in both wall heat transfer and skin friction. These streaks increase in amplitude, develop three-dimensional modulation, and eventually lose coherence as nonlinear interactions intensify, leading to the irregular small-scale structure associated with turbulent breakdown. Increasing the freestream turbulence intensity advances this sequence, whereas reducing the forcing delays breakdown and shifts the principal wall-response features downstream. The transition pathway remains predominantly driven by G{\"o}rtler, while compressible wave activity remains evident in instantaneous fields and contributes to the broader disturbance environment. The present simulations capture weak upstream disturbances, organized streamwise vortices, strongly amplified wall streaks, nonlinear modulation, and the emergence of progressively smaller-scale structures during transition.

\section{Conclusions}

This paper developed a massively parallel HDG framework for direct numerical simulation of compressible transitional and turbulent flows on GPUs. The formulation combines a high-order HDG discretization of the compressible Navier--Stokes equations with robust shock capturing, implicit DIRK time integration, Newton--GMRES method, additive Schwarz preconditioning, and reduced-basis method. The significance of this combination is that accuracy, robustness, and scalability are addressed together rather than as separate algorithmic components. High-order HDG provides geometric flexibility and low-dissipation resolution of transitional and turbulent structures; localized artificial viscosity regularizes shocks and strong compressions; implicit time integration and Newton--Krylov solution permit stable advancement on highly refined wall-bounded meshes; and the GPU implementation targets the dense local algebra and sparse trace communication patterns generated by hybridization. The benchmark problems span subsonic, supersonic, and hypersonic regimes and test the solver over a range of flow physics.

For the Eppler 387 airfoil, the simulations  capture subsonic laminar--turbulent transition over a curved aerodynamic configuration with laminar separation, shear-layer instability, reattachment, and wake development. The computed pressure distributions agree closely with the experimental measurements of McGhee et al. over the four Reynolds numbers considered, including the suction peak, pressure recovery, and Reynolds-number dependence of the suction-side pressure plateau. The boundary-layer integral quantities, skin-friction distributions, pressure and velocity fluctuation fields, turbulent kinetic energy, and instantaneous vortical structures provide a consistent description of the transition process: increasing Reynolds number moves the onset of three-dimensionality upstream, shortens the separated region, and produces a more compact turbulent recovery.

The supersonic Taylor--Green vortex involves strong compressibility, shocklet-like structures, and transition to decaying three-dimensional turbulence. The simulations show the deformation of the initially organized vortex system, the formation of smaller-scale vortical structures, and the coexistence of vortical and compressive dynamics. The mesh-sequence study shows convergence of the kinetic-energy and dissipation histories toward the reference DNS data. Kinetic energy is less resolution-sensitive than the dissipation measures, whereas the total, solenoidal, and especially dilatational dissipation require finer resolution to reproduce the timing and magnitude of the benchmark behavior. 

The Mach-7 transition simulations involve strongly compressible wall-bounded flow with thin thermal and viscous layers, curved geometry, imposed freestream disturbances, G{\"o}rtler-vortex amplification, acoustic/compressive activity, and nonlinear breakdown. The pressure and density fields show the downstream growth of instantaneous disturbances and the emergence of increasingly complex compressible structures. The wall heat-flux and skin-friction distributions provide the clearest signatures of transition: elongated streamwise streaks form over the concave surface, amplify downstream, undergo three-dimensional modulation, and eventually lose coherence as the flow breaks down. Reducing the freestream turbulence intensity from $0.5\%$ to $0.2\%$ and $0.1\%$ shifts the dominant wall-response features and the breakdown process downstream, while preserving the same qualitative G{\"o}rtler-dominated transition pathway.

This work uses carefully designed high-order meshes that concentrate resolution in boundary layers, separated shear layers, vortical structures, and regions where compressive features are expected. A natural next step is to make this resolution placement dynamic through adaptive high-order discretizations. Future work should combine $h$-, $p$-, and $hp$-adaptivity with anisotropic mesh modification and shock- or feature-aware error indicators. The HDG formulation is well suited to such developments because many operations remain element-local after hybridization, and the trace formulation provides a compact global coupling on adapted meshes. The principal challenges are to perform mesh refinement, coarsening, solution transfer, and trace-space repartitioning without degrading nonlinear robustness, load balance, or GPU efficiency. 

The hypersonic flow simulations reported here use the calorically perfect-gas compressible Navier--Stokes equations. Predictive simulation of high-enthalpy hypersonic flight and atmospheric-entry flows requires extending the framework to real-gas thermochemistry, including temperature-dependent thermodynamic and transport properties, multispecies reacting mixtures, finite-rate chemistry, vibrational and thermal nonequilibrium, dissociation, and, where relevant, ionization and nonideal equations of state. 
Real-gas extensions will also introduce substantial numerical challenges. Additional species and energy variables increase the state dimension and strengthen nonlinear coupling; finite-rate chemistry and nonequilibrium relaxation introduce stiffness; positivity and realizability constraints become more restrictive; and shock stabilization must remain compatible with both thermochemical source terms and resolved instability dynamics. These features make scalable implicit solvers and effective preconditioners even more important than in the perfect-gas setting considered in this paper. 

Beyond these two directions, further gains may come from adaptive time stepping, stronger physics-based or data-assisted preconditioners, more aggressive reduced-order acceleration, and larger validation studies on realistic aerodynamic configurations. The broader objective is to make high-fidelity compressible DNS increasingly practical on heterogeneous supercomputers while retaining the resolution, robustness, and physical fidelity required for high-speed turbulent flows.

\section*{Acknowledgements} \label{}

This material is based upon work supported by the U.S. Department of Energy, National Nuclear Security Administration, under Award Number DE-NA0004266. The authors gratefully acknowledge the computational resources provided by the Frontier supercomputer at the Oak Ridge Leadership Computing Facility and the Tuolumne supercomputer at Lawrence Livermore National Laboratory. These high-performance computing resources were instrumental in enabling the numerical simulations presented in this work.

\

\noindent
\textbf{Disclaimer:} This report was prepared as an account of work sponsored by an agency of the United States Government. Neither the United States Government nor any agency thereof, nor any of their employees, makes any warranty, express or implied, or assumes any legal liability or responsibility for the accuracy, completeness, or usefulness of any information, apparatus, product, or process disclosed, or represents that its use would not infringe privately owned rights. Reference herein to any specific commercial product, process, or service by trade name, trademark, manufacturer, or otherwise does not necessarily constitute or imply its endorsement, recommendation, or favoring by the United States Government or any agency thereof. The views and opinions of authors expressed herein do not necessarily state or reflect those of the United States Government or any agency thereof.

\bibliographystyle{elsarticle-num}
\bibliography{library,library2}
\end{document}